\documentclass[12pt]{amsart}
\usepackage{ amsmath, amsthm, amsfonts, amssymb, color}
 \usepackage{mathrsfs}
\usepackage{amsfonts, amsmath}
 \usepackage{amsmath,amstext,amsthm,amssymb,amsxtra}
 \usepackage{txfonts} 
 \usepackage[colorlinks, citecolor=blue,pagebackref,hypertexnames=false]{hyperref}
 \allowdisplaybreaks
 \usepackage{pgf,tikz}
\begin{document}

 \baselineskip 16.6pt
\hfuzz=6pt

\widowpenalty=10000

\newtheorem{cl}{Claim}
\newtheorem{thm}{Theorem}[section]
\newtheorem{prop}[thm]{Proposition}
\newtheorem{cor}[thm]{Corollary}
\newtheorem{lemma}[thm]{Lemma}
\newtheorem{definition}[thm]{Definition}
\newtheorem{assum}{Assumption}[section]
\newtheorem{example}[thm]{Example}
\newtheorem{remark}[thm]{Remark}
\renewcommand{\theequation}
{\thesection.\arabic{equation}}

\def\SL{\sqrt H}

\newcommand{\mar}[1]{{\marginpar{\sffamily{\scriptsize
        #1}}}}

\newcommand{\as}[1]{{\mar{AS:#1}}}

\newtheorem*{fact}{Fact}

\theoremstyle{assumption}
\newtheorem{assumption}[thm]{Assumption}

\theoremstyle{definition}

\newcommand{\cA}{\mathcal{A}}
\newcommand{\cB}{\mathcal{B}}
\newcommand{\cC}{\mathcal{C}}
\newcommand{\cF}{\mathcal{F}}
\newcommand{\cH}{\mathcal{H}}
\newcommand{\cK}{\mathcal{K}}
\newcommand{\cM}{\mathcal{M}}
\newcommand{\cN}{\mathcal{N}}
\newcommand{\cP}{\mathcal{P}}
\newcommand{\cS}{{\mathcal{S}}}
\newcommand{\cU}{\mathcal{U}}
\newcommand{\cX}{\mathcal{X}}
\newcommand{\dsone}{\mathds{1}}

\newcommand{\bB}{\mathbb{B}}
\newcommand{\bC}{\mathbb{C}}
\newcommand{\bD}{\mathbb{D}}
\newcommand{\bE}{\mathbb{E}}
\newcommand{\bF}{\mathbb{F}}
\newcommand{\bG}{\mathbb{G}}
\newcommand{\bN}{\mathbb{N}}
\newcommand{\bR}{\mathbb{R}}
\newcommand{\bZ}{\mathbb{Z}}
\newcommand{\bT}{\mathbb{T}}
\newcommand{\di}{\mathrm{div}}

\newcommand{\ctimes}{\rtimes_\theta}
\newcommand{\mtx}[4]{\left(\begin{array}{cc}#1&#2\\#3&#4\end{array} \right)}
\newcommand{\dvarphi}[2]{D^{(#1)}(\varphi_{#2})}
\newcommand{\xpsi}[2]{\psi^{[#1]}_{#2}(\omega)}
\newcommand{\xPsi}[2]{\Psi^{[#1]}_{#2}}
\newcommand{\dt}[1]{\partial_t^{(#1)}}
\newcommand{\vN}{vN(G)}
\newcommand{\add}[1]{\quad \text{ #1 } \quad}
\newcommand{\xp}{\rtimes_\theta }
\newcommand{\crossa}{\cM \rtimes_\theta G}
\newcommand{\uu}[2]{u_{#1}^{(#2)}}

\newcommand{\Lp}{{L_p(vN(G))}}
\newcommand{\Lplcr}[1]{{L_p(#1;\ell_2^{cr})}}
\newcommand{\Lplc}[1]{{L_p(#1;\ell_2^{c})}}
\newcommand{\Lplr}[1]{{L_p(#1;\ell_2^{r})}}
\newcommand{\Lpinfty}[1]{{L_p(#1;\ell_\infty)}}
\newcommand{\Lpli}[1]{{L_p(#1;\ell_1)}}
\newcommand{\Lpco}[1]{{L_p(#1;c_0)}}
\newcommand{\Lpcor}[1]{{L_p(#1;c_0)}}
\newcommand{\LpR}{{L_p(\bR^d;L_p(\cN))}}

\newcommand\CC{\mathbb{C}}
\newcommand\NN{\mathbb{N}}
\newcommand\ZZ{\mathbb{Z}}
\newcommand\hDelta{{\bf L}}
\renewcommand\Re{\operatorname{Re}}
\renewcommand\Im{\operatorname{Im}}

\newcommand{\mc}{\mathcal}
\newcommand\D{\mathcal{D}}
\def\hs{\hspace{0.33cm}}
\newcommand{\la}{\alpha}
\def \l {\alpha}
\def\ls{\lesssim}
\def\su{{\sum_{i\in\nn}}}
\def\lz{\lambda}
\newcommand{\eps}{\varepsilon}
\newcommand{\pl}{\partial}
\newcommand{\supp}{{\rm supp}{\hspace{.05cm}}}
\newcommand{\x}{\times}
\newcommand{\lag}{\langle}
\newcommand{\rag}{\rangle}

\newcommand\wrt{\,{\rm d}}
\newcommand{\botimes}{\bar{\otimes}}
\def\nn{{\mathbb N}}
\def\bx{{\mathbb X}}
\def\fz{\infty}
\def\r{\right}
\def\lf{\left}
\def\cm{{\mathcal M}}
\def\cs{{\mathcal S}}
\def\rr{{\mathbb R}}
\def\lm{\mathcal{N}}
\def\rd{\mathbb{R}^d_{\theta}}
\def\ri{\mathbb{R}^2_{\theta}}
\def\rn{{{\rr}^n}}
\def\zz{{\mathbb Z}}
\def\cl{{\mathcal L}}
\def\cq{{\mathcal Q}}
\def\cd{{\mathcal D}}

\title[Navier-Stokes equation on quantum Euclidean spaces]{Navier-Stokes equations on quantum Euclidean spaces}

\author{Deyu Chen}
\address{Deyu Chen, Institute for Advanced Study in Mathematics\\ Harbin Institute of Technology\\ Harbin 150001\\China}
\email{1201200317@stu.hit.edu.cn}

\author{Guixiang Hong}
\address{Guixiang Hong, Institute for Advanced Study in Mathematics\\ Harbin Institute of Technology\\ Harbin 150001\\ China}
\email{gxhong@hit.edu.cn}

\author{Liang Wang}
\address{Liang Wang, Department of Mathematics\\ City University of Hong Kong\\Hong Kong SAR\\China
}
\email{wlmath@whu.edu.cn}

\author{Wenhua Wang}
\address{Wenhua Wang, Institute for Advanced Study in Mathematics\\Harbin Institute of Technology\\ Harbin 150001\\China}
\email{whwangmath@whu.edu.cn}

  \date{\today}

\begin{abstract}
We investigate in the present paper the Navier-Stokes equations on quantum Euclidean spaces  $\mathbb{R}^d_{\theta}$ with $\theta$ being a $d\times d$ antisymmetric matrix, which is a standard example of non-compact noncommutative manifolds. The quantum analogues of Ladyzhenskaya and Kato's results are established, that is, we obtain the global well-posedness in the 2D case  and the local well-posedness with solution in $L_d(\mathbb{R}^d)$ in higher dimensions. To achieve these optimal results, we develop the related theory of harmonic analysis and function spaces on $\mathbb{R}^d_{\theta}$, and apply the sharp estimates around noncommutative $L_p$-spaces to quantum Navier-Stokes equations. Moreover, our techniques, which are independent of the deformed parameter $\theta$, allow us to conclude some results on the semiclassical limits. This is the first instance of systematical applications to the theory of quantum partial differential equations
of the powerful real analysis techniques around noncommutative $L_p$-spaces, which date back to the seminal work  \cite{PiXu97} in 1997 on noncommutative martingale inequalities. As in classical case, one may expect numerous similar applications in the future.
\end{abstract}

\maketitle

  \tableofcontents


\section{Introduction}

For a $d\times d $ antisymmetric matrix $\theta$, the quantum Euclidean space  $\rd$ is defined as a von Neumann algebra  generated by a $d$-parameter strongly
continuous unitary family $\{\lambda_{\theta}(t)\}_{t\in\rr^d}$ satisfying the Weyl relation:
\begin{align*}
 \lambda_\theta(t)\lambda_\theta(s)=e^{{\mathrm{i}}(t,\theta s)} \lambda_\theta(s)\lambda_\theta(t),\ \ \ \ \text{for\: any}\: t,s\in\mathbb{R}^d.
 \end{align*}
 There exists a canonical trace $\tau_\theta$ on it, with which one may construct noncommutative $L_p$-spaces. For $1\leq k\leq d$, there is an intrinsic definition of the partial derivative $\partial_k$ on this object with the help of a Fourier-like expansion formula, so many notions such as partial differential equations (abbreviated as PDEs), Fourier multipliers, function spaces  can be naturally defined.
When $\theta=0$, $\rr^d_{\theta}$ will be reduced to $L_{\fz}(\rr^d)$, the quantization of classical Euclidean space $\mathbb{R}^d$. When $d=2n$,
$\theta={\hbar}\left(\begin{array}{ccc}0 & -I_{n} \\ I_{n} & 0 \end{array}\right)$ with the Planck constant $\hbar$, $\rr^d_{\theta}$
is known as the Moyal plane or the phase space.  For more information about quantum Euclidean spaces, we refer the reader to Section \ref{s2} below.

In the literature,
the quantum Euclidean spaces $\rd$ is also regarded as the operator formulation, also called path integral quantization in \cite{BSS2014}, of classical Euclidean spaces but equipped with a noncommutative coordinate system $\{x_j\}_{j=1}^d$ such that $[x_i,x_j]=\mathrm i\theta_{i,j},$ in which case the classical product is replaced by the Moyal product. From this point, Seiberg and Witten \cite{sw99} developed the noncommutative gauge theory, which is an important research object in quantum mechanics (see also \cite{BS1991,dn01,ns98}). As a mathematical object, they are standard model examples of the non-compact manifolds  \cite{cgrs14,ggisv04} in Connes' noncommutative geometry theory. On the other hand, motivated by the noncommutative martingale theory and harmonic analysis \cite{Jun02,JuXu07,Mei07,NeRi11,PiXu97},  there have appeared several fundamental works \cite{FHW25,GJM2022,gjp21,hlw23,hlw25,jmpx21,LSZ20,MSX20} on harmonic analysis over quantum Euclidean spaces since the seminal one \cite{CXY13} in 2013 on quantum tori (see also \cite{XXY18}).

Inspired by the study of condensed matter physics, nuclear physics, and quantum chemistry etc., several PDEs on the phase spaces such as the Bogoliubov-de Gennes equations, the generalized Hartree-Fock equations, and the Hartree-Fock-Bogoliubov equations etc. have naturally appeared  in the literature as the mean-field limit of many body Schr\"odinger equation for wave functions. As PDEs on the Moyal plane, they should be viewed as variants of the von Neumann equation or the Liouville-von Neumann equation, which describe the Schr\"odinger evolution of states; in the Heisenberg picture, the Heisenberg equation describes the Schr\"odinger evolution of observables, which  can be also regarded as one PDE on the Moyal plane. In the last two decades, there have appeared many papers studying the posedness problem of these equations in order to understand the mean-filed limits or the semiclassical limits.  We refer the reader to \cite{NSS18,CLS24,CLS,ESY10} and the references therein for more information on the development of these PDEs.

For a general $d\times d$ antisymmetric matrix $\theta$ there are also many studies on the theory of PDEs. Based on the noncommutative gauge theory, Bahcall and Susskind \cite{BS1991} described the discrete nature of quantum Hall liquid and developed the noncommutative hydrodynamics on quantum Euclidean spaces. Up to now, many authors considered the fluid equations on this noncommutative spaces, such as the Maxwell equations \cite{GJPP2001}, Hamiltonian dynamics \cite{BSS2014,Ma2018}, continuity equations \cite{JPP,Poly} and Euler equations \cite{DG2016}. Besides these fluid equations, Hamanaka and Toda \cite{ht02,ht03,ht06} also derived from the noncommutative Yang-Mills theory and Lax representation the quantization of many integrable systems, such as Burgers equations, KdV equations, mKdV equations and nonlinear Schr{\"o}dinger equations. Although many authors studied these noncommutative equations from different perspectives, there are nearly no result on the theory of well-posedness of these PDEs  due to the bad behaviors of $\theta$-Moyal product on function spaces. These strongly inspire us to study the theory of PDEs on quantum Euclidean spaces.

Besides, with different motivations, several other PDEs have also been investigated in the noncommutative setting.  For instance, Chakraborty, Goswami and Sinha \cite{cgs03} have analyzed diffusion equations on quantum tori;  Rosenberg \cite{r08} developed a theory of nonlinear elliptic PDEs and studied the Laplace equations and its variants over quantum tori; Labuschagne and  Majewski \cite{lm22} studied the quantum Fokker-Planck equations on general von Neumann algebras. It is worthy to note that Carlen and Maas \cite{CM14} as well as Voiculescu \cite{v20} studied
 the hydrodynamic equations on some certain noncommutative algebras from the perspective of noncommutative Wasserstein manifolds,
 which are also important motivations for the study of noncommutative fluid equations.

However, as far as the authors know, almost all the previous investigations on the above mentioned noncommutative PDEs are essentially restricted to $L_p$ for $p=1,2,\infty$ by exploiting the Hilbert space structure $L_2$ or the algebraic structure $L_\infty$. The newly emerging but powerful real analysis over other noncommutative $L_p$-spaces, in particular noncommutative harmonic analysis, seems to have been ignored by the noncommutative PDE community. On the other hand, as some direct applications of noncommutative harmonic analysis, Gonz\'{a}les-P\'{e}rez,
Junge and Parcet \cite{gjp21} studied the $L_p$-regularity of linear elliptic pseudodifferential
equations; Fan, Hong and Wang \cite{FHW25} obtained the sharp endpoint $L_p$ estimates of the free quantum Schr\"odinger equations, and Hong, Lai and Wang \cite{hlw25} also established a local smoothing estimate of the free wave equations on 2D quantum Euclidean space. Motivated by all this, we have a project which aims at exploring real analysis over noncommutative $L_p$-spaces to study noncommutative PDEs, and then finding applications to the theory of mean filed limit and semiclassical limit, and thus facilitating the understanding of condensed matter physics and quantum chemistry etc..

In this paper, we are restricted to considering one of the most fundamental nonlinear PDEs--Navier-Stokes equation on quantum Euclidean spaces. This quantum PDE can be easily derived from \cite{BSS2014} via Weyl quantization. Indeed, it is well-known that the classical incompressible Navier-Stokes equations contains two parts: the momentum conservation equation
 \begin{align*}
   \rho\partial_t \phi-\nu\Delta \phi+\phi\cdot \nabla \phi+\nabla q=0
\end{align*} and the divergence free condition $\mathbf{div} \phi=0,$  where $\rho,\phi,\nu>0$ and $q$ denote the density, the velocity, the viscosity constant and the pressure of the fluid system, respectively. In terms of the Nambu dynamics, these two parts can be rewritten as (cf. \cite{Nambu}):
\begin{align*}
&\rho(\partial_t\{x_i,  \varphi_1,\dots, \varphi_{d-1}\}_N+\{\{x_i, \varphi_1, \dots, \varphi_{d-1}\}_N, \varphi_1, \dots, \varphi_{d-1}\}_N)\\
&+\frac{1}{(d-1)!} \sum_{1\le i_1,\dots,i_d\le d}\epsilon^{i_1,\dots,i_d} \{q, x_{i_1}, \dots,x_{i_d}\}_N - \nu \Delta \{x_i, \varphi_1\dots, \varphi_{d-1}\}_N=0,\\
&\;\mathrm{and}\;\phi_i=\{x_i,\varphi_1, \dots,\varphi_{d-1} \}_N, \ \ \ i=1,\dots,d
\end{align*}
for some stream functions $\varphi_1,\dots,\varphi_{d-1},$ where $\epsilon^{i_1, i_2, \dots, i_d} $ is the Levi-Civita tensor and $\{A_1,\dots,A_d\}_N$ is the Nambu bracket given by $$\{A_1,\dots,A_d\}_N(x,t)=\sum_{1\le i_1, \dots, i_d\le d} \epsilon^{i_1, \dots, i_d} \partial_{i_1} A_1(x,t)\cdots \partial_{i_d} A_d(x,t).$$
The main contribution of Saitou et al. \cite{BSS2014} is that they introduced the Moyal-Nambu bracket in place of the traditional Nambu bracket, where the classical product is replaced by the Moyal product, and exploited the noncommutative Moyal-Nambu dynamics to describe the hydrodynamics of granular materials, resulting in  a striking application to the noncommutative hydrodynamics. Like the classical process, they also deduced the noncommutative Navier-Stokes equation
\begin{align}\label{nsq}
   \lf\{\begin{array}{ll}
\rho\partial_t\phi-\nu\Delta \phi+\phi\star_\theta \nabla \phi+\nabla q=0;\\[5pt]
\mathbf{div} \phi=0;\\[5pt]
\phi(0)=\phi_0,
\end{array}\r.
\end{align}
where $\star_\theta$ is the Moyal product defined by
$$(f\star_\theta g)(x):=\exp(\frac {\mathrm{i}}2\sum_{1\le i,j\le d}\theta_{i,j}\partial_{y_i}\partial_{z_j})f(y)g(z)|_{y=z=x},
$$ $\phi\star_\theta \nabla\phi:=\lf\{\sum_{j=1}^d\phi_j\star_\theta(\partial_j\phi_k)\r\}_{k=1}^d$.
An equivalent definition of $\star_\theta$ is given by Rieffel \cite{Rieffel}, which will be introduced in Section \ref{s8}.
For brevity, we often set $\rho=1$ and $\nu=1$ (when $\nu=0,$ it will become the noncommutative Euler equation).
Let $U_\theta$ be the Weyl transform and  $\mathcal{F}$ be the Fourier transform that can be found in Section \ref{s2}.
Via the Weyl quantization, the operator formulation
$$u:=U_\theta\circ \cF(\phi),\ \ p:=U_\theta\circ \cF(q),$$
then leads to the operator version of Navier-Stokes (abbreviated as NS) equations on quantum Euclidean spaces:
\begin{align}\label{ens}
\lf\{\begin{array}{ll}
\partial_tu-\Delta_{\theta} u+A(u)+\nabla_{\theta} p=0;\\[5pt]
\mathbf{div}\, u=0;\\[5pt]
u(0)=u_0,
\end{array}\r.
\end{align}
where $\Delta_{\theta}=\sum^d_{i=j}\partial^2_j$, $\nabla_{\theta}=(\partial_1,\ldots,\partial_d)$, $\mathbf{div}\, u=\partial_1u_1+\dotsm+\partial_du_d$, $u=(u_1,\ldots, u_d)$, $p$ are unknown operator-valued functions from $[0,T)$ to operators, $T>0$, and $A(u)$ denotes the nonlinear term
$$A(u):=u\cdot\nabla_{\theta}u:
=\lf\{\sum_{j=1}^du_j(\partial_ju_k)\r\}_{k=1}^d,$$
whenever the operations are allowed. {We will give a sufficient large class for $(u,p)$ in Section \ref{s3} such that the above operations are all meaningful}.

Regarding the theory of well-posedness, the advantage of the operator version of Navier-Stokes equations \eqref{ens} is that: the product of two operators behaves better on $L_p$ spaces than the Moyal product of two functions. Hence we concern about the well-posedness of \eqref{ens} in this article. Let us point out that in the course of preparing the present paper and a note on the theory of  function spaces over quantum Euclidean spaces \cite{Hon}, McDonald \cite{m23} and Ruzhansky et al. \cite{RST,RST25} developed some nonlinear estimates or Sobolev type inequalities to study nonlinear evolution equations on quantum Euclidean spaces in an abstract way. Compared to these works, we will deeply exploit various  techniques from noncommutative harmonic analysis and function spaces to get a complete quantum analogues of the classical results; moreover, we not only fully recover the classical results ($\theta=0$), but also obtain some results on the semiclassical limit, that is, by using techniques independent of $\theta$ we will show the quantum solution $u_\theta$ converges to the classical solution $\phi$ in a proper way as $\theta\rightarrow 0,$ {which relates the quantum Navier-Stokes equation to the classical one quantitatively.}

As recalled before, if $\theta=0$, the above Navier-Stokes equation \eqref{ens} reduces to classical
Navier-Stokes equation on Euclidean space $\rr^d$. As is known to all,
classical Navier-Stokes equation is one of the most fundamental equations in the theory of fluid mechanics. In recent years, there has been a substantial amount of literature focusing on the well-posedness theory of the incompressible Navier-Stokes equations, see \cite{agz11,cz07,fjr72,fk64,gm85,hw14,kf62,l02,llt22,ln24,syc17,whhg11,w81,xlz16,z08}. Among them, let us mention the works that are closely related to the present paper, the one by Ladyzhenskaya \cite{l69}  where she established the global well-posedness for 2D Navier-Stokes equations, and the one by Kato  \cite{k84} where he not only provided the semigroup approach to the 2D case but also obtained the local well-posedness of solution in $L_d(\rr^d)$ in high-dimension case.

In the present paper, as the first step to fully understand the Navier-Stokes equations on the quantum Euclidean spaces, we will establish quantum analogues of Ladyzhenskaya and Kato's results.

In what follows, let $\dot{H}^1(\rd)$ denote the homogeneous Sobolev space and $[L_p(\rd)]_0^d$ ($1\leq p\leq\fz$) be the set of tuples of divergence free elements. All the notations below will be rigorously introduced in later sections.

\begin{thm}\label{t3.1}
Let $d\geq2$ and $u_0\in [L_d(\rr_{\theta}^d)]^d_0$. Then we have the following conclusions:
\begin{enumerate}
\item[\rm{(i)}]
There exists a maximal time $T_{u_0}>0$ such that the Navier-Stokes equation \eqref{ens} exists a unique smooth solution $u\in C([0,T_{u_0});[L_d(\rr_{\theta}^d)]^d_0)
\cap L_{d+2}^{\mathrm{loc}}([0,T_{u_0});[L_{d+2}(\rr_{\theta}^d)]^d_0)$;
if $T_{u_0}<\fz$, then we have $\|u\|_{L_{d+2}([0,T_{u_0});[L_{d+2}(\rr_{\theta}^d)]^d)}=\infty$.
Moreover, if $\|u_0\|_{[L_d(\rd)]^d}$ is sufficiently small, then $T_{u_0}=\infty$.

\item[\rm{(ii)}]
If $u_0\in [L_2(\rr_{\theta}^d)]^d_0\cap [L_d(\rr_{\theta}^d)]^d_0$, then the solution obtained in $(i)$ satisfies
$$ u\in C([0,T_{u_0});[L_2(\rr_{\theta}^d)]^d_0)\cap L_{2}^{\mathrm{loc}}([0,T_{u_0});[\dot{H}^1(\rr_{\theta}^d)]^d_0).$$
\end{enumerate}
\end{thm}

Due to the noncommutativity, the nonlinear term $A(u)$ in \eqref{ens} cannot ensure an energy identity and thus a good global well-posedness theory even if the initial datum $u_0$ is self-adjoint. Therefore, we will consider the Navier-Stokes equation with nonlinear term of symmetric form $S(u)$:
\begin{eqnarray}\label{ens1}
\lf\{\begin{array}{ll}
\partial_tu-\Delta_{\theta} u+S(u)+\nabla_{\theta} p=0;\\[5pt]
\mathbf{div}\, u=0;\\[5pt]
u(0)=u_0,
\end{array}\r.
\end{eqnarray}
where
$$S(u):=\frac{1}{2}\lf[u\cdot(\nabla_{\theta} u)+
((\nabla_{\theta} u)^T\cdot u^T)^T\r]:
=\lf\{\frac{1}{2}\sum_{j=1}^du_j(\partial_ju_k)+(\partial_ju_k)u_j\r\}_{k=1}^d.$$
Here, for any vector or matrix $\Lambda$, $\Lambda^T$ denotes the transpose of $\Lambda$.

\begin{thm}\label{t3.1x}
Let $d\geq2$ and $u_0\in [L_d(\rr_{\theta}^d)]^d_0$. Then we have the following conclusions:
\begin{enumerate}
\item[\rm{(i)}]
There exists a maximal time $T_{u_0}>0$ such that the Navier-Stokes equation \eqref{ens1} exists a unique smooth solution $u\in C([0,T_{u_0});[L_d(\rr_{\theta}^d)]^d_0)
\cap L_{d+2}^{\mathrm{loc}}([0,T_{u_0});[L_{d+2}(\rr_{\theta}^d)]^d_0)$;
if $T_{u_0}<\fz$, then we have $\|u\|_{L_{d+2}([0,T_{u_0});[L_{d+2}(\rr_{\theta}^d)]^d)}=\infty$.
Moreover, if $\|u_0\|_{[L_d(\rd)]^d}$ is sufficiently small, then $T_{u_0}=\infty$.

\item[\rm{(ii)}]
If $u_0\in [L_2(\rr_{\theta}^d)]^d_0\cap [L_d(\rr_{\theta}^d)]^d_0$, then the solution in (i) satisfies
$$ u\in C([0,T_{u_0});[L_2(\rr_{\theta}^d)]^d_0)\cap L_{2}^{\mathrm{loc}}([0,T_{u_0});[\dot{H}^1(\rr_{\theta}^d)]^d_0).$$
Moreover, if $u_0$ is self-adjoint, then the solution satisfies the energy identity
\begin{align}\label{energy00}
\frac{1}{2}\|u(t)\|^2_{[L_2(\rd)]^d}+\int^t_0\|\nabla_{\theta} u(s)\|^2_{[L_2(\rd)]^{d^2}}\,ds=\frac{1}{2}\|u_0\|^2_{[L_2(\rd)]^d}, \;\forall\;0<t<T_{u_0}.
\end{align}

\item[\rm{(iii)}]
If $d=2$ and $u_0$ is self-adjoint, then the Navier-Stokes equation \eqref{ens1} is globally well-posed in $C([0,\infty);[L_2(\ri)]^2_0)\cap L_2([0,\infty);$ $[\dot{H}^1(\ri)]^2_0)$, and the solution $u$ is smooth and satisfies the energy identity
\begin{align}\label{energy11}
\frac{1}{2}\|u(t)\|^2_{[L_2(\ri)]^2}+\int^t_0\|\nabla_{\theta} u(s)\|^2_{[L_2(\ri)]^4}\,ds=\frac{1}{2}\|u_0\|^2_{[L_2(\ri)]^2}, \;\forall\;0<t<\fz.
\end{align}
\end{enumerate}
\end{thm}

\begin{remark}{\rm
\begin{enumerate}
\item[\rm{(i)}] When $\theta=0$, Theorem \ref{t3.1x} recovers the classical results of Ladyzhenskaya \cite{l69} and Kato \cite{k84}.
\item[\rm{(ii)}] With the asymmetric nonlinear term $A(u)$, we establish the local well-posedness theory of \eqref{ens} in Theorem \ref{t3.1}, and do not know how to achieve  the global well-posedness theory. In the current paper, we only consider the global well-posedness of equation \eqref{ens1} with symmetric nonlinear term. Let us mention that in Theorem \ref{t3.1}, we do not require the self-adjointness of $u_0$.
\end{enumerate}
}\end{remark}

Since in the quantum setting the notion ``point" is not available anymore, one cannot formally  put forward the Navier-Stokes equations \eqref{ens} \eqref{ens1} in the same way as in classical case. We rigorously formulate the differential equation in \eqref{ns} with the help of the class of quantum tempered distribution and noncommutative $L_p$-spaces, and establish its equivalence with the integral form so that we are reduced to focusing on the mild solutions of \eqref{ens} \eqref{ens1} (see Lemma \ref{lem:mild}). We also introduce the notion of critical spaces in Section \ref{s3}, which is non-trivial in the setting of quantum Euclidean spaces.

As in the special case $\theta=0$, our approach to finding the mild solutions will mainly rely on the contraction mapping principle (Sections \ref{s6} and \ref{s7}). In order to apply this principle, we establish several sharp time-space estimates for both the linear and nonlinear terms in Section \ref{s6.1}. The sharp time-space estimates in turn involve the sharp $(L_r, L_p)$ estimates of the quantum heat semigroup (Proposition \ref{prop:rp}), Besov spaces, and various sharp embeddings between Besov spaces and Sobolev spaces (Section \ref{s5}). Due to the noncommutativity and the lack of the notion of ``points", there need several new techniques to achieve these sharp estimates; among them, a notable one is a transference principle (Theorem \ref{thm:trans}) which will not only show the $L_p$-boundedness of the Leray projection---Lemma \ref{t5.1}\ but also yield the sharp embedding---Lemma \ref{Besov embedding}.

It should be pointed out that McDonald \cite{m23} studied the well-posedness of the  equations \eqref{ens1} when $\mathrm{det}(\theta)\neq0$, the initial datum $u_0$ is self-adjoint and belongs to the noncommutative Sobolev
space $H^2(\rd)$ that does not admit any scaling symmetry; one key fact used by him is that when $\mathrm{det}(\theta)\neq0$, the resulting noncommutative $L_p$-spaces are the Schatten classes, and thus the Sobolev embeddings trivially hold but with the constant depending on $\theta$. Our Theorem \ref{t3.1x} goes further beyond McDonald's result since we work with $u_0\in L_d(\rd)$---the critical space; moreover, our techniques, which are independent of $\theta$, allow to derive the semiclassical limit of quantum Navier-Stokes equations, which reveals below the relationship between the solutions of quantum Navier-Stokes equations and those of classical ones.

Given an antisymmetric matrix $\theta$, let $\|\theta\|$ denote the maximal eigenvalue of $\theta$. In what follows, $\mathcal F$ denotes the Fourier transform on $\mathbb R^d$.

\begin{thm}\label{semiclassical}
 For a fixed $\phi_0\in [L_2(\bR^2)]^2_0,$ let \( u_\theta \in C([0,\fz);[L_2(\ri)]^2_0) \cap L_2([0,\fz); [\dot{H}^1(\ri)]^2_0) \) be the unique smooth solution to equation \eqref{ens1} with initial datum \( u_{\theta,0}=U_\theta\circ \cF(\phi_0) \).  Then as \( \|\theta\| \to 0 \), we have \( \mathcal{F}^{-1} \circ U^{-1}_\theta(u) \to \phi \) in the weak-\( \ast \) topology of \( L_\infty([0,\fz),[L_2(\mathbb{R}^2)]^2_0) \), where \( \phi \in C([0,\fz);[L_2(\mathbb{R}^2)]^2_0) \cap L_2([0,\fz); [\dot{H}^1(\mathbb{R}^2)]^2_0) \) is the unique solution to the classical Navier-Stokes equation with initial datum \( \phi_0\).
\end{thm}

For the equation \eqref{ens} with asymmetric nonlinear term, we have no idea how to produce a similar result due to the lack of energy identity. More explanations can be found in Remark \ref{reason}.

In the case of higher dimensions, we derive the similar conclusion.
 For convenience, we define for $0<T\leq\fz$,
\begin{align*}
    \mathcal{N}_T^{\theta}:=C([0,T);[L_2(\mathbb{R}_{\theta}^d) \cap L_d(\rd)]^d_0) \cap L_{2}^{\mathrm{loc}}([0,T);[\dot{H}^1(\rd)]_0^d)
    \cap L_{d+2}^{\mathrm{loc}}([0,T);[L_{d+2}(\mathbb{R}_{\theta}^d)]^d_0).
\end{align*}

\begin{thm}\label{semi2}
    Let \( d > 2 \) and $\phi_0$ be a function satisfying one of the following assumptions:
    \begin{enumerate}
        \item[\rm{(i)}] \( \phi_0\in [L_2(\bR^d)\cap \cF^{-1}(L_{d'-\varepsilon}(\bR^d)]^d_0\) for some $0<\varepsilon\le d'-1,$ where $d'=\frac d{d-1};$
        \item[\rm{(ii)}]\( \phi_0\in [L_2(\bR^d)\cap \cF^{-1}(L_{d'}(\bR^d)]^d_0\) with $\|\cF\phi_0\|_{[L_{d'}(\bR^d)]_0^d}$ sufficiently small.
    \end{enumerate}
Then there exists one time $0<T_{\phi_0}\leq \infty$ such that as \( \|\theta\| \to 0 \), the unique smooth solution   \( u_\theta \in \mathcal{N}_{T_{\phi_0}}^\theta \) to \eqref{ens} or \eqref{ens1} with initial datum \( u_{\theta,0}=U_\theta\circ \cF(\phi_0) \) converges to $\phi$ in the weak-\( \ast \) topology of \( L_\infty([0,T);[L_2(\mathbb{R}^d) \cap L_d(\mathbb{R}^d)]_0^d)\) for every $T<T_{\phi_0}$, where \(\phi\in \mathcal{N}_{T_{\phi_0}}^0 \) is  the unique solution to the classical Navier-Stokes equation with initial datum \( \phi_0 \).
\end{thm}

\bigskip

\noindent {\bf Notations.}
Conventionally, we set $\nn:=\{0,\,1,\, 2,\,\ldots\}$
and $\rr_+:=(0,\,\infty)$.
Throughout the whole paper, we denote by $C$ a positive constant which is independent of the main parameters, but it may
vary from line to line. We use $A\lesssim B$ to denote the statement that $A\leq CB$ for some constant $C>0$, and $A\thicksim B$ to denote the statement that $A\lesssim B$ and $B\lesssim A$.
For any $1\leq p\leq\infty$, we denote by $p'$ the conjugate of $p$, which satisfies $\frac{1}{p}+\frac{1}{p'}=1$.
We also use $C_{\alpha,\,\beta,\,\ldots}$ to denote a positive
constant depending on the indicated parameters $\alpha,\,\beta,\,\ldots$.
For a Banach space $X$, $C(I;X)$ denotes all continuous functions on the interval $I\subset\mathbb{R}$ with value in $X$.
For the Banach spaces $X$, $Y$, let $\mathcal{B}(X,\,Y)$ denote the all bounded linear operators from $X$ to $Y$.
 If $X=Y$, then we set $\mathcal{B}(X):=\mathcal{B}(X,\,X)$.

\bigskip



\bigskip
\section{Preliminaries} \label{s2}
\setcounter{equation}{0}
\bigskip
 In this section, let us recall some basic definitions and properties of some function spaces, including noncommutative $L_p$-spaces, Lorentz spaces, quantum Euclidean spaces, and among others.
\subsection{Noncommutative $L_p$-spaces, Lorentz spaces}\label{n2.0}
Let us recall the definitions and some basic properties of noncommutative $L_p$-spaces and Lorentz spaces (see \cite{PX} for more details about noncommutative integration theory). Let $\mathcal M$ be a von Neumann algebra equipped with a normal semifinite faithful (abbreviated as {\it n.s.f}) trace $\tau$. To begin with, let $S^+_{\cm}$ be the set of all positive element $x\in\cm$ such that
$$\tau(s(x))<\infty,$$
where $s(x)$ denotes the least projection $e\in\cm$, called the support of $x$, such that
$exe=x.$
Let $S_\cm$ be the linear span of $S_{\cm}^+$. For any $p\in(0,\,\infty)$, we define
$$\|x\|_{L_p(\cm)}:=(\tau(|x|^p))^{\frac 1p}, \ \ \ x\in S_\cm,$$
where $|x|:=(x^*x)^{\frac 12}$. We define the {\it noncommutative $L_p$-space} associated with $(\cm,\,\tau)$, denoted by $L_p(\cm)$, to be the completion
of $(S_{\cm}, \|\cdot\|_{L_p(\cm)})$. For convenience, we usually set $L_{\infty}(\cm)=\cm$ equipped
with the operator norm $\|\cdot\|_{\cm}$.
 As classical $L_p$-spaces, the noncommutative
$L_p$-spaces possess the basic properties such as the duality and the interpolation etc..

In addition to the above definition, we also know that the elements in $L_p(\cm)$ can be described as closed densely defined operators on
$H$, where $H$ is the Hilbert space on which $\cm$ acts. A closed densely defined operator $x$
on $H$ is said to be affiliated to $\cm$ if $ux=xu$ for any unitary $u$ in the commutant
$\cm'$ of $\cm$. We say that $x$ affiliated to $\cm$ is {\it $\tau$-measurable} with respect to
$(\cm,\,\tau)$ (or simply measurable) if for any $\delta>0$, there exists a projection $e\in \mathcal{B}(H)$
such that
$$e(H)\subset \mathrm{Dom}(x)\ \ \ \mathrm{and} \ \ \ \tau(e^{\bot})\leq\delta,$$
where $\mathrm{Dom}(x)$ denotes the domain of the operator $x$.

In what follows,
we denote the $*$-algebra of $\tau$-measurable operators by $L_0(\cm)$. For $0<p<\fz$ and $0<q<\fz$,  the {\it noncommutative Lorentz space} $L_{p,\,q}(\cm)$ is defined as the set of all $x\in L_0(\cm)$ such that
$$\|x\|_{L_{p,\,q}(\cm)}:=\lf[\int_0^{\fz}\lf(t^{\frac 1p}\mu_t(x)\r)
^q\frac{dt}{t}\r]^{\frac 1q}<\fz.$$
For $q=\fz$, the space $L_{p,\,\fz}(\cm)$ is usually called a weak $L_p$-space
with $0<p<\fz$, with its quasi-norm defined as
$$\|x\|_{L_{p,\,\fz}(\cm)}:=\sup_{s>0}s\lambda_s(x)^{\frac 1p}.$$
In the above, for all $t>0$, $s>0$,
$\lambda_s(x):=\tau(e_s^{\bot}(|x|))$, $\mu_t(x):=\inf\{s>0:\lambda_s(x)\leq t\}$, and $e^{\bot}_s(|x|):=\chi_{(s,\,\fz)}(|x|)$ is the spectrum projection of $|x|$ corresponding to the
interval $(s,\,\fz)$.

Like the classical Lorentz spaces, the noncommutative Lorentz space $L_{p,\,q}(\cm)$ also has a number of analogous properties as follows:
\begin{remark}{\rm
\begin{enumerate}
\item[\rm{(i)}]
If $p=q$, then  $L_{p,\,p}(\cm)=L_{p}(\cm)$.
\item[\rm{(ii)}]
Obviously, when $0<q_1\leq q_2<\fz$, we have $L_{p,\,q_1}(\cm)\subset L_{p,\,q_2}(\cm)$.
\end{enumerate}}
\end{remark}

The following result is the real interpolation of noncommutative Lorentz spaces, see \cite{PX}.

\begin{lemma}\label{lor1}
Let $0<p_1,p_2\leq\fz$ with $p_1\neq p_2$, $0<\eta<1$ and $0<q\leq\fz$. Then we have
$$\lf[L_{p_1}(\cm),\,L_{p_2}(\cm)\r]_{\eta,\,q}=
L_{p,\,q}(\cm)$$
with equivalent quasi-norms, where $1/p=(1-\eta)/p_1+\eta/p_2$.
\end{lemma}

\subsection{Quantum Euclidean spaces}\label{n2.1}
As in \cite{hlw23,LSZ20,MSX20} and so on, we will recall the definition of quantum Euclidean spaces $\rd$.

\begin{definition}
Let $\theta$ be a $d\times d$ antisymmetric matrix and $t\in\mathbb{R}^d$. We define the unitary
operator $\lambda_\theta(t)$ on $L_2(\mathbb{R}^d)$:
\begin{align}\label{wll}
		(\lambda_\theta(t)f)(\xi):=e^{-\frac{\mathrm{i}}{2}(t,\theta \xi)}f(\xi-t),\quad f\in{L_2(\mathbb{R}^d)}, \, \xi\in\mathbb{R}^d,
	\end{align}
where $(\cdot,\,\cdot)$ denotes the usual inner product in $\rr^d$, and $\mathrm{i}:=\sqrt{-1}$.
We define the quantum Euclidean space $\rd$ to be
a closed subalgebra generated by  $\{\lambda_\theta(t)\}_{t\in\mathbb{R}^d}$ of $\mathcal{B}(L_2(\mathbb{R}^d))$
 with respect to the weak operator topology.
\end{definition}

It can be demonstrated that the family $\{\lambda_{\theta}(t)\}_{t\in\rr^d}$ satisfies
the following Weyl relation:
$$ \lambda_\theta(t)\lambda_\theta(s)=e^{{\mathrm{i}}(t,\theta s)}\lambda_\theta(s)\lambda_\theta(t),\quad \text{for\: all}\: t,s\in\mathbb{R}^d.$$
The above relation is known as the Weyl representation of the canonical
commutation relation.

\begin{remark}{\rm
\begin{enumerate}
\item[\rm{(i)}]
  In the case $\theta=0$, the quantum Euclidean space $\rd$ is the von Neumann algebra generated by the unitary group of translations on $\mathbb{R}^d$, which is
  $*$-isomorphic to $L_\infty(\mathbb{R}^d)$.
\item[\rm{(ii)}] It is easy to verify that the family $\{\lambda_\theta(t)\}_{t\in\mathbb{R}^d}$ is strongly continuous.
\end{enumerate}}
\end{remark}

Now we introduce a map from $L_1(\mathbb{R}^d)$ to $L_{\fz}(\rd)$ which is denoted by $U_\theta$  as below:
Let $f \in L_1(\mathbb{R}^d)$, one defines  the Weyl transform $U_\theta(f): L_2(\mathbb{R}^d)\rightarrow L_2(\mathbb{R}^d)$ as
\begin{equation}\label{defU}
 U_\theta(f)( g):=\int_{\mathbb{R}^d}f(\xi)(\lambda_\theta(\xi)g)\,d\xi
\end{equation}
for $g\in L_2(\mathbb{R}^d)$.
This $L_2(\mathbb{R}^d)$-valued integral is convergent in the Bochner sense.

In what follows, we normalize
the Fourier transform of a reasonable function $f$ as
$$\mathcal{F}(f)(\xi):=\hat{f}(\xi) := \int_{\mathbb{R}^d}
f(t)e^{-\mathrm{i}(t,\,\xi)}\,dt, \ \ \  \forall\,\xi\in\mathbb{R}^d,$$
where $\mathrm{i}:=\sqrt{-1}$.
Define the inverse Fourier transform of $f$ by
$$\mathcal{F}^{-1}(f)(t):=\check{f}(t):=\frac{1}{(2\pi)^{d}}\int_{\mathbb{R}^d}
f(\xi)e^{\mathrm{i}(\xi,\,t)}\,d\xi, \ \ \  \forall\,t\in\mathbb{R}^d.$$

The image of $\mathcal{S}(\mathbb{R}^d)$, the Schwartz class on $\mathbb R^d$, under the map $U_\theta$ is called the class of Schwartz functions on $\rd $ :
$$ \mathcal{S}(\rd):=\lf\{x\in L_{\fz}(\rd): x=U_\theta(f), \text{\:for\: some\;} f\in\mathcal{S}(\mathbb{R}^d)\r\}.$$
Then $U_\theta$ is a bijection from $\mathcal{S}(\mathbb{R}^d)$ to $ \mathcal{S}(\rd)$,
and thus $\mathcal{S}(\rd)$ is a Fr\'{e}chet topological space equipped with the Fr\'{e}chet topology induced by $U_\theta$.
The topological dual of $\mathcal{S}(\rd)$ is denoted as $\mathcal{S}'(\rd)$, then $U_\theta$ extends to a bijection from $\mathcal{S}'(\mathbb{R}^d)$  to $\mathcal{S}'(\rd)$: for $f\in\mathcal{S}^\prime(\mathbb{R}^d)$,
$$ \lf\langle U_\theta(f),\,U_\theta(g)\r\rangle:=\langle f,\,\tilde{g}\rangle, \ \ \ \text{for\:all\;}g\in \mathcal{S}(\mathbb{R}^d),$$
where $\widetilde{g}(\cdot):=g(-\cdot)$.

If $x\in \mathcal{S}(\rd)$ is given by $x=U_\theta(f)$ for $f\in \mathcal{S}(\mathbb{R}^d)$, we define $\tau_\theta(x):=f(0)$, then $\tau_\theta$ extends to a $n.s.f.$ trace on $\rd$. The noncommutative $L_p$-space associated to $(\rd, \tau_\theta)$
 is denoted by $L_p(\rd)$. The space $ \mathcal{S}(\rd)$ is dense in $L_p(\rd)$ for $1\leq p<\infty$ with respect to the norm $\|\cdot\|_{L_p(\rd)}$, and dense in $L_\infty(\rd)$ in the weak-$*$ topology.
 We refer the reader to \cite{gjp21, MSX20} for more information.

\begin{remark}{\rm
When $\det(\theta)\neq0$, it is known that $L_p(\rd)$ coincides with the Schatten $p$-class and we have the following embedding:
$$L_p(\rd)\subset L_q(\rd), \ \ \  \mathrm{if} \ p\le q.$$}
\end{remark}

The following Hausdorff-Young inequality should be well-known to the experts (cf. e.g. \cite[Lemma 2.7]{hlw23} or \cite[Proposition 2.10]{MSX20}).
\begin{lemma}\label{l2.a1}
Let $f\in\cs(\rr^d)$. Then we have
\begin{enumerate}
\item[\rm{(i)}]
$$\|U_{\theta}(f)\|_{L_2(\rd)}=\|f\|_{L_2(\rr^d)};$$
\item[\rm{(ii)}] for $p\in[1,\,2)$,
$$\|U_{\theta}(f)\|_{L_{p'}(\rd)}\leq\|f\|_{L_p(\rr^d)}.$$
\end{enumerate}
Therefore, $U_\theta$ extends to a contraction from $L_p(\rr^d)$ $(p\in[1,\,2))$ to $L_{p'}(\rd)$ and
an isometry on $L_2(\rd)$.
\end{lemma}

For a reasonable function  $\psi:\rr^d\rightarrow\mathbb{C}$ and $x=U_{\theta}(f)$, $f\in\cs(\rr^d)$, define the Fourier multiplier $T_\psi(U_{\theta}(f)):=U_{\theta}(\psi f)$. The following Young inequality will be instrumental in establishing the embedding properties.

\begin{lemma}\label{lem:psi}
Let $1\leq p, q, r\leq \infty$ satisfy $1+\frac 1r=\frac1q+\frac1p$. If the function $\psi$ satisfies that $\check{\psi}\in L_q(\rr^d)$, then we have
\begin{align}\label{psi}
\lf\|T_{\psi}x\r\|_{L_r(\rd)}\leq \lf\|\check{\psi}\r\|_{L_q(\rr^d)}\|x\|_{L_p(\rd)},\ \ x\in L_p(\rd).
\end{align}
\end{lemma}
\begin{proof}
If $(p,q,r)=(1,\fz,\fz)$, then the inequality \eqref{psi} is just from \cite[Lemma 5.1]{hlw23}.
If $(p,q,r)=(1,1,1)$ or $(\fz,1,\fz)$, from the triangle inequality, we immediately obtain the inequality \eqref{psi}.
By the above fact and the multi-linear interpolation theory (see \cite[Theorem 4.4.2]{bl76}), we conclude the inequality \eqref{psi} in the whole scales.
\end{proof}

\begin{remark}
{\rm A completely bounded version of Lemma \ref{lem:psi} holds still true (cf. \cite{Hon}), which allows us to deduce the Sobolev embedding in the operator space category; this provides further applications to functional inequalities and quantum information (cf. \cite{GJX}).}
\end{remark}

\subsection{Sobolev spaces on quantum Euclidean spaces}\label{n2.2}

Now, let us recall the differential structure on quantum Euclidean space $\rd$ (cf. \cite{LSZ20,MSX20, Hon, gjp21}). For $1\leq j\leq d$, let $\mathbf{D}_j$ be the multiplication operator defined as
$$\mathbf{D}_j(f)(t) = t_jf(t), \ \ \ t=(t_1,\,t_2,\ldots,t_d)\in\rr^d$$
on the domain $\mathrm{dom}(\mathbf{D}_j)=\{f\in L_2(\rr^d):f\in L_2(\rr^d,\,t_j^2dt)\}$. Given a fixed $s=(s_1,\,s_2,\ldots,s_d)\in\rr^d$,
it is straightforward to verify that the unitary generator $\lambda_{\theta}(s)$ preserves $\mathrm{dom}(\mathbf{D}_j)$. We can compute the following:
$$\lf[\mathbf{D}_j,\,\lambda_{\theta}(s)\r]=s_j\lambda_{\theta}(s)$$ and
$$e^{\mathrm{i}t\mathbf{D}_j}\lambda_{\theta}(s)e^{-\mathrm{i}t\mathbf{D}_j}
=e^{\mathrm{i}ts_j}\lambda_{\theta}(s)\in L_{\fz}(\rd), \ \ \ t>0.$$
For a general element $x\in L_{\fz}(\rd)$, if $[\mathbf{D}_j,\,x]$ extends to a bounded operator on $L_2(\rr^d)$, then we write
$$\mathrm{i}\lf[\mathbf{D}_j,\,x\r]=\lim_{t\rightarrow0}
\frac{e^{\mathrm{i}t\mathbf{D}_j} x e^{-\mathrm{i}t\mathbf{D}_j}-x}{t}$$
with respect to the strong operator topology. Consequently, $\mathrm{i}[\mathbf{D}_j,\,x]\in L_{\fz}(\rd)$. Furthermore, with respect to the
$x$ mentioned above, the operator $\mathrm{i}[\mathbf{D}_j,\,x]$ is defined as the partial derivative of $x$, which is denoted as $\partial_jx$.

For a multi-index $\alpha:=(\alpha_1,\,\alpha_2,\ldots,\alpha_d)\in \zz_+^d$ and $x\in L_{\fz}(\rd)$, if every iterated commutator $[\mathbf{D}^{\alpha_j}_j,\,[\mathbf{D}^{\alpha_{j+1}}_{j+1},\,
\ldots,[\mathbf{D}_d
^{\alpha_d},\,x]]], j=1,\ldots,d$ extends to a bounded operator on $L_2(\rr^d)$, then the mixed
partial derivative $\partial^{\alpha}x$ is defined as
$$\partial^{\alpha}x=\mathrm{i}^{|\alpha|}\lf[\mathbf{D}_1^{\alpha_1},\,\lf[\mathbf{D}_2
^{\alpha_2},\ldots,\lf[\mathbf{D}_d^{\alpha_d},\,x\r]\r]\r].$$

\begin{remark}{\rm
Notice that, on Schwartz
functions, the partial derivative $\partial_{j}$ can be defined as a multiplier in terms of the map $U_{\theta}$. It can be easily verified that
$$\partial_{j}U_{\theta}(f)=U_{\theta}(\mathrm{i}\mathbf{D}_j(f))$$
for any $f \in\cs(\rr^d)$ and $j=1,\ldots,d.$}
\end{remark}

Define the quantum Laplace operator $\Delta_{\theta}$ by
$$\Delta_{\theta}:=\sum_{j=1}^d\partial_{j}^2,$$
and the quantum gradient operator $\nabla_{\theta}$ by
$$\nabla_{\theta}:=(\partial_1,\partial_2,\ldots,\partial_d).$$

Based on the above definitions, now we can define the derivative over $\cs'(\rd)$, the space of tempered distributions.
\begin{definition} For any muti-index $\alpha\in\zz^d_+$ and $T\in\cs'(\rd)$,
we define the partial derivative $\partial^{\alpha}$ of the distribution $T\in\cs'(\rd)$ as follows:
$$\langle\partial^{\alpha}T,\,\varphi\rangle=(-1)^{|\alpha|}\langle T,\,\partial^{\alpha}\varphi\rangle,
\ \ \ \mathrm{for\ any} \  \varphi\in\cs(\rd).$$

\end{definition}

\begin{remark}{\rm It should be noted that the operators $\Delta_{\theta}$ and $\nabla_{\theta}$ are independent of the matrix $\theta$.
Nevertheless, within the context of this paper, we choose to employ the notation with $\theta$ in order to distinguish $\Delta_{\theta}$, $\nabla_{\theta}$ from the classical
Laplacian $\Delta$ and gradient operator $\nabla$, respectively.}
\end{remark}

In what follows, we will frequently make reference to the operators $(1-
\Delta_{\theta})^{\frac{1}{2}}$ and $(-\Delta_{\theta})^{\frac{1}{2}}$.
Specifically, they are the operators on $L_2(\rr^d)$ corresponding to pointwise multiplication by $(1+|t|^2)^{\frac{1}{2}}$ and $|t|$, respectively, where $t\in\rr^d$.
Classically, the operator $(1-\Delta_{\theta})^{\frac{1}{2}}$ and $(-\Delta_{\theta})^{\frac{1}{2}}$ are known as the Bessel potential and the Riesz potential respectively.
\begin{definition} Let $1\leq  p<\fz$, $k\in\nn$ and $s\in\rr_+$.
The $k$-th order Sobolev space on $\rd$  is defined as
$$ W_p^k(\rd):=\lf\{x\in\mathcal{S}'(\rd): \partial^mx\in L_p(\rd), \ \mathrm{for \  any} \  m\in\zz_+^d \ \mathrm{with} \ |m|\leq k\r\}$$
equipped with the norm
$$\|x\|_{W_p^k(\rd)}:=\lf(\sum_{0\leq|m|\leq k}\lf\|\partial^mx\r\|_{L_p(\rd)}^p\r)^{\frac{1}{p}}.$$
The Bessel potential Sobolev space $H_p^s(\rd)$ is defined as the subset of $x\in \mathcal{S}'(\rd)$ such that
$(1-\Delta_{\theta})^{\frac s2}x\in L_p(\rd)$, equipped with the norm
$$\|x\|_{H_p^s(\rd)}:=\lf\|(1-\Delta_{\theta})^{\frac s2}x\r\|_{L_p(\rd)}.$$
The Riesz potential Sobolev space $\dot{H}_p^s(\rd)$ is defined as the subset of $x\in \mathcal{S}'(\rd)$ such that
$(-\Delta_{\theta})^{\frac s2}x\in L_p(\rd)$, equipped with the norm
$$\|x\|_{\dot{H}_p^s(\rd)}:=\lf\|(-\Delta_{\theta})^{\frac s2}x\r\|_{L_p(\rd)}.$$
\end{definition}
The following lemma concerns the relationship between $H_p^k(\rd)$ and $W_p^k(\rd)$. Since its proof is similar to that of \cite[Theorem 2.9]{XXY18}, we omit the details here.
\begin{lemma}\label{hw}
Let $1<p<\fz$ and $k\in\nn$. Then $$H_p^k(\rd)=W_p^k(\rd)$$ with equivalent norms.
\end{lemma}

\begin{remark}{\rm
Throughout this paper, for the sake of brevity and convenience, we denote $$\dot{H}^s(\rd):=\dot{H}_2^s(\rd),$$
for any $s>0.$
}
\end{remark}

Here we briefly present Sobolev spaces needed for the present paper's purpose, and a rather complete investigation of Sobolev spaces over quantum Euclidean spaces can be found in \cite{Hon}.

\section{Navier-Stokes equations on quantum Euclidean spaces}\label{s3}
\setcounter{equation}{0}

\bigskip
In this section, we will introduce and study the Cauchy problem for the incompressible Navier-Stokes equations in the  framework of quantum Euclidean spaces.

\subsection{The definition of Navier-Stokes equations on $\rd$}\label{n2.3}

The following notation will be frequently used when referring to time-dependent objects. Let
$\mathfrak{M}$
denote the set of all operator-valued functions $\nu$ from $\mathbb{R}_+$ to $ L_1(\rd)+L_\infty(\rd)$ that satisfy the following condition: for any $t\in\rr_+$ and  $j=1,\dots,d$, the partial derivative $\partial_j\nu(t)\in L_1(\rd)+L_\infty(\rd)$.

Now, we turn our attention to the Cauchy problem for the incompressible Navier-Stokes equation
\begin{eqnarray}\label{ns}
\lf\{\begin{array}{ll}
\partial_tu-\Delta_{\theta} u+F(u)+\nabla_{\theta} p=0;\\[5pt]
\mathbf{div}\, u=0;\\[5pt]
u(0)=u_0,
\end{array}\r.
\end{eqnarray}
where $u=(u_1,\ldots, u_d)$ with $u_j\in \mathfrak M,j=1,\dots,d$, $\Delta_{\theta}=\sum^d_{i=j}\partial^2_j$, $\nabla_{\theta}=(\partial_1,\ldots,\partial_d)$, $\mathbf{div}\, u=\partial_1u_1+\dotsm+\partial_du_d$, $p$ is unknown operator-valued function from $[0,T)$ to $ \cS'(\rd)$, $T>0$, and
$F(u)$ denotes the nonlinear term of asymmetric form:
$$A(u):=u\cdot\nabla_{\theta} u:=
\lf\{\sum_{j=1}^du_j(\partial_ju_k)\r\}_{k=1}^d$$
or the nonlinear term of symmetric form:
$$S(u):=\frac{1}{2}\lf[u\cdot\nabla_{\theta} u+
((\nabla_{\theta} u)^T\cdot u^T)^T\r]:=
\lf\{\sum_{j=1}^du_j(\partial_ju_k)+(\partial_ju_k)u_j\r\}_{k=1}^d,$$
where, for any vector or matrix $\Lambda$, $\Lambda^T$ denotes the transpose of $\Lambda$.
\begin{remark}\label{r0}
{\rm
Note that, under the condition $\mathbf{div}\, u=0$, the nonlinear term
$A(u)$ and $S(u)$ can be expressed in the following forms:
$$ A(u):=u\cdot\nabla_{\theta}u={\bf div}\,[u\otimes u]
$$
and
$$ S(u):=\frac{1}{2}\lf[u\cdot(\nabla_{\theta} u)+
((\nabla_{\theta} u)^T\cdot u^T)^T\r]
=\frac{1}{2}{\bf div}\lf[u\otimes u+(u\otimes u)^T\r],
$$
where ${\bf div}\,\Lambda:=\nabla_{{\theta}}\cdot\Lambda$, for any matrix-valued function $\Lambda$.
For the sake of convenience, these notations will be frequently used in the following content.
}\end{remark}

Like the classical Navier-Stokes equations, we introduce the definitions of strong solutions and  weak solutions.
\begin{definition}For an operator-valued function $u\in \mathfrak M$,
\begin{enumerate}
\item[\rm{(i)}]$u$ is called a strong solution of \eqref{ns} on the Banach space $X\subset L_1(\rd)+L_\fz(\rd)$, if $u$ satisfies the equality in the following sense:
\begin{align*}
\lim_{h\to 0}\lf\|\frac{u(t+h)-u(t)}{h}-\Delta_{\theta} u+F(u)+\nabla_{\theta} p\r\|_X=0;
\end{align*}
\item[\rm{(ii)}]$u$ is called a weak solution of \eqref{ns} if for every $\phi\in \cS(\rd)$, the function $\lag u,\phi\rag:t\mapsto \lag u(t),\phi\rag$ is differentiable and
\begin{align*}
\partial_t\lag u,\phi\rag-\lag u,\Delta_{\theta} \phi\rag+\lag F(u),\phi\rag-\lag p,\nabla_{\theta} \phi\rag=0,
\end{align*}
\end{enumerate}
where $F(u)=A(u)$ or $S(u)$.
\end{definition}

 Now we introduce the definition of smoothness on quantum Euclidean spaces.
 \begin{definition}
Let $u\in L_\infty((0,T);L_d(\rd))$. We say that $u$ is infinitely smooth if
$$\partial_t^{n}u\in C((0,T);H_d^k(\rd)),\ \ \mathrm{for\ all} \  n\in\mathbb{N}, \ k\in\nn.$$
\end{definition}
In terms of special structure, we claim that, if the initial datum $u_0$ is self-adjoint,
then the smooth solution $u$ of \eqref{ens1} satisfies the following conservation law:
\begin{align}\label{energy}
\frac{1}{2}\|u(t)\|^2_{[L_2(\rd)]^d}+\int^t_0\lf\|\nabla_{\theta} u(s)\r\|^2_{[L_2(\rd)]^{d^2}}\,ds=\frac{1}{2}\|u_0\|^2_{[L_2(\rd)]^d}.
\end{align}
Indeed, by multiplying the first equation of \eqref{ens1} by $u$ and applying integration by parts, we derive the following results:
$$\tau_\theta\int^t_0(\partial_su)(s) u(s)\,ds=\frac{1}{2}\lf(\|u(t)\|^2_{[L_2(\rd)]^d}-\|u_0\|^2_{[L_2(\rd)]^d}\r),$$
$$\tau_\theta\int^t_0\lf(-\Delta_{\theta} u(s)\r) u(s)\,ds=\int^t_0\lf\|\nabla_{\theta} u(s)\r\|^2_{[L_2(\rd)]^{d^2}}\,ds,$$
$$\tau_\theta\int^t_0\lf(u(s)\cdot\nabla_{\theta} u(s)\r) u(s)\,ds=-\tau_\theta\int^t_0\lf[\lf(\nabla_{\theta} u(s)\r)^T\cdot u(s)^T\r]^T u(s)\,ds=0,$$
and
$$\tau_\theta\int^t_0(\nabla_{\theta} p(s)) u(s)\,ds=-\int^t_0p(s)(\mathbf{div} \,u(s))\,ds=0.$$
Then, by using the tracial property in the third equality, we obtain the claim \eqref{energy}.

\subsection{An equivalent form of Navier-Stokes equations on $\rd$}
Let $(u,p)$ be a solution of the NS equation. Given the condition $\mathbf{div}\, u=0$, when we take the divergence of the first equation in \eqref{ns}, we obtain
\begin{align}\label{aft div}
\Delta_{\theta} p+{\bf div}\lf(A(u)\r)=0.
\end{align}
This, in turn, implies
\begin{align*}
\nabla_{\theta} p=\lf(-\Delta_{\theta}\r)^{-1}\nabla_{\theta}\,{\bf div} \lf(A(u)\r).
\end{align*}
Then the resulting operator is known as the Leray projection:
\begin{align}\label{def of P}
\mathbb P:=I+(-\Delta_{\theta})^{-1}\nabla_{\theta}\,{\bf div}.
\end{align}
Let $u$ be the solution of \eqref{aft div}. Substituting it into \eqref{ns}, we get
\begin{eqnarray} \label{ns1}
\left\{
\begin{array}{ll}
    \partial_tu-\Delta_{\theta} u+\mathbb P\lf(A(u)\r)=0;
  \\[5pt]
  u(0)=u_0;\\[5pt]
  \mathbf{div}\, u=0.
\end{array}
\right.
\end{eqnarray}

At this point, it becomes clear that the NS equation is essentially a nonlinear parabolic equation. Moreover, if
we replace the nonlinear term $A(u)$ with $S(u)$, the above conclusion also holds.

\subsection{The equivalence between differential and integral forms of Navier-Stokes equations}

In this subsection, to prove Theorem \ref{t3.1}, we will deduce an equivalent integral form of Navier-Stokes equations.

Firstly, let us recall certain definitions and establish several technical lemmas.
We begin with the definition of the heat semigroup on $\rd$.
Our focus will be on semigroups of operators on $\cs(\rd)$, specifically the heat semigroups denoted by
$t\rightarrow e^{t\Delta_{\theta}}$.
These operators can be defined either through functional calculus on $\rd$, or equivalently as Fourier multipliers. For any $x=U_{\theta}(f)$, $f\in\cs(\rr^d)$,  we define
$$e^{t\Delta_{\theta}}U_{\theta}(f):=U_{\theta}(e^{-t|\cdot|^2}f).$$

Now, we introduce some properties of the heat semigroup $\{H(t)\}_{t>0}$.

\begin{lemma}\label{heat}
Let $p\in[1,\,\fz]$.
Then the operator $\{H(t)\}_{t>0}$ has the following properties:
\begin{enumerate}
    \item[\rm{(i)}] Let $t>0$. $H(t)$ is bounded on $L_p(\rd)$, with norm
$$\|H(t)\|_{L_p(\rd)\rightarrow L_p(\rd)}\leq1.$$
For $p\in[1,\,\fz)$, $H(t)$ is strongly continuous on $L_p(\rd)$, in the sense that the mapping
$$[0,\,\fz)\times L_p(\rd)\rightarrow L_p(\rd),\ \ (t,\,x)\mapsto H(t)x$$
is continuous.
\item[\rm{(ii)}]If $p\in(1,\,\fz)$, then for every $x\in L_p(\rd),t>0,$ the integral $\int_0^t H(s)x\, ds$ satisfies \begin{align*}
    \Delta_{\theta}\int_0^t H(s)x\, ds=\int_0^t \Delta_{\theta} H(s)x\, ds=H(t)x-x.
\end{align*}
\item[\rm{(iii)}]If $x\in W_p^2(\rd),$ then
\begin{align*}
    \lim_{t\to 0^+}\lf\|\frac{H(t)x-x}{t}-\Delta_{\theta} x\r\|_{L_p(\rd)}=0.
\end{align*}
\end{enumerate}
\end{lemma}
\begin{proof}
From \cite[p.51]{jmpx21}, we know that the heat semigroup $\{H(t)\}_{t>0}$ is a Markov semigroup. Therefore, the conclusion (i) holds true.
Furthermore, combining Proposition 2.1.4 in \cite{Lunardi} with Proposition 5.4 in \cite{jlx06},
one deduces the conclusions (ii) and (iii).

\end{proof}

To solve the NS equation, we first consider the Cauchy problem for the heat equation:
$$\partial_tu-\Delta_{\theta} u=\widetilde{F}(u), \;u(0)=u_0.$$
By solving the corresponding ordinary differential equation, we can obtain the mild form of the heat equation. This mild form is given by:
\begin{align}\label{mild form}
u(t)=H(t)u_0+\int^t_0H(t-s)\widetilde{F}(u)(s)\,ds,
\end{align}
where $H(t)x:=e^{t\Delta_{\theta}}x=U_{\theta}(e^{-t|\cdot|^2}f)$ for all $x=U_{\theta}(f)$ with $f\in\cs(\rr^d)$.

\begin{definition}\label{def:mild}
We say that $u$ is a mild solution of \eqref{ns} if $u$ satisfies the integral equality
 \begin{align}\label{sw}
u(t)=H(t)u_0-\int_0^t H(t-s)\mathbb P(F(u))(s)\,ds,\;\forall t\in [0,T], T>0,
\end{align}
where $F(u)=A(u)$ or $S(u)$.
\end{definition}

\begin{lemma}\label{lem:mild}
Let $u\in C([0,T];[L_d(\rr_{\theta}^d)]^d_0)$ with $0<T<\fz$. Then $u$ is a weak solution of \eqref{ns} if and only if $u$ is a mild solution of \eqref{ns}. The definition of $[L_d(\rr_{\theta}^d)]^d_0$ can be found in Section \ref{s6.1}.
\end{lemma}

\begin{proof}
Firstly, we show ``$\Longleftarrow$". Suppose $u$ is a solution to the equation:
$$u(t)=H(t)u_0-\int^t_0H(t-s)\mathbb P \lf(F(u)\r)(s)\,ds.$$ For any $\varphi\in\cs(\rd)$, we have
\begin{align*}
\lim_{h\rightarrow0}\lf\langle \frac{H(t+h)u_0-H(t)u_0}{h},\,\varphi\r\rangle
=&\lim_{h\rightarrow0}\lf\langle u_0,\,\frac{H(t+h)\varphi-H(t)\varphi}{h}\r\rangle\\
=&\lf\langle u_0,\,\Delta_{\theta}H(t)\varphi\r\rangle\\
=&\lf\langle \Delta_{\theta}H(t)u_0,\,\varphi\r\rangle.
\end{align*}
Similarly, for the integral term, we have
\begin{align*}
&\lim_{h\rightarrow0}\lf\langle\frac{\int^{t}_0\lf(H(t+h-s)-H(t-s)\r)
\mathbb{P} \lf(F(u)\r)(s)\,ds}{h},\,
\varphi\r\rangle+\lim_{h\rightarrow0}\lf\langle\frac{\int^{t+h}_tH(t+h-s)
\mathbb{P} \lf(F(u)\r)(s)\,ds}{h},\,
\varphi\r\rangle\\
=&\lf\langle \Delta_{\theta}\int_0^tH(t-s)\mathbb{P} \lf(F(u)\r)(s)\,ds,\,\varphi\r\rangle+\lim_{h\rightarrow0}\lf\langle\frac{\int^{t+h}_tH(t+h-s)
\mathbb{P} \lf(F(u)\r)(s)-\mathbb{P} \lf(F(u)\r)(t)\,ds}{h},\,
\varphi\r\rangle\\&+\big\langle \mathbb{P} \lf(F(u)\r)(t),\,\varphi\big\rangle\\
=&
\lf\langle \Delta_{\theta}\int_0^tH(t-s)\mathbb{P} \lf(F(u)\r)(s)\,ds,\,\varphi\r\rangle
+\big\langle \mathbb{P} \lf(F(u)\r)(t),\,\varphi\big\rangle.
\end{align*}
This implies that
\begin{align*}
\partial_tu&=\Delta_{\theta} H(t)u_0-\Delta_{\theta}\int^t_0H(t-s)\mathbb{P} \lf(F(u)\r)(s)\,ds
-\mathbb{P} \lf(F(u)\r)(t)\\
&=\Delta_{\theta} u-
\mathbb{P} \lf(F(u)\r)(t).
\end{align*}

Next, we prove
``$\Longrightarrow$".
Suppose $u\in C([0,T];[L_d(\rr_{\theta}^d)]^d_0)$ satisfies
$$\partial_tu=\Delta_{\theta} u-\mathbb{P} \lf(F(u)\r).$$
Define
 $$F(u)(t):=H(t)u_0-\int^t_0H(t-s)\mathbb{P} \lf(F(u)\r)(s)\,ds.$$
Then we have
\begin{align*}
\partial_tF(u)(t)=&\Delta_{\theta} H(t)u_0-\mathbb{P} \lf(F(u)\r)(t)
 -\Delta_{\theta}\int^t_0H(t-s)\mathbb{P} \lf(F(u)\r)(s)\,ds\\
=& \Delta_{\theta} F(u)-\mathbb{P} \lf(F(u)\r)(t).
\end{align*}
Therefore,
$$\partial_t(u-F(u))=\Delta_{\theta}(u-F(u)),\ \ \ u(0)-F(u)(0)=0.
$$
Define the equation
\begin{eqnarray}\label{e2.a20}
\lf\{\begin{array}{ll}
\partial_tv=\Delta_{\theta} v;\\[5pt]
v(0)=0.
\end{array}\r.
\end{eqnarray}
Next, we will show that the equation \eqref{e2.a20}  only has the zero solution.
Let $v$ be a solution to \eqref{e2.a20}, and define
$$\widetilde{v}(t):=e^{(T^*-t)\Delta_{\theta}}v(t), \ \ \  0\leq t\leq T^*\leq T.$$
Then, for any $\varphi\in\cs(\rd)$, we have
\begin{align*}
\partial_t\lf\langle\widetilde{v}(t),\,\varphi\r\rangle
=&\lf\langle\partial_t\widetilde{v}(t),\,\varphi\r\rangle\\
=&\lf\langle-\Delta_{\theta}e^{(T^*-t)\Delta_{\theta}}v(t)
+e^{(T^*-t)\Delta_{\theta}}\partial_tv(t),\,\varphi\r\rangle\\
=&\lf\langle\partial_tv(t)-\Delta_{\theta}v(t)
,\,e^{(T^*-t)\Delta_{\theta}}\varphi\r\rangle\\
=&0.
\end{align*}
From the above calculation, we conclude that $\lf\langle\widetilde{v}(t),\,\varphi\r\rangle$ is constant in $t$. Hence, we deduce
$$\widetilde{v}(T^*)=\widetilde{v}(0)=0.$$
Thus, the equation \eqref{e2.a20} only has  the zero solution.
Therefore, we have $u-F(u)=0$, which implies that
$$u=H(t)u_0-\int^t_0H(t-s)\mathbb{P} \lf(F(u)\r)(s)\,ds.$$
\end{proof}

This lemma demonstrates that, to consider the Navier-Stokes equation \eqref{ns1}, we only need to study the integral equation \eqref{sw}.

\subsection{Critical spaces for Navier-Stokes equations on $\rd$}

The critical space is the function space whose norm is invariant under the intrinsic scaling symmetry of a PDE. It represents the precise threshold regularity that determines the transition between ``tame" (subcritical spaces) and ``wild" (supercritical spaces) behavior, making it the central object of study in the analysis of nonlinear PDEs. Problems posed in critical spaces are often the most challenging and interesting, as they are on the verge of ill-posedness and can exhibit phenomena like blow-up (solutions becoming infinite in finite time) or loss of regularity.

We next introduce the critical spaces for Navier-Stokes equations on the quantum Euclidean
space $\rd$. Recall that for $f\in \cs(\mathbb R^d)$, the dilation is defined as $\delta_\varepsilon(f)(s)=f(\varepsilon s)$ for all $s\in\rr^d$ and $\varepsilon>0$. On quantum Euclidean space $\rd$, the dilation is defined via the Weyl transform: for any $\xi\in\rr^d$ and $\varepsilon>0$, define $$\delta_\varepsilon(\lambda_{\theta}(\xi)):=
\varepsilon^{-d}\lambda_{\varepsilon^{-2}\theta}(\varepsilon^{-1}\xi),$$ then for any $x=U_\theta(f)$ with $f\in\cs(\mathbb R^d)$, we have $$\delta_\varepsilon(x)=U_{\varepsilon^{-2}\theta}(\delta_{\varepsilon}(f)),$$ which recovers the usual definition when $\theta=0$.

If $u$ is a solution of equation \eqref{ns1} on $\mathbb R^d_\theta$, so does $u_\varepsilon(t):=\varepsilon\delta_{\varepsilon}(u)(\varepsilon^2t)$ with the initial datum $u_{\varepsilon}(0)=\varepsilon \delta_\varepsilon(u_0)$ on $\mathbb R^d_{\varepsilon^{-2}\theta}$.
Then as usual, one can define the critical spaces in the quantum setting. A space $X=X(\mathbb R^d_\theta)$ is said to be a critical space for the NS equation if the norm $u_{\varepsilon}(0)$ in $X(\mathbb R^d_{\varepsilon^{-2}\theta})$ is invariant for all $\varepsilon>0$.  Similar to the case $\theta=0$, there are corresponding critical spaces for the quantum NS equation.
However, verifying the norm invariance is more intricate. The following result, which was established in \cite{hlw23, MSX20}, is presented here with a proof for the sake of completeness.

\begin{prop}\label{prop:dilation}
Let $1\leq p\leq\infty$. Then for any $x\in L_p(\rd)$, we have
$$\|\delta_\varepsilon(x)\|_{L_p(\rr_{\varepsilon^{-2}\theta}^d)}=\varepsilon^{-\frac dp}\|x\|_{L_p(\rd)}.$$
\end{prop}

\begin{proof}
It suffices to show that, for any $x=U_{\theta}(f)$, $f\in\cs(\rr^d)$,
\begin{align}\label{e2.a1}
\|\delta_\varepsilon(x)\|_{L_p(\rr_{\varepsilon^{-2}\theta}^d)}\leq\varepsilon^{-\frac dp}\|x\|_{L_p(\rd)}.
\end{align}
In fact, the reverse inequality can be obtained via the map $\delta_{\varepsilon^{-1}}$, so the concrete details are omitted.

If $p=2$, for any $x\in\cs(\rd)$, there exist a $f\in\cs(\rr^d)$ such that $x=U_{\theta}(f)$. It follows from Lemma \ref{l2.a1} (i) that
\begin{align*}
\|\delta_\varepsilon(x)\|_{L_2(\rr_{\varepsilon^{-2}\theta}^d)}&
=\lf\|U_{\varepsilon^{-2}\theta}(\delta_{\varepsilon}(f))\r\|
_{L_2(\rr_{\varepsilon^{-2}\theta}^d)}\\
&=\lf\|U_{\varepsilon^{-2}\theta}(f(\varepsilon\cdot))\r\|
_{L_2(\rr_{\varepsilon^{-2}\theta}^d)}\\
&=\lf\|f(\varepsilon\cdot)\r\|_{L_2(\rr^d)}\\
&=\varepsilon^{-\frac{d}{2}}\lf\|x\r\|_{L_2(\rd)}.
\end{align*}
If $p=\fz$, we can check that $\delta_\varepsilon$ is a $*$-isomorphism from
$L_{\fz}(\rd)$ to $L_{\fz}(\rr^d_{\varepsilon^{-2}\theta})$
 (see e.g. \cite{MSX20}). Therefore, we see that $\delta_\varepsilon$ is an isometry. This, in combination with the interpolation theory of $L_p(\rd)$, yields
that \eqref{e2.a1} holds true for $p\in[2,\,\fz]$.

Next we consider the case: $p\in[1,\,2]$. For any $x\in\cs(\rd)$, from the duality theory, the density argument and the H\"{o}lder inequality, we deduce that
\begin{align*}
\|\delta_\varepsilon(x)\|_{L_p(\rr_{\varepsilon^{-2}\theta}^d)}
=&\sup_{y\in\cs(\rr^d_{\varepsilon^{-2}\theta}),\,\|y\|
_{L_{p'}(\rr^d_{\varepsilon^{-2}\theta})}\leq1}
\lf|\tau_{\varepsilon^{-2}\theta}(\delta_\varepsilon(x)y)\r|\\
=&\varepsilon^{-d}\sup_{y\in\cs(\rr^d_{\varepsilon^{-2}\theta}),\,\|y\|_{L_{p'}
(\rr^d_{\varepsilon^{-2}\theta})}\leq1}
\lf|\tau_{\theta}(x\delta_{{\varepsilon}^{-1}}(y))\r|\\
\leq&\varepsilon^{-d}\sup_{y\in\cs(\rr^d_{\varepsilon^{-2}\theta}),\,\|y\|_{L_{p'}
(\rr^d_{\varepsilon^{-2}\theta})}\leq1}
\|x\|_{L_{p}(\rr^d_{\theta})}
\lf\|\delta_{{\varepsilon}^{-1}}(y)\r\|_{L_{p'}(\rr^d_{\theta})}\\
\leq&\varepsilon^{-\frac{d}{p}}\|x\|_{L_{p}(\rr^d_{\theta})}.
\end{align*}
Here we applied the fact that $\tau_{\varepsilon^{-2}\theta}(\delta_\varepsilon(x)y)
=\varepsilon^{-d}\tau_{\theta}(x\delta_{\varepsilon^{-1}}(y))$. Therefore, we complete the proof of Proposition \ref{prop:dilation}.
\end{proof}

\begin{thm}\label{thm:critical}
For any $u_0\in L_d(\rd)$, we have
\begin{align}\label{critical}
\|u_\varepsilon(0)\|_{L_d(\rr^d_{\varepsilon^{-2}\theta})}
=\|u_0\|_{L_d(\rd)},
\end{align}
and thus $L_d(\rd)$ is a critical space of the NS equation  \eqref{ns1}.
\end{thm}
\begin{proof}
For any $r>0$,
$$\|u_\varepsilon(0)\|_{L_r(\rr^d_{\varepsilon^{-2}\theta})}\
=\varepsilon\|\delta_{\varepsilon}(u_0)\|_{L_r(\rr^d_{\varepsilon^{-2}\theta})}
=\varepsilon^{1-\frac{d}{r}}
\|u_0\|_{L_r(\rd)},
$$
which follows from Proposition \ref{prop:dilation}.
\end{proof}

\begin{remark}{\rm
Let us point out that when defining the ``dilation by $\varepsilon>0$" map on $\rd$,  we need to consider it as a mapping between two different quantum Euclidean spaces unless $\theta=0$ so that one cannot say  the norm of $u_{\varepsilon}(0)$ in the critical space $L_d(\rd)$ is invariant for all $\varepsilon>0$. However, a posteriori, we find that this is reasonable due to the following isomorphism (cf. e.g.  \cite{gv98,MSX20}):
$$\rd\cong L_{\fz}(\rr^{\mathrm{dim}(\mathrm{ker}(\theta))})\overline{\otimes} \mathcal{B}(L_2(\rr^{\mathrm{rank}(\theta)/2})),$$
where $\overline{\otimes}$ is the von Neumann tensor product, $\mathrm{ker}(\theta)$ and $\mathrm{rank}(\theta)$ denote the kernel space and rank of matrix $\theta$, respectively.
In particular, if $\det(\theta)\neq0$, then the above isomorphism
reduces to
$$\rd\cong \mathcal{B}(L_2(\rr^{d/2})). $$
Thus all the resulting spaces $\rr^d_{\varepsilon^{-2}\theta}$ for $\varepsilon>0$ are isomorphic to the matrix algebra $\mathcal{B}(L_2(\rr^{d/2}))$.
}
\end{remark}

\section{A transference technique on quantum Euclidean spaces}\label{s4}

In this section, we will establish a transference principle for Fourier multipliers on quantum Euclidean spaces, which will in turn yield the $L_p$-boundedness of the Leray projection \eqref{def of P} in Lemma \ref{t5.1}. Moreover, the approximation process appearing in the proof of the transference principle will deduce the Besov embedding in Lemma \ref{Besov embedding}. Other variants or applications of the transference principle and its proof can be found in \cite{Hon}, \cite[Section 6]{hlw25}. Given a symbol $m:\mathbb R^d\rightarrow\mathbb C$, we abuse the use of the same notation $T_m$ to denote the Fourier multiplier on $\mathbb R^d_\theta$ for all $\theta$.

To state the following result of the present section, we need a notion of ``completely bounded" and the basic properties in operator space theory, which should be well-known to noncommutative analysts and can be found in e.g. \cite[(3.1), p.39]{Pisier}.

\begin{thm}\label{thm:trans}
Let $1<p<\infty$, and $m\in L_\infty(\mathbb R^d)$ be a continuous symbol such that $T_m$ is completely bounded on $L_p(\mathbb R^d)$. Then $T_m$ is completely bounded on $L_q(\mathbb R_\theta^d)$ for any $\min(p,p')\leq q\leq \max(p,p')$. Moreover,
$$\lf\|T: L_q(\mathbb R_\theta^d)\rightarrow L_q(\mathbb R_\theta^d)\r\|\leq \lf\|T_{m}: L_q(\mathbb R^d)\rightarrow L_q(\mathbb R^d)\r\|_{CB}.$$
\end{thm}

The above result still hold in the case $p=1,\infty$, see \cite{Hon}; moreover, in \cite{Hon} the author establishes the transference principle for maximal inequalities and various square function inequalities.

The main ingredient in the proof of Theorem \ref{thm:trans} is the following intertwining identity
\begin{align}\label{e5.21}
\sigma_{\theta}\circ T_m=\lf(T_m \otimes \mathrm{Id}_{\rd}\r)\circ\sigma_{\theta},
\end{align}
where  $\sigma_{\theta}: \rd\rightarrow L_{\fz}(\rr^d)\overline{\otimes}\rd$ is a normal injective $*$-homomorphism  (cf. \cite[Corollary 1.4]{gjp21})  defined as
$$\sigma_{\theta}(\lambda_{\theta}(\xi)):=\exp_{\xi}\otimes \lambda_{\theta}(\xi),$$
where $\exp_{\xi}$ represents the character $x\rightarrow e^{\mathrm{i}(x,\,\xi)}$ in $L_{\fz}(\rr^d)$.

Since $\mathbb{R}^d$ equipped with the Lebesgue measure $dt$ is not a probability space, the $*$-homomorphism $\sigma_\theta$ is not trace-preserving, whence the intertwining identity can not work directly except $p=\infty$. This induces a lot of additional work and an involved approximation argument has to be taken into account. We choose the heat kernels $h_{\epsilon}(\eta)=(\frac{\epsilon}{\pi})^{\frac{d}{2}}e^{-\epsilon|\eta|^2}$ with $\epsilon>0$. Then $\mathbb{R}^d$ equipped with the Gaussian measure $h_{\epsilon}(\eta)d\eta$ is a probability space.

\begin{lemma}[\cite{hlw25}, Lemma 6.2]\label{l5.1}
Let $1\leq p\leq \fz$. Then we have, for any $\epsilon>0, x \in L_p(\rd)$,
$$\|x\|_{L_p(\rd)}=\lf\|h_{\epsilon}^{\frac{1}{p}}\sigma_{\theta}(x)\r\|_{L_p(L_{\fz}(\rr^d)\overline{\otimes} \rd)}.$$
\end{lemma}

\begin{proof}[Proof of Theorem \ref{thm:trans}]
By duality, we may assume $p\geq2$. It suffices to consider the case $q=p$, other cases follow from the interpolation. By Lemma \ref{l5.1}, we obtain that, for any  fixed $x=U_{\theta}(f)$, $f\in\cs(\rr^d)$,
\begin{align*}
\lf\|T_mx\r\|_{L_p(\rd)}
=\lf\|h_{\epsilon}^{\frac{1}{p}}
\sigma_{\theta}(T_mx)\r\|_{L_p(L_{\fz}(\rr^d)\overline{\otimes}\rd)}.
\end{align*}
It follows from the intertwining identity \eqref{e5.21} that
\begin{align*}
h_{\epsilon}^{\frac{1}{p}}(\sigma_{\theta}(T_mx))
=h_{\epsilon}^{\frac{1}{p}}\lf(T_m\otimes\mathrm{Id}_{\rd}\r)\sigma_{\theta}(x).
\end{align*}
From this, we conclude that
\begin{align*}
\lf\|T_mx\r\|_{L_p(\rd)}
=&\lf\|h_{\epsilon}^{\frac{1}{p}}\lf(T_m\otimes\mathrm{Id}_{\rd}\r)\sigma_{\theta}(x)\r\|
_{L_p(L_{\fz}(\rr^d)\overline{\otimes}\rd)}\\
\leq&\lf\|h_{\epsilon}^{\frac{1}{p}}\lf(T_m\otimes\mathrm{Id}_{\rd}\r)\sigma_{\theta}(x)
-\lf(T_m\otimes\mathrm{Id}_{\rd}\r)\lf(h_{\epsilon}^{\frac{1}{p}}\sigma_{\theta}(x)\r)\r\|
_{L_p(L_{\fz}(\rr^d)\overline{\otimes}\rd)}\\
& \ \ \ +\lf\|\lf(T_m\otimes\mathrm{Id}_{\rd}\r)\lf(h_{\epsilon}^{\frac{1}{p}}\sigma_{\theta}(x)\r)\r\|
_{L_p(L_{\fz}(\rr^d)\overline{\otimes}\rd)}
\\
=&:\mathrm{I}_{\epsilon}+\mathrm{J}_{\epsilon}.
\end{align*}
For the term $\mathrm{J}_{\epsilon}$, by the assumption,
we have
\begin{align*}
\mathrm{J}_{\epsilon}&\lesssim
\lf\|T_m\otimes\mathrm{Id}_{\rd}\r\|
_{L_p(L_{\fz}(\rr^d)\overline{\otimes}\rd)\rightarrow L_p(L_{\fz}(\rr^d)\overline{\otimes}\rd)}
\lf\|h_{\epsilon}^{\frac{1}{p}}\sigma_{\theta}(x)\r\|
_{L_p(L_{\fz}(\rr^d)\overline{\otimes}\rd)}\\
&=
\lf\|\mathrm{Id}_{\rd}\otimes T_m\r\|
_{L_p(\rd,L_p(\bR^d))\rightarrow L_p(\rd,L_p(\bR^d))}
\lf\|h_{\epsilon}^{\frac{1}{p}}\sigma_{\theta}(x)\r\|
_{L_p(L_{\fz}(\rr^d)\overline{\otimes}\rd)}\\
&\leq\|T_{m}: L_q(\mathbb R^d)\rightarrow L_q(\mathbb R^d)\|_{CB}
\lf\|x\r\|_{L_p(\rd)},
\end{align*}
where for the first identity we use the noncommutative Fubini theorem (see e.g. \cite[(3.6')]{Pisier}) and for the second inequality, we use the fact that $\rd$ is injective (see e.g. \cite[Section 6]{LSZ20}).
Therefore, to complete the proof of Theorem \ref{thm:trans}, we only need to show that
\begin{align}\label{e5.23}
\lim_{\epsilon\rightarrow0} \,\mathrm{I}_{\epsilon}=0.
\end{align}
By easy calculation, we get
\begin{align}\label{e5.24}
&h_{\epsilon}^{\frac{1}{p}}(t)\lf(T_m\otimes\mathrm{Id}_{\rd}\r)\sigma_{\theta}(x)(t)
-\lf(T_m\otimes\mathrm{Id}_{\rd}\r)\lf(h_{\epsilon}^{\frac{1}{p}}\sigma_{\theta}(x)\r)(t)\\\nonumber
=&\int_{\rr^d}\lf(h_{\epsilon}^{\frac{1}{p}}(t)m(\xi)
-T_{m(\cdot+\xi)}(h_{\epsilon}^{\frac{1}{p}})(t)
\r)\exp_{\xi}(t)f(\xi)\lambda_{\theta}(\xi)\,d\xi.\nonumber
\end{align}
Then by the quantum Hausdorff-Young inequality---Lemma \ref{l2.a1}, the Minkowski inequality and the classical Hausdorff-Young inequality, we obtain
\begin{align}
&\lf\|h_{\epsilon}^{\frac{1}{p}}\lf(T_m\otimes\mathrm{Id}_{\rd}\r)\sigma_{\theta}(x)
-\lf(T_m\otimes\mathrm{Id}_{\rd}\r)\lf(h_{\epsilon}^{\frac{1}{p}}\sigma_{\theta}(x)\r)\r\|^p
_{L_p(L_{\fz}(\rr^d)\overline{\otimes}\rd)}\\\nonumber
=&\tau_\theta\int_{\rr^d}\lf|\int_{\rr^d}\lf(h_{\epsilon}^{\frac{1}{p}}(t)m(\xi)
-T_{m(\cdot+\xi)}(h_{\epsilon}^{\frac{1}{p}})(t)
\r)\exp_{\xi}(t)f(\xi)\lambda_{\theta}(\xi)\,d\xi\r|^p\,dt\\ \nonumber
\leq&\int_{\rr^d}\lf(\int_{\rr^d}\lf|h_{\epsilon}^{\frac{1}{p}}(t)m(\xi)
-T_{m(\cdot+\xi)}(h_{\epsilon}^{\frac{1}{p}})(t)
\r|^{p'}\lf|f(\xi)\r|^{p'}\,d\xi\r)^\frac{p}{p'}\,dt\\ \nonumber
\leq&\lf[\int_{\rr^d}\lf(\int_{\rr^d}\lf|h_{\epsilon}^{\frac{1}{p}}(t)m(\xi)
-T_{m(\cdot+\xi)}(h_{\epsilon}^{\frac{1}{p}})(t)
\r|^{p}\lf|f(\xi)\r|^{p}\,dt\r)^\frac{p'}{p}\,d\xi\r]^{\frac{p}{p'}}\\ \nonumber
=&\lf[\int_{\rr^d}\lf(\int_{\rr^d}\lf|h_{1}^{\frac{1}{p}}(t)m(\xi)
-T_{m(\sqrt{\epsilon}\cdot+\xi)}(h_{1}^{\frac{1}{p}})(t)
\r|^{p}\lf|f(\xi)\r|^{p}\,dt\r)^\frac{p'}{p}\,d\xi\r]^{\frac{p}{p'}}\\ \nonumber
\leq &\lf[\int_{\rr^d}\lf(\int_{\rr^d}\lf|(h_{1}^{\frac{1}{p}})^{\wedge}(\eta)(m(\xi)
-m(\sqrt{\epsilon}\eta+\xi))
\r|^{p'}\,d\eta\r)\,\lf|f(\xi)\r|^{p'}d\xi\r]^{\frac{p}{p'}}.
\end{align}
By the dominated convergence theorem, we conclude assertion \eqref{e5.23}.

\end{proof}

\section{Besov spaces on quantum Euclidean spaces}\label{s5}

In this section, we will introduce Besov spaces on quantum Euclidean spaces, and establish some fundamental properties such as the interpolation and Sobolev embedding theorems, a product estimate and so on. These properties play a crucial role in deducing the sharp time-space estimates in next section, which are instrumental in carrying out the contraction mapping principle to solve the quantum Navier-Stokes equations. More properties on quantum Besov spaces can be found in \cite{Hon}.

As in classical case, Besov spaces on $\rd$ are also defined via the following Littlewood-Paley decomposition.
Let $\phi$ be a smooth radial function on $\rr^d$ such that $\supp \phi\subset B(0,\,2)$ and $\phi=1$ on $B(0,\,1)$. Define
$$\varphi(\xi)=\phi(\xi)-\phi(2\xi).$$
For any $k\in\zz$, set
$$\varphi_k(\xi):=\varphi(2^{-k}\xi),\ \ \ k\in\zz.$$
Then we have $\supp \varphi_k\subset B(0,\,2^{k+1})\backslash B(0,\,2^{k-1})$, $k\in\zz$
 and
$$\sum_{k\in\zz}\varphi_k(\xi)=1,\ \ \ \xi\in\rr^d\backslash\{0\}.$$ We define the homogeneous Littlewood-Paley operator as
$$\dot{\triangle}_kx=\dot{\triangle}_kU_{\theta}(f):=U_{\theta}(\varphi_kf),\ \ \ \mathrm{for \, any}\
x=U_{\theta}(f), f\in\cs(\rr^d), k\in\zz.$$
 The inhomogeneous Littlewood-Paley operator is defined as
$$
\triangle_j:=\lf\{
\begin{array}{cl}
\dot{\triangle}_j,\,\,\, & \mathrm{if} \ j\geq 1,\\
1-\sum_{j\geq 1}\dot{\triangle}_j, \,\,\,& \mathrm{if} \ j =0.
\end{array}\r.
$$

\begin{definition}
Let $\alpha\in\rr$ and $p,\,q\in[1,\,\fz]$. The homogeneous Besov class $\dot{B}^{\alpha}_{p,q}(\rd)$ consists of all distributions $x\in\cs'(\rd)$ such that
$$\|x\|_{\dot{B}^{\alpha}_{
p,q}(\rd)}:=\lf(\sum_{j\in \zz}
2^{j\alpha q}\lf\|\dot{\triangle}_jx\r\|^q_{L_p(\rd)}\r)^{\frac1q}<\fz, \ \ \ \mathrm{for}\ \ q<\fz,$$
and
$$
\|x\|_{\dot{B}^{\alpha}_{p,\fz}(\rd)}:=\sup_{j\in\zz}2^{j\alpha }\lf\|\dot{\triangle}_jx\r\|_{L_{p}(\rd)}.$$
The inhomogeneous Besov class $B^{\alpha}_{p,q}(\rd)$ consists of all distributions $x\in\cs'(\rd)$ such that
$$\|x\|_{B^{\alpha}_{
p,q}(\rd)}:=\lf(\sum_{j\in\nn}
2^{j\alpha q}\lf\|\triangle_jx\r\|^q_{L_p(\rd)}\r)^{\frac1q}<\fz, \ \ \ \mathrm{for}\ \ q<\fz,$$
and
$$
\|x\|_{B^{\alpha}_{p,\fz}(\rd)}:=\sup_{j\in\nn}2^{j\alpha }\lf\|\triangle_jx\r\|_{L_{p}(\rd)}.$$
\end{definition}

The following lemma is noncommutative Bernstein's inequality.
\begin{lemma}\label{lber}
Let $1\leq p\leq q\leq \fz$. Then, for any $f\in \cs(\rr^d)$ with
$\supp f\subset B_r:=\{t\in\rr^d:|t|\leq r\}$, we have
$$\|U_{\theta}(f)\|_{L_q(\rd)}\leq C_dr^{d(\frac{1}{p}-\frac{1}{q})}\|U_{\theta}(f)\|_{L_p(\rd)}.$$
\end{lemma}
\begin{proof}
  First, we prove the case when $q=\infty$ and $p=2^n$ with $n\in \nn.$ When $n=1,$ let $\phi$ be a bump function such that $\phi=1$ on $B_r$ and vanishes outside $B_{2r}.$ Then, by Lemma \ref{l2.a1}, we have
    $$\|U_{\theta}(f)\|_{L_{\infty}(\rd)}=\lf\|U_{\theta}(\phi f)\r\|_{L_{\infty}(\rd)}\le \|\phi f\|_{L_1(\rr^d)}\le \|\phi\|_{L_2(\rr^d)}\|f\|_{L_2(\rr^d)}\le C_dr^{\frac{d}{2}}\|U_{\theta}(f)\|_{L_2(\rd)}.$$
    When $n=2$,  note that $$\lf|U_{\theta}(f)\r|^2=U_{\theta}(f)^*U_{\theta}(f)=U_{\theta}(g),$$
     where for any $\xi\in\rr^d$, $$g(\xi):=\int_{\rr^d}e^{-\frac{\mathrm{i}}{2}(t,\,\theta \xi)}f(\xi-t)\overline{f(t)}\,dt.$$ Hence we have $$\supp g\subset \overline{\supp f-\supp f}\subset B_{2r}$$ and $$\|U_{\theta}(f)\|_{L_{\infty}(\rd)}=\lf\|\lf|U_{\theta}(f)\r|^2\r\|_{L_{\infty}(\rd)}^{\frac 12}\lesssim r^{\frac d 4}\lf\|\lf|U_{\theta}(f)\r|^2\r\|_{L_2(\rd)}^{\frac 12}=r^{\frac d 4}\lf\|U_{\theta}(f)\r\|_{L_4(\rd)}.$$ By induction, we can conclude this for all $n\in \bN.$
Next, we consider the case when $q=\infty$ and $1\leq p<\infty.$ Notice that we can find an integer $n\in \bN$ such that $2^{n-1}<p\le 2^n,$ then
    $$\|U_{\theta}(f)\|_{L_{\infty}(\rd)}^{2^n}\lesssim r^d\tau_{\theta}(|U_{\theta}(f)|^{2^n})\le r^d\tau_{\theta}(|U_{\theta}(f)|^{p})\,\|U_{\theta}(f)\|_{L_{\infty}(\rd)}^{2^n-p},$$ which implies the desired conclusion.

    Finally, for the general case $1\leq p\le q\le \infty,$ by the above conclusion and the H\"{o}lder inequality, we get
    \begin{align*}
        \|U_{\theta}(f)\|_{L_q(\rd)}^q=\tau_{\theta}(|U_{\theta}(f)|^q)\le \|U_{\theta}(f)\|_{L_{\infty}(\rd)}^{q-p}\tau_{\theta}(|U_{\theta}(f)|^p)\lesssim r^{d(q/p-1)}\|U_{\theta}(f)\|_{L_p(\rd)}^q,
    \end{align*}
   which implies the desired inequality.
\end{proof}

\subsection{The interpolation and some embeddings on $\rd$}

In this subsection, we will establish some interpolation and embedding theorems on Besov spaces.
The following lemma outlines some fundamental properties of Besov space $\dot{B}^{\alpha}_{p,q}(\rd).$ We first recall a lemma.

\begin{lemma}\cite[Theorem 6.4.2]{bl76}\label{l2.31}
Let $0<\theta<1, 1\leq q\leq\fz$, $X_0$, $X_1$ and $Y_0$, $Y_1$ be two interpolation couples such that
there exist two operators $\mathbf{S}\in \mathcal{B}(Y_i,\,X_i)$ and $\mathbf{Q}\in \mathcal{B}(X_i,\,Y_i)$ with $\mathbf{S}\mathbf{Q}(f)=f$ for any $f\in X_i$, $i=0,\,1$. Then we have
$$[X_0,\,X_1]_{\theta}=\mathbf{S}[Y_0,\,Y_1]_{\theta},$$
and
$$[X_0,\,X_1]_{\theta,q}=\mathbf{S}[Y_0,\,Y_1]_{\theta,q}.$$
\end{lemma}
\begin{lemma}\label{l2.81}
The Besov spaces have the following elementary properties:
\begin{enumerate}
\item [\rm{(i)}] For $0<\eta<1$, $\alpha_0, \alpha_1\in\rr$ with $\alpha_0\neq\alpha_1$,
$1\leq p,q_0,q_1,q\leq\fz$, we have the following real interpolation and complex
interpolation for $\alpha=(1-\eta)\alpha_0+\eta\alpha_1$,
$$\lf[\dot{B}^{\alpha_0}_{p,q_0}(\rd),\,\dot{B}^{\alpha_1}_{p,q_1}(\rd)
\r]_{\eta,\,q}=\dot{B}^{\alpha}_{p,q}(\rd),$$
and for $1/q=(1-\eta)/{q_0}+{\eta}/{q_1}$
$$\lf[\dot{B}^{\alpha_0}_{p,q_0}(\rd),\,\dot{B}^{\alpha_1}_{p,q_1}(\rd)
\r]_{\eta}=\dot{B}^{\alpha}_{p,q}(\rd).$$

\item [\rm{(ii)}] For $1\leq p\leq p_1<\fz$, $1\leq q\leq q_1\leq\fz$, $\alpha,\alpha_1\in\rr$ satisfy $\alpha-d/p=\alpha_1-d/p_1$, then we have
   $$\dot{B}^{\alpha}_{p,q}(\rd)\subset\dot{B}^{\alpha_1}_{p_1,q_1}(\rd).$$

\item[\rm{(iii)}] For any $1\leq p< p_1\leq\fz$, $1\leq q\leq\fz$ and
$\alpha=d(1/p-1/p_1)$, then we have
$$\dot{B}_{p,q}^{\alpha}(\rd)\subseteq L_{p_1,q}(\rd), \ \ \ \mathrm{if} \ p_1<\fz$$
and
$$\dot{B}_{p,1}^\alpha\subseteq L_{\fz}(\rd), \ \ \ \ \mathrm{if} \ p_1=\fz.$$

\item[\rm{(iv)}] For any $d\geq2$, we have
$$L_d(\rd)\subseteq\dot{B}^0_{d,d}(\rd).$$
\end{enumerate}

\end{lemma}

\begin{proof}
To demonstrate (i), let us consider a fixed Banach space $X$ and denote by $\dot{\ell}^{\alpha}_q(X)$ the weighted $\ell_q$-direct sum of $(X,\,X,\ldots)$, equipped with the norm defined as
$$
\lf\|\{x_j\}_{j\in\zz}\r\|_{\dot{\ell}^{\alpha}_q(X)}:=
 \lf(\sum_{j\in\zz}2^{jq \alpha}\lf\|x_j\r\|_{X}^q\r)^{\frac{1}{q}}.$$
The interpolation properties of the space $\dot{\ell}^{\alpha}_q(X)$ are well-established and can be found in \cite[Theorem 5.6.1]{bl76}. Here, we recall the relevant results. For a given Banach space $X$, $0<\eta<1$, $\alpha_0, \alpha_1\in\rr$ with $\alpha_0\neq\alpha_1$,
$1\leq p,q_0,q_1,q\leq\fz$, we have
\begin{align*}
\lf[\dot{\ell}^{\alpha_0}_{q_0}(X),\,\dot{\ell}^{\alpha_1}_{q_1}(X)
\r]_{\eta,\,q}=\dot{\ell}^{\alpha}_q(X),
\end{align*}
where $\alpha=(1-\eta)\alpha_0+\eta\alpha_1$.

Next, we define the map $\mathbf Q$ as follows
$$\mathbf{Q}:x\mapsto\lf\{\dot{\triangle}_jx\r\}_{j\in\zz}.$$
By its definition, $\mathbf{Q}$ is an isometry from $\dot{B}_{p,q}^{\alpha}(\rd)$ to $\dot{\ell}^{\alpha}_q(L_p(\rd))$. Additionally, we define the map $\mathbf{S}$ on
$\dot{\ell}^{\alpha}_q(L_p(\rd))$ as
$$\mathbf{S}:\{x_i\}_{i\in\zz}\mapsto\sum_{j\in\zz}
\lf(\dot{\triangle}_{j-1}+\dot{\triangle}_j+\dot{\triangle}_{j+1}\r)(x_j).$$
The map $\mathbf{S}$ is bounded from  $\dot{\ell}^{\alpha}_q(L_p(\rd))$ to
$\dot{B}_{p,q}^{\alpha}(\rd)$ and
satisfies $\mathbf{S}\mathbf{Q}(x)=x$, for any $x\in\dot{B}_{p,q}^{\alpha}(\rd)$.
From this and Lemma \ref{l2.31}, we deduce that
$$\lf[\dot{B}^{\alpha_0}_{p,q_0}(\rd),\,\dot{B}^{\alpha_1}_{p,q_1}(\rd)
\r]_{\eta,\,q}
=\mathbf{S}\lf[\dot{\ell}^{\alpha_0}_{q_0}(L_p(\rd)),\,\dot{\ell}^{\alpha_1}_{q_1}(L_p(\rd))
\r]_{\eta,\,q}=\mathbf{S}(\dot{\ell}^{\alpha}_q(L_p(\rd)))=\dot{B}^{\alpha}_{p,q}(\rd).$$
 The  complex interpolation conclusion can be proved by a slight modification
 of the above process, and the details are omitted for brevity.

For (ii), since $\supp \varphi_j\subset\{t\in\rr^d:|t|\leq2^j\} $, for any $1\leq p\leq p_1\leq\fz$, by Lemma \ref{lber}, we obtain
\begin{align*}
\lf\|\dot{\triangle}_jx\r\|_{L_{p_1}(\rd)}
\leq C_d2^{dj(\frac{1}{p}-\frac{1}{p_1})}\lf\|x\r\|_{L_{p}(\rd)}.
\end{align*}
Therefore, it follows from this and $\alpha_1-d/p=\alpha-d/p_1$ that
\begin{align*}
2^{\alpha_1j}\lf\|\dot{\triangle}_jx\r\|_{L_{p_1}(\rd)}
\lesssim2^{\alpha j}\lf\|\dot{\triangle}_jx\r\|_{L_{p}(\rd)}.
\end{align*}
Taking the $\ell_q$-norm in both sides of the above inequality, we have
\begin{align*}
\lf\|x\r\|_{\dot{B}^{\alpha_1}_{p_1,q}(\rd)}
\lesssim\lf\|x\r\|_{\dot{B}^{\alpha}_{p,q}(\rd)},
\end{align*}
which, combined with $\ell_q\subset\ell_{q_1}$, implies that
$$\dot{B}^{\alpha}_{p,q}(\rd)\subset\dot{B}^{\alpha_1}_{p_1,q_1}(\rd).$$

To prove (iii),
we apply (ii) with $\alpha_1= 0$ and $q=q_1=1$. This yields the embeddings:
$$\dot{B}^{\alpha}_{p,1}(\rd)\subset \dot{B}^0_{p_1,1}(\rd)
\subset L_{p_1}(\rd).$$
When $p_1=\fz$, this directly implies
$\dot{B}^{d/p}_{p,1}(\rd)
\subset L_{\fz}(\rd).$
 For the case where $p_1<\fz$, we combine the above embedding with Lemmas \ref{lor1} and (ii) to conclude
 $$\lf[\dot{B}^{0}_{p,1}(\rd),\,\dot{B}^{d/p}_{p,1}(\rd)
\r]_{\eta,\,q}=\dot{B}^{\alpha}_{p,q}(\rd) \ \  \ \mathrm{and} \ \ \
\lf[L_p(\rd),\,L_{\fz}(\rd)
\r]_{\eta,\,q}=L_{p_1,q}(\rd),$$
and thus
$\dot{B}^{\alpha}_{p,q}(\rd) \subset L_{p_1,q}(\rd).$

Finally, we prove (iv).
Notice that for any $x\in L_2(\rd)$,
$$\lf(\sum_{j\in\zz}\lf\|\dot{\triangle}_jx\r\|_{L_2(\rd)}^2\r)^{\frac{1}{2}}
\lesssim\lf\|x\r\|_{L_2(\rd)}.$$
Moreover, for any $x\in L_{\fz}(\rd)$, it follows from Lemma \ref{lem:psi} that
$$\sup_{j\in\zz}\lf\|\dot{\triangle}_{j}x\r\|_{L_{\fz}(\rd)}
\leq\sup_{j\in\zz}\lf\|\mathcal{F}^{-1}\varphi_j\r\|_{L_1(\rr^d)}\lf\|x\r\|_{L_{\fz}(\rd)}
\lesssim\lf\|x\r\|_{L_{\fz}(\rd)}.$$
Using these results and the interpolation theorem from \cite[Theorem 5.6.2]{bl76}, which states:
$$\lf[\dot{\ell}_{q_0}^{\alpha_0}(X_0),\,\dot{\ell}_{q_1}^{\alpha_1}(X_1)\r]_{\eta}
=\dot{\ell}_{q}^{\alpha}([X_0,\,X_1]_{\eta})$$ with $\alpha=(1-\eta)\alpha_0+\eta\alpha_1$  and $1/q=(1-\eta)/q_0+\eta/q_1,$
we conclude that for any $d\geq2$,
$$\lf(\sum_{j\in\zz}\lf\|\dot{\triangle}_jx\r\|_{L_d(\rd)}^d\r)^{\frac{1}{d}}
\lesssim\lf\|x\r\|_{L_d(\rd)}.$$
This implies the embedding stated in (iv).
\end{proof}

The following is the so-called reduction theorem of the Besov space $\dot{B}^{\alpha}_{p,q}(\rd)$, whose proof is similar to the quantum tori case (cf. \cite[Theorem 3.7(ii)]{XXY18}), and we omit the proof.
\begin{lemma}\label{bd}
 Let $a\in\zz^d_+$, $\alpha\in\rr$ and
$1\leq p,q\leq\fz$. If $x\in\dot{B}^{\alpha}_{p,q}(\rd)$, then $\partial^ax\in \dot{B}^{\alpha-|a|}_{p,q}(\rd)$ and there exists a positive constant $C$ such that
 $$\|\partial^ax\|_{\dot{B}^{\alpha-|a|}_{p,q}(\rd)}\leq C\|x\|_{\dot{B}^{\alpha}_{p,q}(\rd)}.$$
\end{lemma}

The following lemma is a Sobolev embedding theorem on $\rd$ which can be deduced
from the noncommutative Young inequality and Bernstein's inequality.
\begin{lemma}\label{se1}
    For any $\alpha>\frac dp-\frac dq,1\le p<q\le \fz,$ we have the Sobolev embedding $H_p^{\alpha}(\rd)\subset L_q(\rd).$
\end{lemma}
\begin{proof}
For any $x\in H_p^{\alpha}(\rd)$, from Lemmas \ref{lber} and \ref{lem:psi}, one concludes that
\begin{align*}
\|x\|_{L_q(\rd)}&\leq \sum_{j\in\nn}\lf\|\triangle_jx\r\|_{L_q(\rd)}\
\lesssim\sum_{j\in\nn}2^{jd(\frac{1}{p}-\frac{1}{q})}
\lf\|\triangle_jx\r\|_{L_p(\rd)}\\
&\lesssim\sum_{j\in\nn}2^{jd(\frac{1}{p}-\frac{1}{q})}2^{-j\alpha}
\lf\|x\r\|_{H_p^{\alpha}(\rd)}
\thicksim\lf\|x\r\|_{H_p^{\alpha}(\rd)}.
\end{align*}
\end{proof}

\begin{lemma}\label{bec}
Let $\alpha\in\rr$, $\varepsilon>0$ and $1\leq p\leq \fz$. Then we have
$B^{\alpha+\varepsilon}_{p,\infty}(\rd)\subseteq H_p^{\alpha}(\rd).$
\end{lemma}
\begin{proof}
For any $x\in B_{p,\fz}^{\alpha+\varepsilon}(\rd)$, by the definitions and applying Lemma \ref{lem:psi}, we have
\begin{align*}
\|x\|_{H_p^{\alpha}(\rd)}&\leq \sum_{j\in\nn}\lf\|\triangle_j(1-\Delta_{\theta})^{\frac{\alpha}{2}}x\r\|_{L_p(\rd)}\\
&\lesssim\sum_{j\in\nn}2^{j\alpha}
\lf\|\triangle_jx\r\|_{L_p(\rd)}\\
&\lesssim\sum_{j\in\nn}2^{-j\varepsilon}\sup_{j\in\nn}2^{j(\alpha+\varepsilon)}
\lf\|\triangle_jx\r\|_{L_p(\rd)}\\
&\thicksim\lf\|x\r\|_{B_{p,\fz}^{\alpha+\varepsilon}(\rd)}.
\end{align*}
\end{proof}

\begin{lemma}\label{hom and inhom}
    For $\alpha>0,1\le p,q\le \infty,$ we have
    \begin{align*}
        B_{p,q}^{\alpha}(\rd)=\dot{B}_{p,q}^{\alpha}(\rd)\cap L_p(\rd).
    \end{align*}
\end{lemma}
\begin{proof}
Let $\alpha>0,1\le p,q\le \infty.$ By the inequality $$\|x\|_{B_{p,q}^\alpha(\rd)}^q\le \|x\|_{L_p(\rd)}^q+\sum_{j\ge 1}2^{j\alpha q}\lf\|\dot\triangle_j x\r\|_{L_p(\rd)}\le \|x\|_{L_p(\rd)}^q+\|x\|_{\dot B_{p,q}^\alpha(\rd)}^q,$$ we have
 $$\dot{B}_{p,q}^{\alpha}(\rd)\cap L_p(\rd)\subset B_{p,q}^{\alpha}(\rd).$$
Next, we prove that
$$ B_{p,q}^{\alpha}(\rd)\subset\dot{B}_{p,q}^{\alpha}(\rd)\cap L_p(\rd).$$
For any $x\in B_{p,q}^{\alpha}(\rd)$, we have
\begin{align*}
\|x\|_{L_p(\rd)}\leq\sum_{j\in\nn}\lf\|\triangle_jx\r\|_{L_p(\rd)}
\leq\lf(\sum_{j\in\nn}2^{-\alpha q}\r)^{\frac{1}{q}}
\lf(\sum_{j\in\nn}2^{\alpha q}\lf\|\triangle_jx\r\|_{L_p(\rd)}^q\r)^{\frac{1}{q}}
\thicksim\|x\|_{B_{p,q}^{\alpha}(\rd)}.
\end{align*}
From this result and Lemma \ref{lem:psi}, we conclude that
\begin{align*}
\|x\|_{\dot{B}_{p,q}^{\alpha}(\rd)}^q
=&\sum_{j=-\fz}^{0}2^{j\alpha q}
\lf\|\dot{\triangle}_jx\r\|_{L_p(\rd)}^q
+\sum_{j=1}^{\fz}2^{j\alpha q}
\lf\|\triangle_jx\r\|_{L_p(\rd)}^q\\
\lesssim&\sum_{j=-\fz}^{0}2^{j\alpha q}
\lf\|x\r\|_{L_p(\rd)}^q
+\sum_{j=1}^{\fz}2^{j\alpha q}
\lf\|\triangle_jx\r\|_{L_p(\rd)}^q\\
\lesssim&\lf\|x\r\|_{L_p(\rd)}^q+\|x\|_{B_{p,q}^{\alpha}(\rd)}\\
\lesssim&\|x\|_{B_{p,q}^{\alpha}(\rd)},
\end{align*}
which implies the desired conclusion.
\end{proof}

At the end of this subsection, we will establish an embedding theorem about Besov spaces and noncommutative $L_p$-spaces on $\rd$.
To prove our desired conclusion, we need to recall the definition of the Hilbert space-valued noncommutative $L_p$-spaces. In what follows,
we will only introduce the following concrete representations of these spaces; for a more
general description we refer to the papers \cite{jlx06,lpp91}.

For $1\leq p\leq\infty$, consider a finite sequence $\{x_n\}_n$ in $L_p(\mathcal{M})$. We define the norms as follows:
 \begin{equation*}
   \|\{x_n\}_n\|_{L_p(\mathcal{M};\ell_2^c)}:
   =\lf\|\lf(\sum_n|x_n|^2\r)^{\frac12}\r\|_{L_p(\mathcal{M})},\quad  \|\{x_n\}_n\|_{L_p(\mathcal{M};\ell_2^r)}:=\lf\|\{x_n^*\}_n\r\|_{L_p(\mathcal{M};\ell_2^c)}.
 \end{equation*}
 For $1\leq p<\infty$, the space $L_p(\mathcal{M};\ell_2^c)$ (respectively $L_p(\mathcal{M};\ell_2^r ) )$ is defined as the completion of the set of all finite sequences in $L_p(\mathcal{M})$ with respect to $\|\cdot\|_{L_p(\mathcal{M};\ell_2^c)}$ (respectively $\|\cdot\|_{L_p(\mathcal{M};\ell_2^r)} )$. For $p=\infty$, the Banach space $L_\infty(\mathcal{M};\ell_2^c)$ (respectively $L_\infty(\mathcal{M};\ell_2^r ))$ consists of all sequences in $L_\infty(\mathcal{M})$ such that the series $\sum_n x_n^*x_n$ (respectively $\sum_n x_nx_n^*$) converges in the weak-$*$ topology.

 \begin{lemma}\label{column}
Let $2\leq p\leq\infty$ and $\{x_n\}_n\subset L_p(\mathcal{M})$, we have
\begin{equation*}
  \max\lf\{\|\{x_n\}_n\|_{L_p(\mathcal{M};\ell_2^c)},\|\{x_n\}_n
  \|_{L_p(\mathcal{M};\ell_2^r)}\r\}\leq\lf(\sum_n\|x_n\|_{L_p(\mathcal{M})}^2\r)^{\frac12}.
\end{equation*}
\end{lemma}
\begin{proof}
  By the triangle inequality, we deduce
  \begin{equation*}
    \|\{x_n\}_n\|_{L_p(\mathcal{M};\ell_2^c)}^2=
    \lf\|\sum_n|x_n|^2\r\|_{L_{p/2}(\cm)}\leq\sum_n\lf\||x_n|^2\r\|_{L_{p/2}(\cm)}
    =\sum_n\|x_n\|_{L_p(\cm)}^2
  \end{equation*}
  and
   \begin{equation*}
    \lf\|\{x_n\}_n\r\|_{L_p(\mathcal{M};\ell_2^r)}^2
    =\lf\|\sum_n\lf|x_n^*\r|^2\r\|_{L_{p/2}(\cm)}
    \leq\sum_n\lf\|\lf|x_n^*\r|^2\r\|_{L_{p/2}(\cm)}=\sum_n\|x_n\|_{L_p(\cm)}^2.
  \end{equation*}
\end{proof}
In what follows, the above lemma will be frequently used by choosing
$\cm=\rd$.
Now, we can establish the following embedding property:
\begin{lemma}\label{Besov embedding}
    Let $p\in[2,\,\fz)$. Then we have
    \begin{align*}
        \dot B_{p,2}^0(\rd)\subset L_p(\rd).
    \end{align*}
\end{lemma}

\begin{proof}
Given $x\in \dot B_{p,2}^0(\rd)$, by the density argument, we may assume that $x=U_{\theta}(f)$, $f\in\cs(\rr^d)$. By Lemma \ref{column}, we have
$$\max \lf\{\lf\|\lf\{\dot{\triangle}_jx\r\}_j\r\|
_{L_p(\rd;\ell_2^c)},\lf\|\lf\{\dot{\triangle}_jx\r\}_j\r\|_{L_p(\rd;\ell_2^r)}\r\}
\le\|x\|_{\dot B_{p,2}^0(\rd)}.$$
Hence, it suffices to show that
    \begin{align}\label{Besov embedding:eq1}
        \max\lf\{\lf\|\lf\{\dot{\triangle}_jx\r\}_j\r\|
        _{L_p(\rd;\ell_2^c)},\lf\|\lf\{\dot{\triangle}_jx\r\}_j\r\|
        _{L_p(\rd;\ell_2^r)}\r\}\sim\|x\|_{L_p(\rd)}.
    \end{align}
By Lemma \ref{l5.1} and the operator-valued Littlewood-Paley theorem (see \cite[Section 2.4.2]{mp2009}), we  get
\begin{align*}
  \|x\|_{L_p(\rd)}&=\lf\|h_\epsilon^{\frac 1p}\sigma_{\theta}(x)\r\|_{L_p(L_{\fz}(\rr^d)\overline{\otimes}\rd)}\\
  &\thicksim \max \lf\{\lf\|\lf\{\dot{\triangle}_j(h_\epsilon^{\frac 1p}\sigma_{\theta}(x))\r\}_j\r\|_{L_p(L_{\fz}(\rr^d)\overline{\otimes}\rd;\ell_2^c)},
  \lf\|\lf\{\dot{\triangle}_j(h_\epsilon^{\frac 1p}\sigma_{\theta}(x))\r\}_j\r\|_{L_p(L_{\fz}(\rr^d)\overline{\otimes}\rd;\ell_2^r)}\r\}.
\end{align*}
Thus, it suffices to show that
\begin{align}\label{Besov embedding:eq2}
        \lf\|\lf\{\dot{\triangle}_jx\r\}_j\r\|_{L_p(\rd;\ell_2^c)}
        =\lim_{\epsilon\to 0}\lf\|\lf\{\dot{\triangle}_j(h_\epsilon^{\frac 1p}\sigma_{\theta}(x))\r\}_j\r\|_{L_p(L_{\fz}(\rr^d)\overline{\otimes}\rd;\ell_2^c)}
    \end{align}
    and
    \begin{align}\label{Besov embedding:eq3}
      \lf\|\lf\{\dot{\triangle}_jx\r\}_j\r\|_{L_p(\rd;\ell_2^r)}=\lim_{\epsilon\to 0}\lf\|\lf\{\dot{\triangle}_j(h_\epsilon^{\frac 1p}\sigma_{\theta}(x))\r\}_j\r\|_{L_p(L_{\fz}(\rr^d)\overline{\otimes}\rd;\ell_2^r)}.
    \end{align}

By Lemma \ref{l5.1}, we have
\begin{equation*}
  \lf\|\lf\{\dot{\triangle}_jx\r\}_j\r\|_{L_p(\rd;\ell_2^c)}
  =\lf\|h_\epsilon^{\frac 1p}\sigma_{\theta}\lf(\sum_{j\in\mathbb{Z}}\lf|\dot{\triangle}_jx\r|^2\r)^{\frac12}
  \r\|_{L_p(L_{\fz}(\rr^d)\overline{\otimes}\rd)}
  =\lf\|\lf\{h_\epsilon^{\frac 1p}\sigma_\theta(\dot{\triangle}_jx)\r\}_j\r\|_{L_p(L_{\fz}(\rr^d)
  \overline{\otimes}\rd;\ell_2^c)}.
\end{equation*}
Then, by the triangle inequality and Lemma \ref{column}, we deduce
\begin{align*}
  &\lf|\lf\|\lf\{\dot{\triangle}_j(h_\epsilon^{\frac1p}\sigma_{\theta}(x))\r\}_j\r\|
  _{L_p(L_{\fz}(\rr^d)\overline{\otimes}\rd;\ell_2^c)}-
   \lf\|\lf\{\dot{\triangle}_jx\r\}_j\r\|_{L_p(\rd;\ell_2^c)}\r|^2\\
   \leq& \lf\|\lf\{\dot{\triangle}_j( h_\epsilon^{\frac 1p}(\sigma_{\theta}(x)))-h_\epsilon^{\frac 1p}\sigma_{\theta}(\dot{\triangle}_jx)\r\}_{j\in \bZ}\r\|_{L_p(L_{\fz}(\rr^d)\overline{\otimes}\rd;\,\ell_2^c)}^2\\
   \leq&\sum_{j\in\mathbb{Z}}\lf\|\dot{\triangle}_j (h_\epsilon^{\frac 1p}(\sigma_{\theta}(x)))-h_\epsilon^{\frac 1p}\sigma_{\theta}(\dot{\triangle}_jx)\r\|_{L_p(L_{\fz}(\rr^d)\overline{\otimes}\rd)}^2=:\mathrm{K}_{\epsilon}.
\end{align*}
Analogous to the proof of Theorem \ref{thm:trans}, we derive
\begin{align*}
  \mathrm{K}_{\epsilon}^{\frac{p^\prime}{2}}&\leq\lf[\sum_{j\in\mathbb{Z}}\lf(\int_{\rr^d}\int_{\rr^d}\lf|(h_{1}^{\frac{1}{p}})^{\wedge}(\eta)(\varphi(2^{-j}\xi)
-\varphi(2^{-j}(\sqrt{\epsilon}\eta+\xi)))
\r|^{p'}\lf|f(\xi)\r|^{p'}\,d\eta d\xi\r)^{\frac{2}{p'}}\r]^{\frac{p^\prime}{2}}\\
&\leq\sum_{j\in\mathbb{Z}}\int_{\rr^d}\int_{\rr^d}\lf|(h_{1}^{\frac{1}{p}})^{\wedge}(\eta)(\varphi(2^{-j}\xi)
-\varphi(2^{-j}(\sqrt{\epsilon}\eta+\xi)))
\r|^{p'}\lf|f(\xi)\r|^{p'}\,d\eta d\xi\\
&=\int_{\rr^d}\int_{\rr^d}\sum_{j\in\mathbb{Z}}\lf|(h_{1}^{\frac{1}{p}})^{\wedge}(\eta)(\varphi(2^{-j}\xi)
-\varphi(2^{-j}(\sqrt{\epsilon}\eta+\xi)))
\r|^{p'}\lf|f(\xi)\r|^{p'}\,d\eta d\xi.
\end{align*}
Then, $\lim_{\epsilon\rightarrow0}\mathrm{K}_{\epsilon}=0$ by the dominated convergence theorem. Therefore, we obtain the equality \eqref{Besov embedding:eq2}. The proof of \eqref{Besov embedding:eq3} is similar, so we omit the details here.
\end{proof}

\subsection{A product estimate}
In this subsection, we provide a product estimate of the intersection of the same quantum Besov space with $L_\infty(\rd)$ (see also \cite[Corollary 5.4]{m23}), which is quite useful to deal with the nonlinear term in solving PDEs. A more complete picture of the product estimates between Besov spaces, Sobolev spaces and more generally Trieble-Lizorkin spaces on quantum Euclidean spaces can be found in the second author's note \cite{Hon}.

\begin{lemma}\label{McDonald}
Let $s>0$ and $p, q\in[1,\,\fz]$. Then, for any $u, v\in B_{p,q}^s(\rd)\cap L_{\fz}(\rd)$, we have
$uv\in B_{p,q}^s(\rd)\cap L_{\fz}(\rd)$ and $$\|uv\|_{B_{p,q}^s(\rd)}
\lesssim_{s,p,q}\|u\|_{B_{p,q}^s(\rd)}\|v\|_{L_{\fz}(\rd)}+\|u\|_{L_{\fz}(\rd)}
\|v\|_{B_{p,q}^s(\rd)}.$$
\end{lemma}
\begin{proof}
  Let $S_j=\sum_{k=0}^j\triangle_k.$ Recall for the paraproduct $uv,$ we have the Bony decomposition (see e.g \cite[Section 4.4.1]{RS})
    \begin{align*}
        uv=\Pi_1(u,v)+\Pi_2(u,v)+\Pi_3(u,v),
    \end{align*}
    where \begin{align*}
        &\Pi_1(u,v):=\sum_{j\ge 3}S_{j-3}u\triangle_j v=:\Pi_2(v,u),\\
        &\Pi_3(u,v):=\sum_{k=-1}^1\sum_{j\in \bN}\triangle_{j+k}u\triangle_jv.
    \end{align*}
    For $\Pi_1(u,v),$ by using the fact that $$\triangle_n(\triangle_ku\triangle_jv)=0, \ \ \ n> \max\{j,k\}+3,$$ we can obtain
    \begin{align*}
\|\Pi_1(u,v)\|_{B_{p,q}^s(\rd)}^q&=\sum_{n\in\bN}2^{snq}\lf\|\sum_{j\ge n-3}\triangle_n(S_{j-3}u\triangle_jv)\r\|_{L_p(\rd)}^q\\
&\lesssim \sum_{n\in\bN}\sum_{j\ge n-3}2^{snq}\lf\|S_{j-3}u\r\|_{L_{\fz}(\rd)}^q\lf\|\triangle_jv\r\|_{L_p(\rd)}^q\\
&\lesssim \|u\|_{L_{\fz}(\rd)}^q\sum_{j\in\bN}\sum_{n\le j+3}2^{snq}\lf\|\triangle_jv\r\|_{L_p(\rd)}^q\\
&\lesssim \|u\|_{L_{\fz}(\rd)}^q\sum_{j\in\bN}2^{sjq}\lf\|\triangle_jv\r\|_{L_p(\rd)}^q\\
&=\|u\|_{L_{\fz}(\rd)}^q\|v\|_{B_{p,q}^s(\rd)}^q.
    \end{align*}
Therefore, we have
 \begin{align*}
    &\|\Pi_1(u,v)\|_{B_{p,q}^s(\rd)}\lesssim \|u\|_{B_{p,q}^s(\rd)}\|v\|_{L_{\fz}(\rd)},\\
    &\|\Pi_2(u,v)\|_{B_{p,q}^s(\rd)}\lesssim \|v\|_{L_{\fz}(\rd)}\|u\|_{B_{p,q}^s(\rd)}.
\end{align*}
For $\Pi_3(u,v),$ by a similar proof, we have
\begin{align*}
\lf\|\Pi_3(u,v)\r\|_{B_{p,q}^s(\rd)}^q&=\sum_{n\in\bN}2^{snq}\lf\|\sum_{k=-1}^1\sum_{j\ge n-3}\triangle_n(\triangle_{j+k}u\triangle_jv)\r\|_{L_p(\rd)}^q\\
&\lesssim \max_{|v|\le 1}\sum_{n\in\bN}2^{snq}\lf\|\sum_{j\ge n-3}\triangle_n(\triangle_{j+k}u\triangle_jv)\r\|_{L_p(\rd)}^q\\
&\lesssim \|u\|_{L_{\fz}(\rd)}^q\sum_{j\in\bN}\sum_{n\le j+3}2^{snq}\lf\|\triangle_jv\r\|_{L_p(\rd)}^q\\
&=\|u\|_{L_{\fz}(\rd)}^q\|v\|_{B_{p,q}^s(\rd)}^q,
\end{align*}
which implies the desired conclusion.
\end{proof}

\section{The sharp time-space estimates}\label{s6.1}

In this section, in order to prepare for the proof of main results, as in the classical case (see e.g. \cite{k84,whhg11}), we focus intently on proving the estimates for the quantum heat semigroup, including the sharp time-space estimates for both the linear terms and the non-linear terms, among others.

Hereafter, for the sake of simplicity, let $X(\rd)$ denote the Banach function spaces on quantum Euclidean space $\rd$, such as $L_p$-spaces $L_p(\rd)$, Besov spaces $B_{p,q}^{\alpha}(\rd)$, Sobolev spaces $H_p^k(\rd)$ and so on. For any vector $u=(u_1,\ldots,u_d)\in [X(\rd)]^d$, we define the norms as follows:
$$\|u\|_{[X(\rd)]^d}=\lf(\sum_{j=1}^d\lf\|u_j\r\|^2_{X(\rd)}\r)^{\frac{1}{2}} \ \ \mathrm{and} \ \ \|\nabla_{\theta} u\|_{[X(\rd)]^{d^2}}=\lf(\sum_{j,\,k=1}^d\lf\|\partial_ju_k\r\|^2_{X(\rd)}\r)^{\frac{1}{2}}.$$
For vector-valued function $u(\cdot)=(u_1(\cdot),\ldots,u_d(\cdot))$ defined on the interval $I\subset\mathbb R$ with values in $[X(\rd)]^d$, we define
\begin{align}\label{lqx}
\|u\|_{L_q(I;[X(\rd)]^d)}=\lf(\int_I\|u(t)\|^q_{[X(\rd)]^d}\,dt\r)^{\frac1q}.
\end{align}
We then introduce the following Banach spaces:
$$\lf[X(\rd)\r]^d_0=\lf\{u\in \lf[X(\rd)\r]^d:{\bf div}\,u=0\r\}$$ and
$$L_q(I;[X(\rd)]^d_0)=\lf\{u\in L_q(I;[X(\rd)]^d):{\bf div}\,u(t)=0,\,\mathrm{a.e.}\,t\in I\r\}.$$

In what follows, for convenience, sometimes we write $\|\nabla_{\theta}x\|_{[X(\rd)]^d}$ as $\|\nabla_{\theta}x\|_{X(\rd)}.$

\subsection{The $(L_r, L_p)$-type estimate for the heat semigroup}

In this subsection, we are devoted to establishing the ($L_r$, $L_p$) estimate for the heat semigroup,
where $1\leq r\leq p\leq\infty$.

\begin{prop}\label{prop:rp}
Let $1\leq r\leq p\leq\infty$. Then, for $k=0,1$, $t>0$ and $x\in L_r(\rd)$, we have the following estimates:
\begin{align}\label{rp}
\lf\|\nabla_{\theta}^kH(t)x\r\|_{L_p(\rd)}\lesssim t^{-\frac k2-\frac d2(\frac 1r-\frac 1p)}\|x\|_{L_r(\rd)}.
\end{align}
\end{prop}

\begin{proof} By the density argument, we may assume that $x=U_{\theta}(f)$, $f\in\cs(\rr^d)$.
Now we consider Proposition \ref{prop:rp} by two cases: $k=0$ and $k=1$.

\textbf{Case 1.} If $k=0$, by applying Lemma \ref{lem:psi}, we obtain
$$\lf\|H(t)x\r\|_{L_r(\rd)}=\lf\|U_{\theta}\lf(e^{-t|\cdot|^2}f\r)\r
\|_{L_r(\rd)}
\leq\lf\|\mathcal{F}^{-1}e^{-t|\cdot|^2}\r\|_{L_1(\rr^d)}
\|x\|_{L_r(\rd)}\leq\|x\|_{L_r(\rd)}.$$
Using Lemma \ref{lem:psi} again, we have
\begin{align*}
\lf\|H(t)x\r\|_{L_{\fz}(\rd)}&=\lf\|U_{\theta}\lf(e^{-t|\cdot|^2}
f\r)\r
\|_{L_{\fz}(\rd)}\\
&\leq\lf\|\mathcal{F}^{-1}e^{-t|\cdot|^2}\r\|_{L_{r'}(\rr^d)}
\|x\|_{L_r(\rd)}\\
&\lesssim t^{-\frac {d}{2r}}\|x\|_{L_r(\rd)}.
\end{align*}
 Combining these estimates with the noncommutative H\"{o}lder inequality, we have
$$\lf\|H(t)x\r\|_{L_{p}(\rd)}\leq\lf\|H(t)x\r\|_{L_{\fz}(\rd)}^{1-\frac rp}
\lf\|H(t)x\r\|_{L_{r}(\rd)}^{\frac rp}
\lesssim t^{-\frac d2(\frac1r-\frac1p)}\|x\|_{L_r(\rd)}.
$$

\textbf{Case 2.} If $k=1$, by Lemma \ref{lem:psi}, we have
\begin{align*}
\lf\|\nabla_{\theta} H(t)x\r\|_{L_p(\rd)}
&=\lf(\sum_{j=1}^d\lf\|\partial_j H(t)x\r\|_{L_p(\rd)}^2\r)^{\frac12}\\
&=\lf(\sum_{j=1}^d\lf\|U_{\theta}(\xi_je^{-t|\xi|^2}f)\r\|
_{L_p(\rd)}^2\r)^{\frac12}\\
&\leq\lf(\sum_{j=1}^d\lf\|\mathcal{F}^{-1}(\xi_je^{-t|\xi|^2})\r\|
_{L_{1}(\rr^d)}^2\r)^{\frac 12}
\|x\|_{L_p(\rd)}\\
&\lesssim t^{-\frac12}\|x\|_{L_p(\rd)}.
\end{align*}
From this and \textbf{Case 1}, we further conclude that, for any $t>0$,

\begin{align*}
\lf\|\nabla_{\theta} H(t)x\r\|_{L_p(\rd)}
=\lf\|H(\frac t2)\nabla_{\theta} H(\frac t2)x\r\|_{L_p(\rd)}
\lesssim t^{-\frac{d}2(\frac1r-\frac1p)}\lf\|\nabla_{\theta} H(\frac t2)x\r\|_{L_r(\rd)}
\lesssim t^{-\frac12-\frac{d}2(\frac1r-\frac1p)}\|x\|_{L_r(\rd)}.
\end{align*}

\end{proof}

The following lemma is the noncommutative Gagliardo-Nirenberg inequality.
\begin{lemma}\label{gn1}
    Let $1<p,q<\infty$ and $0<\vartheta<1$ satisfy $\frac 1q=\frac 1p-\frac{\vartheta}{d}$, Then for any $u\in W_p^1(\rd),$
    \begin{align*}
        \|u\|_{L_q(\rd)}\lesssim \|u\|_{L_p(\rd)}^{1-\vartheta}\|\nabla_{\theta} u\|_{L_p(\rd)}^{\vartheta}.
    \end{align*}
\end{lemma}
\begin{proof}
Assume $u=U_{\theta}(f), f\in \cs(\rr^d).$ When $\|u\|_{L_p(\rd)}=0$, the inequality holds trivially. When $\|\nabla_{\theta} u\|_{L_p(\rd)}=0,
\nabla_{\theta} u=\{U_{\theta}(\mathrm{i}\mathbf D_j f)\}_{1\le j \le d}=0.$ We recall that $U_{\theta}$ is an isomorphism between $\cS(\bR^d)$ and $\cS(\rd),$ hence $\xi_j f(\xi)=0,1\le j\le d$ for any $\xi\in\bR^d,$ which implies $f=0$ and $u=0.$

Now, assume $u\ne 0,$ which implies $\|u\|_{L_p(\rd)}\ne 0$ and  $\|\nabla_{\theta} u\|_{L_p(\rd)}\ne 0.$ For any $t>0,$ $$\int_0^t \Delta_{\theta} e^{s\Delta_{\theta}}u\, ds=e^{t\Delta_{\theta}}u-u.$$
Therefore, by the triangle inequality and Proposition \ref{prop:rp}, we derive
    \begin{align*}
        \|u\|_{L_q(\rd)}&\le \lf\|e^{t\Delta_{\theta}}u\r\|_{L_q(\rd)}+\int_0^t\lf\|\Delta_{\theta} e^{s\Delta_{\theta}}u\r\|_{L_q(\rd)}\, ds\\
        &\lesssim t^{-\frac d2(\frac 1p-\frac 1q)}\|u\|_{L_p(\rd)}+\int_0^ts^{-\frac 12-\frac d2(\frac 1p-\frac 1q)}\|\nabla_{\theta} u\|_{L_p(\rd)}\,ds\\
        &\lesssim t^{-\frac{\vartheta}{2}}\|u\|_{L_p(\rd)}
        +t^{\frac{1-\vartheta}{2}}\|\nabla_{\theta} u\|_{L_p(\rd)}.
    \end{align*}
Setting  $t={\|u\|_{L_p(\rd)}^2}/{\|\nabla_{\theta} u\|_{L_p(\rd)}^2},$ we obtain the desired inequality for any $u\in \cs(\rd).$

    Now, for any fixed $u\in W_p^1(\rd),$ we can find a sequence $\{v_n\}\subset \cS(\rd)$ such that
    $$\|u-v_n\|_{W_p^1(\rd)}\le \frac 1n.$$
     By the Sobolev embedding $W_p^1(\rd)\subset L_q(\rd)$, which can be deduced by Lemmas \ref{hw} and \ref{se1}, we know that $$\lim_{n\to\fz}\|v_n-u\|_{L_q(\rd)}=0.$$ Since $$\|v_n\|_{L_q(\rd)}\lesssim \|v_n\|_{L_p(\rd)}^{1-\vartheta}\|\nabla_{\theta} v_n\|_{L_p(\rd)}^{\vartheta},$$ taking the limit as $n\to\fz$ yields the desired  inequality.
\end{proof}

\subsection{Time-space estimates for linear terms}
We delve into the study of the heat semigroup in the mixed space $L_{\gamma}(\rr_+; L_p(\rd))$ by utilizing the Littlewood-Paley decomposition alongside the exponential decay property of $e^{-t|\cdot|^2}$.
\begin{prop}\label{prop:gamma B}
Let $a\geq0$, $1\leq r\leq p\leq\infty$, $0<\lambda\leq\infty$ and $\frac2\gamma=a+d(\frac1r-\frac1p)$. Then one has
\begin{align}\label{gamma B}
\|Hx\|_{L_\gamma(\mathbb R_+;\dot{B}^0_{p,\lambda}(\rd))}\lesssim \|x\|_{\dot{B}^{-a}_{r,\lambda\wedge\gamma}(\rd)}
\end{align}
\end{prop}

\begin{proof}
By the density argument, we may assume $x=U_{\theta}(f)$, where $f\in\cs(\rr^d)$.
According to Lemma \ref{lem:psi}, we have:
\begin{align*}
\lf\|\dot{\triangle}_jH(t)x\r\|_{L_r(\rd)}= \lf\|U_{\theta}(\varphi_je^{-t|\cdot|^2}
f)\r\|_{L_r(\rd)}\leq\lf\|\mathcal{F}^{-1}(\varphi_je^{-t|\cdot|^2})
\r\|_{L_1(\rr^d)}\|x\|_{L_r(\rr^d)}.
\end{align*}
By the fact that $\varphi_j(\xi)=\varphi(2^{-j}\xi)$, $\supp \varphi\subset \{\xi\in\rr^d:\frac12\leq|\xi|\leq2\}$, one gets
\begin{align*}
\lf\|\mathcal{F}^{-1}(\varphi_je^{-t|\cdot|^2})\r\|_{L_1(\rr^d)}
\lesssim e^{-ct2^{2j}},
\end{align*}
where $c$ is an absolute constant.
 Combining these estimates, we obtain
\begin{align*}
\lf\|\dot{\triangle}_jH(t)x\r\|_{L_r(\rd)}\leq e^{-ct2^{2j}}\|x\|_{L_r(\rd)}.
\end{align*}
Let $\widetilde{\dot{\triangle}}_j={\dot{\triangle}}_{j-1}+{\dot{\triangle}}_j+{\dot{\triangle}}_{j+1}$.
Since $\dot{\triangle}_j=\dot{\triangle}_j\widetilde{\dot{\triangle}}_j$ and $\widetilde{\dot{\triangle}}_j$ shares the similar properties as ${\dot{\triangle}}_j$, the above estimate implies
 \begin{align}\label{ex decay}
\lf\|\dot{\triangle}_jH(t)x\r\|_{L_r(\rd)}\leq e^{-ct2^{2j}}\lf\|\dot{\triangle}_jx\r\|_{L_r(\rd)}.
\end{align}
Then taking the $\ell_\lambda$-norm over $j\in\zz$ in the inequality \eqref{ex decay}, we get
\begin{align}\label{inter}
\lf\|H(t)x\r\|_{\dot{B}^0_{r,\lambda}(\rd)}\lesssim \left(\sum_{j\in\zz}e^{-ct\lambda2^{2j}}\lf\|\dot{\triangle}_jx\r\|_{L_r(\rd)}^\lambda\right)^{\frac1\lambda}.
\end{align}
Now we divide the proof of the estimate \eqref{gamma B} into two cases: $\gamma\geq\lambda$ and $\gamma<\lambda$.

\textbf{Case 1.} In the first case $\gamma\geq\lambda$, by taking $L_\gamma(\mathbb R_+)$-norm of the inequality \eqref{inter} and applying the Minkowski inequality, we obtain:
\begin{align*}
\|Hx\|_{L_\gamma(\mathbb R_+;\dot{B}^0_{p,\lambda}(\rd))}
&\lesssim\left\|\sum_{j\in\zz}e^{-c\lambda t2^{2j}}\lf\|\dot{\triangle}_jx\r\|_{L_p(\rd)}^\lambda\right\|
_{L_{\gamma/\lambda}(\rr_+)}^{\frac1\lambda}\\
&\lesssim \left(\sum_{j\in\zz}\lf\|e^{-ct2^{2j}}\r\|_{L_{\gamma/\lambda}(\rr_+)}
\lf\|\dot{\triangle}_jx\r\|_{L_p(\rd)}^\lambda\right)^{\frac1\lambda}\\
&\lesssim \left(\sum_{j\in\zz}2^{-\frac{2\lambda j}{\gamma}}\lf\|\dot{\triangle}_jx\r\|_{L_p(\rd)}^\lambda\right)^{\frac1\lambda}
=\|x\|_{\dot{B}^{-2/\gamma}_{p,\lambda}(\rd)}.
\end{align*}
 By Lemma \ref{l2.81} (ii), we get $$\dot{B}^{-2/\gamma+d(1/r-1/p)}_{r,\lambda}(\rd)
\subset\dot{B}^{-2/\gamma}_{p,\lambda}(\rd).$$ Combining this inclusion with the above estimate, we conclude the desired inequality.

\textbf{Case 2.}
Now we consider the case where $\gamma<\lambda$.
Observe that, for any $\{\lambda_i\}_i\subset\mathbb{C}$ and $\vartheta\in(0,\,1]$, the inequality
 $$\lf(\sum_i|\lambda_i|\r)^{\vartheta}\leq\sum_i|\lambda_i|^{\vartheta}$$
 holds.
Utilizing this inequality, and taking the power of $\gamma$ and integrating in the inequality \eqref{inter}, we deduce
\begin{align*}
\int_{\mathbb R_+}\|H(t)x\|^\gamma_{\dot{B}^0_{p,\lambda}(\rd)}\,dt\leq\int_{\mathbb R_+}\sum_{j\in\zz}e^{-ct\gamma2^{2j}}\lf\|\dot{\triangle}_jx\r\|_{L_p(\rd)}^\gamma\,dt
\lesssim\sum_{j\in\zz}2^{-2j}\lf\|\dot{\triangle}_jx\r\|^\gamma_{L_p(\rd)}
=\|x\|_{\dot{B}^{-2/\gamma}_{p,\gamma}(\rd)}^{\gamma}.
\end{align*}
Again using Lemma \ref{l2.81} (ii), we have:$$\dot{B}^{-2/\gamma+d(1/r-1/p)}_{r,\gamma}(\rd)
\subset\dot{B}^{-2/\gamma}_{p,\gamma}(\rd).$$
Combining this inclusion with the above estimate, we conclude the desired result.
\end{proof}

 By taking  $\lambda=r=2$ in Proposition \ref{prop:gamma B} and utilizing the embedding
 $\dot{B}_{p,2}^0(\rd)\subset L_p(\rd)$ (see Lemma \ref{Besov embedding}), we can further get the following result.
\begin{cor}\label{cor:ts}
Let $2\leq p<\infty$ and $\frac2{\gamma(p)}=d(\frac12-\frac1p)$ with $\gamma(p)>2$. Then, we have
\begin{align}\label{ts1}
\|Hx\|_{L_{\gamma(p)}(\mathbb R_+;L_p(\rd))}\lesssim \|x\|_{L_2(\rd)}
\end{align}
and
\begin{align}\label{ts2}
\|\nabla_{\theta} Hx\|_{L_2(\mathbb R_+;L_2(\rd))}\lesssim \|x\|_{L_2(\rd)}.
\end{align}
\end{cor}

\bigskip
\subsection{Time-space estimates for the non-linear terms}

Now, we will establish time-space estimates for the non-linear terms.
To simplify the notations, we define
\begin{align}\label{def of H}
\mathcal H(f)(t):=\int^t_0H(t-s)f(s)\,ds.
\end{align}

\begin{prop}\label{prop:Hts1}
Let $1\leq r\leq p\leq\infty$ and $1<\gamma,\gamma_1<\infty$ be such that for $k=0,1$, the following conditions hold:
\begin{align}\label{expo i1}
\frac1\gamma=\frac1{\gamma_1}+\frac k2+\frac d2\lf(\frac 1r-\frac1p\r)-1,\ \ \frac k2+\frac d2\lf(\frac1r-\frac1p\r)<1.
\end{align}
Then, we have
\begin{align}\label{pHts1}
\lf\|\nabla_{\theta}^k\mathcal H(f)\r\|_{L_\gamma(\mathbb R_+;L_p(\rd))}\lesssim\|f\|_{L_{\gamma_1}(\mathbb R_+;L_r(\rd))}.
\end{align}
\end{prop}

\begin{proof}
By Proposition \ref{prop:rp}, we obtain
$$\lf\|\nabla_{\theta}^k\mathcal H(f)\r\|_{L_p(\rd)}
\lesssim\int_0^t(t-s)^{-\frac k2-\frac d2(\frac 1r-\frac 1p)}
\|f(s)\|_{L_r(\rd)}ds.$$
Utilizing this estimate along with the classical Hardy-Littlewood-Sobolev inequality, we can deduce:
\begin{align*}
\lf\|\nabla_{\theta}^k\mathcal H(f)\r\|_{L_\gamma(\mathbb R_+;L_p(\rd))}
&\lesssim\lf\|\int_0^t(t-s)^{-\frac k2-\frac d2(\frac 1r-\frac 1p)}
\|f(s)\|_{L_r(\rd)}\,ds\r\|_{L^t_\gamma(\mathbb R_+)}\\
&\lesssim\lf\|f\r\|_{L_{\gamma_1}(\mathbb R_+;L_r(\rd))},
\end{align*}
which implies the desired estimate.

\end{proof}

Noticing that Proposition \ref{prop:Hts1} does not include the endpoint case $\gamma=\infty$, and we can fill out the gap.

\begin{prop}\label{prop:Hts2}
Let $1\leq r<\infty$ and $1<q'\leq\lambda\leq\infty$.
Then one has
\begin{align}\label{pHts2}
\lf\|\mathcal Hf\r\|_{L_\infty(\mathbb R_+;\dot{B}^0_{r,\lambda}(\rr^d_{\theta}))}\lesssim\|f\|_{L_{q'}(\mathbb R_+;\dot{B}^{-2/q}_{r,\lambda}(\rd))}.
\end{align}
\end{prop}

\begin{proof}
Applying  \eqref{ex decay} and the classical Young inequality, we have
\begin{align*}
\lf\|\dot{\triangle}_j\mathcal Hf\r\|_{L_r(\rd)}&\lesssim \int^t_0 e^{-c(t-s)2^{2j}}\lf\|\dot{\triangle}_jf(s)\r\|_{L_r(\rd)}\,ds\\
&\lesssim 2^{-\frac{2j}q}\lf\|\dot{\triangle}_jf\r\|_{L_{q'}(\mathbb R_+;L_r(\rd))}.
\end{align*}
Then taking $\ell_\lambda$-norm on both sides and by the Minkowski inequality, one obtains
\begin{align*}
\lf\|\mathcal Hf\r\|_{\dot{B}^0_{r,\lambda}(\rd)}
&\lesssim\left(\sum_{j\in\zz}2^{-\frac{2\lambda j}q}\lf\|\dot{\triangle}_jf\r\|_{L_{q'}(\mathbb R_+;L_r(\rd))}^\lambda\right)^{\frac1\lambda}\\
&\lesssim \|f\|_{L_{q'}(\mathbb R_+;\dot{B}^{-2/q}_{r,\lambda}(\rd))},
\end{align*}
which yields the desired result.
\end{proof}

\begin{cor}\label{cor:Hts}
Let $2\leq p<\infty$ and $\frac2{\gamma(p)}=d(\frac12-\frac1p)$ with $\frac2{\gamma(p)}<1$. Then, we have
\begin{align}\label{Hts1}
\|\mathcal Hf\|_{L_\infty(\mathbb R_+;L_2(\rd))\cap L_{\gamma(p)}(\mathbb R_+;L_p(\rd))}\lesssim\|f\|_{L_{\gamma(p)'}(\mathbb R_+;L_{p'}(\rd))}
\end{align}
and
\begin{align}\label{Hts2}
\lf\|\nabla_{\theta}\mathcal H(f)\r\|_{L_2(\mathbb R_+;L_2(\rd))}\lesssim\|f\|_{L_{\gamma(p)'}(\mathbb R_+;L_{p'}(\rd))}.
\end{align}
\end{cor}
\begin{remark}{\rm
We point out that all the above estimates also hold when $\mathbb R_+$ is replaced by the interval $[0,T]$, $T>0$.}
\end{remark}

By Corollaries \ref{cor:ts} and \ref{cor:Hts}, we can get the following estimates.
\begin{prop}\label{prop:timespace0}
One has:
\begin{align}\label{timespace10}
\|Hx\|_{L_{2}([0,T];\dot{H}^1(\rd))}\leq C\|x\|_{L_2(\rd)};
\end{align}
\begin{align}\label{timespace20}
\|Hx\|_{L_{2+4/d}([0,T];L_{2+4/d}(\rd))}\leq C\|x\|_{L_2(\rd)};
\end{align}
\begin{align}\label{timespace00}
\|Hx\|_{L_{\infty}([0,T];L_{2}(\rd))}\leq C\|x\|_{L_2(\rd)};
\end{align}
\begin{align}\label{timespace30}
\lf\|\nabla_{\theta} \mathcal Hf\r\|_{L_{2}([0,T];L_2(\rd))}\leq C\|f\|_{L_{(2+4/d)'}([0,T];L_{(2+4/d)'}(\rd))};
\end{align}
\begin{align}\label{timespace40}
\| \mathcal Hf\|_{L_{\infty}([0,T];L_{2}(\rd))\cap {L_{2+4/d}([0,T];L_{2+4/d}(\rd))}}
\leq C\|f\|_{L_{(2+4/d)'}([0,T];L_{(2+4/d)'}(\rd))}.
\end{align}
\end{prop}

Moreover we have the following time-space estimates.
\begin{prop}\label{prop:timespace}
We have the following estimates:
\begin{align}\label{timespace1}
\|Hx\|_{L_{d+2}(\mathbb R_+;L_{d+2}(\rd))}\leq C\|x\|_{L_d(\rd)};
\end{align}
\begin{align}\label{timespace2}
\|Hx\|_{L_{\infty}(\mathbb R_+;L_{d}(\rd))}\leq C\|x\|_{L_d(\rd)};
\end{align}
\begin{align}\label{timespace3}
\lf\|\nabla_{\theta} \mathcal Hf\r\|_{L_{d+2}(\mathbb R_+;L_{d+2}(\rd))}\leq C\|f\|_{L_{(d+2)/2}
([0,\fz];[L_{(d+2)/2}(\rd)]^d)};
\end{align}
\begin{align}\label{timespace4}
\lf\|\nabla_{\theta} \mathcal Hf\r\|_{L_{\infty}(\mathbb R_+;L_{d}(\rd))}\leq C\|f\|_{L_{(d+2)/2}([0,\fz];[L_{(d+2)/2}(\rd)]^d)}.
\end{align}
\end{prop}

\begin{proof}
We set $p=r=\lambda=d$, $\gamma=2+d$ and $a=2/(d+2)$ in Proposition \ref{prop:gamma B}, we get
\begin{align*}
\|Hx\|_{L_{2+d}(\mathbb R_+;\dot{B}^{2/(2+d)}_{d,d}(\rd))}\lesssim\|x\|_{\dot{B}^0_{d,d}(\rd)}.
\end{align*}
Utilizing the embeddings $L_d(\rd)\subset \dot{B}^0_{d,d}(\rd)$ and $\dot{B}^{2/(2+d)}_{d,d}(\rd)\subset\dot{B}^{2/(2+d)}_{d,d+2}(\rd)\subset L_{d+2}(\rd)$ (see Lemma \ref{l2.81} (iii)), we obtain the desired estimate \eqref{timespace1}.

The estimate \eqref{timespace2} follows easily from the contraction of heat semigroup. The estiamte \eqref{timespace3} is a special case of Proposition \ref{prop:Hts1}. In fact, we only need to take $\gamma=p=2+d$ and $\gamma_1=r={(2+d)}/2$
in Proposition \ref{prop:Hts1}.

To obtain the estimate \eqref{timespace4}, we set $r=\lambda=(d+2)/2$ and $q'=(d+2)/2$ in Proposition \ref{prop:Hts2}, then,
\begin{align*}
\lf\|\mathcal Hf\r\|_{L_{\infty}(\mathbb R_+;\dot{B}^{2d/(2+d)}_{(d+2)/2,(d+2)/2}(\rd))}\lesssim \|f\|_{L_{(d+2)/2}(\mathbb R_+;\dot{B}^{0}_{(d+2)/2,(d+2)/2}(\rd))}.
\end{align*}
Combining this with Lemma \ref{bd}, we derive
\begin{align*}
\lf\|\nabla_{\theta}\mathcal Hf\r\|_{L_{\infty}(\mathbb R_+;\dot{B}^{(d-2)/(2+d)}_{(d+2)/2,(d+2)/2}(\rd))}\lesssim \|f\|_{L_{(d+2)/2}(\mathbb R_+;\dot{B}^{0}_{(d+2)/2,(d+2)/2}(\rd))}.
\end{align*}
Finally, utilizing the embeddings $L_{(d+2)/2}(\rd)\subset \dot{B}^0_{(d+2)/2,(d+2)/2}(\rd)$ and $\dot{B}^{(d-2)/(2+d)}_{(d+2)/2,(d+2)/2}(\rd)\subset L_{d}(\rd)$ (see Lemma \ref{l2.81} (iii)), we obtain the desired estimate \eqref{timespace4}.

\end{proof}

\section{The proof of Theorem \ref{t3.1}}\label{s6}
The present section is devoted to the proof of Theorem \ref{t3.1}.
We need two estimates involving the Leray projection \eqref{def of P}.
As in classical case, the Leray projection $\mathbb P=I+(-\Delta_{\theta})^{-1}\nabla_{\theta}\,{\bf div}$ can be regarded as a sequence of Fourier multipliers  $\mathbb{P}=(T_{m_{j,\,k}})_{1\leq j,\,k\leq d}$ in terms of the Weyl transform $U_{\theta}$ with symbols

\begin{align}\label{e5.32}
m_{j,\,k}(\xi)=\delta_{j,\,k}-\frac{\xi_j\xi_k}{
|\xi|^2}, \ \ \  \xi=(\xi_1,\xi_2,\ldots,\xi_d)\in \rr^d \setminus\{0\}, \ j,\,k\in[1,\,d]\cap\nn.
\end{align}

\begin{lemma}\label{map:P1}
For $\alpha\in\mathbb{R},1\leq p\leq\infty$, given $u\in [\cS'(\rd)]_0^d$ such that $u\otimes u\in [\dot{B}^\alpha_{p,\fz}(\rd)]^{d^2}$, we have
\begin{equation*}
\|\mathbb{P}A(u)\|_{[\dot B^{\alpha-1}_{p,\fz}(\rd)]^d}\lesssim\|u\otimes u\|_{[\dot{B}^\alpha_{p,\fz}(\rd)]^{d^2}}.
\end{equation*}
\end{lemma}

\begin{proof}
By the definition of $\mathbb{P}$ and the fact that ${\bf div} u=0$, for each $1\leq\ell\leq d$, we have
\begin{equation*}
(\mathbb{P}A(u))_\ell=\sum_{n=1}^d T_{m_{l,n}}\lf(\sum_{i=1}^du_i\partial_iu_n\r)=\sum_{n=1}^d T_{m_{l,n}}\lf(\sum_{i=1}^d\partial_i(u_iu_n)\r)=\sum_{1\leq i,n\leq d}T_{m_{l,n}}(\partial_i(u_iu_n)).
\end{equation*}
Then, by applying Lemma \ref{lem:psi}, we deduce
\begin{align*}
\|\mathbb{P}A(u)\|_{[\dot B^{\alpha-1}_{p,\fz}(\rd)]^d}^2
&=\sum_{\ell=1}^d\lf(\sup_{j\in \bZ}2^{j(\alpha-1)}\lf\|\dot\triangle_j(\mathbb{P}A(u))_\ell\r\|_{L_{p}(\rd)}\r)^2\\
&\lesssim\sum_{\ell=1}^d \lf(\sup_{j\in \bZ}2^{j(\alpha-1)}\sum_{1\leq i,n\leq d}\lf\|\cF^{-1}(\xi_im_{l,n}\varphi_j)\r\|_{L_1(\mathbb{R}^d)}
\lf\|\dot\triangle_j (u_iu_n)\r\|_{L_{p}(\rd)}\r)^2\\
&\lesssim\sum_{\ell=1}^d\lf(\sum_{1\leq i,n\leq d}\sup_{j\in \bZ}2^{j\alpha}\lf\|\dot\triangle_j (u_iu_n)\r\|_{L_{p}(\rd)}\r)^2\\
&\lesssim\sum_{1\leq i,n\leq d}\lf(\sup_{j\in \bZ}2^{j\alpha}\lf\|\dot\triangle_j (u_iu_n)\r\|_{L_{p}(\rd)}\r)^2=\|u\otimes u\|_{[\dot{B}^\alpha_{d,\fz}(\rd)]^{d^2}}^2.
 \end{align*}
\end{proof}

\begin{lemma}\label{t5.1}
Let $p\in(1,\,\fz). $   Then the Leray projection $\mathbb{P}$ is bounded from
$[L_p(\rd)]^d$ to itself.
\end{lemma}
\begin{proof}
Since the Leray projection  $\mathbb{P}$ is a sequence of Fourier multipliers with symbols $(m_{j,k})_{j,k}$ in \eqref{e5.32}, all of which are the linear combination of the identity operator and (composite) Riesz transforms, the assertion follows then from the transference principle---Theorem \ref{thm:trans} and noncommutative Calder\'on-Zygmund theory (cf. e.g. \cite[Theorem 6.4]{Mei07} or \cite{{MR2476951}}).
\end{proof}

With the previous preparatory work, we now start the proof of Theorem \ref{t3.1}, where the contraction mapping principle will be frequently exploited.

\begin{proof}[Proof of Theorem \ref{t3.1}]
{\bf Part (i):}
To establish part (i) of Theorem \ref{t3.1}, we initially seek a solution $u\in L^{\mathrm{loc}}_{d+2}([0,T_{u_0});$ $ [L_{d+2}(\rd)]^d_0)\cap L^{\mathrm{loc}}_{\fz}([0,T_{u_0});$ $[L_d(\rd)]_0^d)$ for some maximal time $T_{u_0}$. Subsequently, we demonstrate that this solution actually resides in  $C([0,T_{u_0});[L_{d}(\rd)]^d_0)$. For clarity, we present the proof through the following seven steps:

\textbf{Step 1.} {\bf(Existence)} In this step, we aim to demonstrate the existence of a solution to the Navier-Stokes equation \eqref{ens} given that $u_0\in[L_d(\rd)]_0^d$.
Let $\delta>0$ be a constant to be determined later, and by the inequality \eqref{timespace1}, we can find $T_{\delta}\triangleq T>0$ such that
$$\|Hu_0\|_{L_{d+2}([0,T];[L_{d+2}(\rd)]^d)}\leq\frac{\delta}{4}.$$
We will utilize the contraction mapping principle in the following set:
$$\mathcal D_\delta=\lf\{u\in L_{d+2}([0,T];[L_{d+2}(\rd)]^d):\;\|u\|_{{L_{d+2}([0,T];[L_{d+2}(\rd)]^d)}}\leq \delta,\ \ \|u\|_{L_{\fz}([0,T];[L_d(\rd)]^d)}\leq\|u_0\|_{[L_d(\rd)]^d}+\delta\r\},$$
equipped with the distance
$$  d(u,v):=\|u-v\|_{L_{d+2}([0,T];[L_{d+2}(\rd)]^d)},$$
for some sufficiently small $\delta>0$.
We claim that the set $(\mathcal D_\delta,\,d)$ is complete. Indeed,
let $\{u_j\}_j$ be a Cauchy sequence in $(\mathcal D_\delta,\,d)$. Then, there exists $u\in L_{d+2}([0,T];[L_{d+2}(\rd)]^d)$ such that
\begin{align*}
u_j\rightarrow u \ \ \mathrm{in} \ \ L_{d+2}([0,T];[L_{d+2}(\rd)]^d),\ \
\mathrm{as} \ \ j\rightarrow\fz,
\end{align*}
and $\|u\|_{L_{d+2}([0,T];[L_{d+2}(\rd)]^d)}\leq\delta$.
On the other hand, $\{u_j\}_j$ is bounded in $L_{\fz}([0,T];[L_d(\rd)]^d)$, then $\exists\,\widetilde{u}\in L_{\fz}([0,T];[L_d(\rd)]^d)$ and a subsequence $\{u_{j_n}\}_n$, such that
 $u_{j_n}\rightarrow\widetilde{u}$  in the weak-$*$ topology on $L_{\fz}([0,T];[L_d(\rd)]^d)$ as $n\rightarrow\fz$.
 It is easy to see that $u=\widetilde{u}$ and
 $$\|u\|_{L_{\fz}([0,T];[L_d(\rd)]^d)}\leq\liminf_{n\rightarrow\fz}
 \lf\|u_{j_n}\r\|_{L_{\fz}([0,T];[L_d(\rd)]^d)}\leq\|u_0\|_{[L_d(\rd)]^d}+\delta.$$
 Therefore, we obtain that $u\in \mathcal{D}_{\delta}$ and $d(u_j,\,u)\rightarrow0$ as $j\rightarrow\fz$, which implies that the metric space $(\mathcal{D}_\delta,\,d)$ is complete.

Next, we consider the solution map:
\begin{align}\label{map1}\mathfrak Mu:=Hu_0-\mathcal H\mathbb P\lf(A(u)\r).
\end{align}
By Lemma \ref{t5.1}, the inequality \eqref{timespace3}, Remark \ref{r0} and the H\"{o}lder inequality, we obtain
\begin{eqnarray}\label{sm1.1}
\|\mathfrak Mu\|_{{L_{d+2}([0,T];[L_{d+2}(\rd)]^d)}}
&\leq&\|Hu_0\|_{{L_{d+2}([0,T];[L_{d+2}(\rd)]^d)}}+
\lf\|\mathcal H\mathbb P\lf(A(u)\r)\r\|_{{L_{d+2}([0,T];[L_{d+2}(\rd)]^d)}}\\\nonumber
&\overset{\eqref{timespace3}}{\leq}& \lf\|Hu_0\r\|_{{L_{d+2}([0,T];[L_{d+2}(\rd)]^d)}}
+C\|u\otimes u\|_{{L_{(d+2)/2}([0,T];[L_{(d+2)/2}(\rd)]^{d^2})}} \\ \nonumber
&\leq&\frac{\delta}{4}+C\|u\|^2_{{L_{d+2}([0,T];[L_{d+2}(\rd)]^d)}} \\ \nonumber
&\leq&\frac{\delta}{4}+C\delta^2,
\end{eqnarray}
and
\begin{eqnarray}\label{sm1.2}
\|\mathfrak Mu\|_{{L_{\fz}([0,T];[L_{d}(\rd)]^d)}}
&\leq&\|Hu_0\|_{{L_{\fz}([0,T];[L_{d}(\rd)]^d)}}+
\lf\|\mathcal H\mathbb P\lf(A(u)\r)\r\|_{{L_{\fz}([0,T];[L_{d}(\rd)]^d)}}\\ \nonumber
&\overset{\eqref{timespace3}}\leq&\|u_0\|_{[L_d(\rd)]^d}+C\|u\|^2_{{L_{d+2}([0,T];[L_{d+2}(\rd)]^d)}}\\ \nonumber
&\leq&\|u_0\|_{[L_d(\rd)]^d}+C\delta^2. \nonumber
\end{eqnarray}
From this, we can further prove that $\mathfrak{M}$ is a contraction mapping from $(\mathcal{D}_{\delta},\,d)$ to itself. In fact, it follows from the H\"{o}lder inequality that
\begin{align*}
\lf\|\mathfrak Mu-\mathfrak M v\r\|_{{L_{d+2}([0,T];[L_{d+2}(\rd)]^d)}}&\leq C\lf(\|u\|_{{L_{d+2}([0,T];[L_{d+2}(\rd)]^d)}}+\|v\|_
{{L_{d+2}([0,T];[L_{d+2}(\rd)]^d)}}\r)\\
&\ \ \ \ \times\|u-v\|_{{L_{d+2}([0,T];[L_{d+2}(\rd)]^d)}}\\
&\leq 2C\delta\|u-v\|_{{L_{d+2}([0,T];[L_{d+2}(\rd)]^d)}}.
\end{align*}
Choosing the constant $\delta$ small enough such that $C\delta< \frac18$, we get
$$\|\mathfrak Mu\|_{{L_{d+2}([0,T];[L_{d+2}(\rd)]^d)}}\leq \delta,\ \ \ \ \|\mathfrak Mu\|_{{L_{\fz}([0,T];[L_{d}(\rd)]^d)}}
\leq \|u_0\|_{[L_d(\rd)]^d}+\delta$$
and $$d(\mathfrak Mu,\mathfrak Mv)\leq\frac14d(u,v).$$
Thus, there exists a unique solution in $ (\mathcal{D}_\delta,\,d)$ which satisfies the NS equation \eqref{ens}. The fact ${\bf div}\, u=0$ can be verified easily from ${\bf div}\, u_0=0$.

\textbf{Step 2.} {(\bf The solution $u\in C([0,T];[L_d(\rd)]_0^d)$)}
In this step, we will show the solution $u$ found above is actually in $C([0,T];[L_d(\rd)]_0^d)$.

Now, we will show that $u\in C([0,T];[L_d(\rd)]^d_0)$.
For any $t_1,\,t_2\in[0,T]$ with $t_1<t_2$, we see that
\begin{align*}
u(t_2)-u(t_1)&=H(t_2)u_0-H(t_1)u_0
+
\int_0^{t_1}H(t_1-s)\mathbb P\lf(A(u)\r)(s)\,ds
 -\int_0^{t_2}H(t_2-s)\mathbb P\lf(A(u)\r)(s)\,ds\\
&=(H(t_2-t_1)-I)H(t_1)u_0
 +
\int_{0}^{t_1}(H(t_2-t_1)-I)H(t_1-s)\mathbb P\lf(A(u)\r)(s)\,ds\\
&\ \ \ \ \  +\int_{t_1}^{t_2}H(t_2-s)\mathbb P\lf(A(u)\r)(s)\,ds\\
&=(H(t_2-t_1)-I)u(t_1)
+\int_0^{t_2-t_1}H(t_2-t_1-s)\mathbb P\lf(A(u)\r)(t_1+s)\,ds,
\end{align*}
which, combined with the inequality \eqref{timespace40}, implies that
\begin{align*}
\|u(t_2)-u(t_1)\|_{[L_d(\rd)]^d}
&\leq\lf\|(H(t_2-t_1)-I)u(t_1)\r\|_{[L_d(\rd)]^d}\\
&\ \ \ \ \ +\lf\|
\int_0^{t_2-t_1}H(t_2-t_1-s)\mathbb P\lf(A(u)\r)(t_1+s)\,ds\r\|_{[L_d(\rd)]^d}.
\end{align*}
By Lemma \ref{heat}, Lemma \ref{t5.1} and the equality \eqref{timespace4}, we have
\begin{align*}
\lim_{t_2\rightarrow t_1}\|u(t_1)-u(t_2)\|_{[L_d(\rd)]^d}
\lesssim&\lim_{t_2\rightarrow t_1}\lf\|(H(t_2-t_1)-I)u(t_1)\r\|_{[L_d(\rd)]^d}\\
& \ \ \ +\lim_{t_2\rightarrow t_1}
\lf\|u \otimes u \r\|_{L_{(d+2)/2}([t_1,t_2];[L_{(d+2)/2}(\rd)]^{d^2})}\\
\lesssim&\lim_{t_2\rightarrow t_1}
\|u\|^2_{L_{d+2}([t_1,t_2]; [L_{d+2}(\rd)]^d)}\\
=&0,
\end{align*}
which implies that $u\in C([0,T];[L_d(\rd)]^d_0)$.

\textbf{Step 3.} {\bf(Uniqueness)}
Now, we show that the solution $u$ is unique in
$C([0,T];[L_d(\rd)]^d_0) \cap L_{d+2}([0,T];[L_{d+2}(\rd)]_0^d)$.
Let us suppose that  $v$ is another solution and
 $$v\in C([0,T];[L_d(\rd)]^d_0) \cap L_{d+2}([0,T];[L_{d+2}(\rd)]_0^d).$$
 Consequently, there exists a $T_1<T$ for which
$$\|v\|_{{L_{d+2}([0,T_1];[L_{d+2}(\rd)]^d)}}\leq \delta$$
and
$$ \|v\|_{L_{\fz}([0,T_1];[L_d(\rd)]^d)}\leq\|u_0\|_{[L_d(\rd)]^d}+\delta.$$
Utilizing the results from {\bf Step 1}, we deduce that
\begin{align*}
v(t)=u(t), \ \ \ \mathrm{when} \ \ 0\leq t \leq T_1.
\end{align*}
Thus, there exists a maximal time $T^*\in[0,T]$ such that
\begin{align*}
v(t)=u(t), \ \ \ \mathrm{for \ \ all} \ \ 0\leq t \le T^*.
\end{align*}
If $T^*=T$, our proof is complete. Alternatively, if $T^*<T$,
we consider the following NS equation starting with the initial datum
$u(T^*)$.
\begin{align}\label{e3.aa}
w(t)=H(t-T^*)u(T^*)-\int^t_{T^*}H(t-s)\mathbb P\lf(A(w)\r)(s)\,ds,\ \ \ w(T^*)=u(T^*).
\end{align}
By replicating the procedures outlined in {\bf Step 1} and {\bf Step 2}, we know that there exists a small positive number $\epsilon$ such that $T^*<T^*+\epsilon <T$ and the equation \eqref{e3.aa} has a unique solution $w$ in
$C([T^*,T^*+\epsilon];[L_d(\rd)]^d_0)\cap L_{d+2}([T^*,T^*+\epsilon];[L_{d+2}(\rd)]^d)$.
Moreover, the solutions $u,v$ satisfy
$$\|u\|_{{L_{d+2}([T^*,T^*+\epsilon];[L_{d+2}(\rd)]^d)}}\leq \delta, \ \ \
 \|u\|_{L_{\fz}([T^*,T^*+\epsilon];[L_d(\rd)]^d)}\leq\|u(T^{*})\|_{[L_d(\rd)]^d}+\delta$$
 and
 $$\|v\|_{{L_{d+2}([T^*,T^*+\epsilon];[L_{d+2}(\rd)]^d)}}\leq \delta, \ \ \
 \|v\|_{L_{\fz}([T^*,T^*+\epsilon];[L_d(\rd)]^d)}\leq\|u(T^{*})\|_{[L_d(\rd)]^d}+\delta.$$
Therefore,
\begin{align*}
v(t)=u(t)=w(t), \ \ \ \mathrm{when} \ \ T^*\leq t < T^*+\epsilon,
\end{align*}
which contradicts the maximality of $T^*$. Therefore, we complete the proof of the above claim.

\textbf{Step 4.} {\bf(Existence of maximal time  $T_{u_0}$)}
In this step, we show that there exists a maximal time $T_{u_0}$ such that the above solution $u\in C([0,T_{u_0}); [L_d(\rd)]^d_0)
\cap{L_{d+2}^{\mathrm{loc}}([0,T_{u_0});[L_{d+2}(\rd)]_0^d)}$.
Given $u_0\in [L_d(\rd)]^d_0$, from the previous arguments, we know there exist $T$ and a unique solution $$u\in C([0,T]; [L_d(\rd)]^d_0)
\cap{L_{d+2}([0,T];[L_{d+2}(\rd)]_0^d)}.$$ We now extend $T$ to a maximal time, denoted by $T_{u_0}$. Specifically, we consider the NS equation with initial datum $u(T)$,
\begin{align*}
v(t)=H(t-T)v(T)-\int^t_{T}H(t-s)\mathbb P\lf(A(v)\r)(s)\,ds,\ \ \ v(T)=u(T).
\end{align*}
As established earlier, there exists a small positive number $\varepsilon$ and a unique solution $v\in C([T, T+\varepsilon]; [L_d(\rd)]^d_0)\cap{L_{d+2}([T, T+\varepsilon]; [L_{d+2}(\rd)]^d_0)}.$ Then, it is straightforward to verify that
\begin{eqnarray*}\tilde{u}(t):=
\lf\{\begin{array}{ll}
u(t),  \ \ \mathrm{when} \ \ t\le T;\\[5pt]
v(t), \ \ \mathrm{when} \ \  T\leq t\le T+\varepsilon.
\end{array}\r.
\end{eqnarray*}
 is the unique solution to the NS equation up to time $T+\varepsilon,$ that is $u\in C([0, T+\varepsilon]; [L_d(\rd)]^d_0)\cap{L_{d+2}([0, T+\varepsilon]; [L_{d+2}(\rd)]^d_0)}.$ Therefore, we can extend the solution step by step and identify a maximal time $T_{u_0}$ such that $u\in C([0,T_{u_0}); [L_d(\rd)]^d_0)\cap{L_{d+2}^{\mathrm{loc}}([0,T_{u_0});[L_{d+2}(\rd)]_0^d)}.$

Now we claim that if $T_{u_0}<\fz,$ then we have $$\|u\|_{L_{d+2}([0,T_{u_0});[L_{d+2}(\rr_{\theta}^d)]^d)}=\infty.$$ If not, then we can find $0\le t^*<T_{u_0}$ such that $$\|u\|_{L_{d+2}([t^*,T_{u_0});[L_{d+2}(\rr_{\theta}^d)]^d)}\le \frac \delta {100}.$$
By the triangle inequality, the inequality \eqref{timespace3}, Lemma \ref{t5.1} and the H\"{o}lder inequality, we obtain
\begin{align*}
\|Hu_0\|_{L_{d+2}([t^*,T_{u_0});[L_{d+2}(\rr_{\theta}^d)]^d)}
&\leq \|u\|_{L_{d+2}([t^*,T_{u_0});[L_{d+2}(\rr_{\theta}^d)]^d)}+
\lf\|\mathcal H\mathbb P\lf(A(u)\r)\r\|_{L_{d+2}([t^*,T_{u_0});[L_{d+2}(\rr_{\theta}^d)]^d)}\\
&\le\|u\|_{L_{d+2}([t^*,T_{u_0});[L_{d+2}(\rr_{\theta}^d)]^d)}+C\|u\|_{L_{d+2}([t^*,T_{u_0});[L_{d+2}(\rr_{\theta}^d)]^d)}^2\\
&\le \frac \delta {50}.
\end{align*}
Notice $$\|Hu_0\|_{L_{d+2}(\bR^+;[L_{d+2}(\rr_{\theta}^d)]^d)}\le \|u_0\|_{[L_d(\rd)]^d},$$ then we can find $t^{**}>T_{u_0}$ such that $$\|Hu_0\|_{L_{d+2}([t^*,t^{**}];[L_{d+2}(\rr_{\theta}^d)]^d)}\le\frac \delta 4.$$ Via the same process as {\bf Step 1}, we can extend $u$ up to time $t^{**},$ which is in contradiction with the maximality of $T_{u_0}.$

\textbf{Step 5.} {\bf(The solution $u$ belongs to $\bigcap_{d\leq p\leq\infty}C((0,T_{u_0});[L_p(\rd)]_0^d)$)} In this step, we will demonstrate that the solution $u$ belongs to the space $C((0,T_{u_0});[L_p(\rd)]_0^d)$ for any $d<p\le \fz$.
Firstly, we will prove that for any $w\in L_{d}(\rd)$ and $d<p\le \fz,$ the following conclusion holds:
\begin{align}\label{elim1}
\lim_{t\rightarrow0}t^{\frac{1}{2}-\frac d{2p}}\|H(t)w\|_{L_{p}(\rd)}=0.
\end{align}
When $w\in L_{d}(\rd)\cap L_{p}(\rd)$, we have
\begin{align*}
\lim_{t\rightarrow0}t^{\frac{1}{2}-\frac d{2p}}\|H(t)w\|_{L_{p}(\rd)}\leq\lim_{t\rightarrow0}t^{\frac{1}{2}-\frac d{2p}}\|w\|_{L_{p}(\rd)}=0.
\end{align*}
Observe that for any $w_1,w_2\in L_{d}(\rd)$,  by Proposition \ref{prop:rp}, we have
\begin{align*}
  \Big| t^{\frac{1}{2}-\frac d{2p}}\|H(t)w_1\|_{L_{p}(\rd)}-t^{\frac{1}{2}-\frac d{2p}}\|H(t)w_2\|_{L_{p}(\rd)}\Big|\le t^{\frac{1}{2}-\frac d{2p}}\|H(t)(w_1-w_2)\|_{L_{p}(\rd)}\le \left \| w_1-w_2\right \|_{L_d(\rd)}.
\end{align*}
Given this inequality and the fact that $L_{d}(\rd)\cap L_{p}(\rd)$ is dense in $L_{d}(\rd)$, we obtain \eqref{elim1} for all $w\in L_d(\rd)$.

Next, as in {\bf Step 1}, for $d<p\leq \infty$, let $\delta>0$ be a constant to be determined later. By the equality \eqref{elim1}, we can find $T_{\delta}\triangleq T>0$ such that
$$\|Hu_0\|_{L_{d+2}([0,T];[L_{d+2}(\rd)]^d)}+\lf\|(\cdot)^{\frac{1}{2}-\frac d{2p}}Hu_0\r\|_{L_{\infty}([0,T];[L_{p}(\rd)]^d)}\leq\frac{\delta}{2}.$$
Define the set
$$\mathcal{D}^\prime_\delta=\lf\{u: \|u\|_{L_{d+2}([0,T]; [L_{d+2}(\rd)]^d)}+\lf\|(\cdot)^{\frac{1}{2}-\frac d{2p}}u(\cdot)\r\|_{L_{\infty}([0,T];[L_{p}(\rd)]^d)}\leq \delta, \
\|u\|_{L_{\fz}([0,T];[L_d(\rd)]^d)}\leq\|u_0\|_{[L_d(\rd)]^d}+\delta\r\}$$
equipped with the distance
$$d^\prime(u,v):=\|u-v\|_{L_{d+2}([0,T];[L_{d+2}(\rd)]^d)}
+\lf\|(\cdot)^{\frac{1}{2}-\frac d{2p}}u(\cdot)-(\cdot)^{\frac{1}{2}-\frac d{2p}}v(\cdot)\r\|_{L_{\infty}([0,T];[L_{p}(\rd)]^d)}
.$$
We now show that the solution map $\mathfrak{M}$ (see \eqref{map1}) is a contraction on $(\mathcal{D}'_\delta,d^\prime)$. As in {\bf Step 1}, it suffices to show the following estimates:
\begin{align}\label{set close}\lf\|(\cdot)^{\frac{1}{2}-\frac d{2p}}\mathcal H\mathbb P(A(u))(\cdot)\r\|_{ L_\infty([0,T];[L_{p}(\rd)]^d)}\leq\frac \delta4\end{align} and
\begin{align}\label{contraction2}\lf\|(\cdot)^{\frac{1}{2}-\frac d{2p}}(\mathfrak Mu(\cdot)-\mathfrak Mv(\cdot))\r\|_{L_\infty([0,T];[L_{p}(\rd)]^d)}\leq\frac14 d'(u,\,v).
\end{align}
To prove the inequality \eqref{set close}, we can write
\begin{align*}
&\mathcal H\mathbb P(A(u))(t)=
\int^t_{0}H(t-s)\mathbb P(A(u))(s)\,ds\\
=&\int^{\frac t2}_{0}H(t-s)\mathbb P(A(u))(s)\,ds+\int^t_{\frac t2}H(t-s)\mathbb P(A(u))(s)\,ds\\
=&:\mathrm{I}(t)+\mathrm{II}(t).
\end{align*}
For the term $\mathrm{I}(t)$, we further express it as
\begin{align*}
\mathrm{I}(t):=H(\frac t2)\mathcal H\mathbb P(A(u))(\frac t2).
\end{align*}
By Proposition \ref{prop:rp} and the inequality \eqref{sm1.2},
\begin{align*}
\sup_{0\le t\le T}t^{\frac{1}{2}-\frac d{2p}}\|\mathrm{I}(t)\|_{[L_{p}(\rd)]^d}&\leq  \sup_{0\le t\le T}\Big\|\mathcal H\mathbb P(A(u))(\frac t2)\Big\|_{[L_{d}(\rd)]^d}\leq C\lf\|u\r\|_{L_{d+2}([0,T];[L_{d+2}(\rd)]^d)}^2\leq C\delta^2.
\end{align*}
For the second term $\mathrm{II}(t)$, applying Lemma \ref{lem:psi} to the Fourier multiplier $H(t-s)\mathbb P{\bf div}$ on $[L_{p}(\rd)]^{d^2}$, we obtain
\begin{align*}
\sup_{0\le t\le T}t^{\frac{1}{2}-\frac d{2p}}\|\mathrm{II}(t)\|_{[L_{p}(\rd)]^d}
&\leq C \sup_{0\le t\le T}t^{\frac{1}{2}-\frac d{2p}}\int^t_{\frac t2}\lf\|H(t-s)\mathbb P{\bf div}\, (u\otimes u)(s)\r\|_{[L_{p}(\rd)]^d} \,ds\\
&\leq C \sup_{0\le t\le T}t^{\frac{1}{2}-\frac d{2p}}\int^t_{\frac t2}\lf\|u(s)\otimes u(s)\r\|_{[L_{p}(\rd)]^{d^2}}\frac{ds}{(t-s)^{\frac 12+\frac d{2p}}}\\
&\leq C \sup_{0\le t\le T}\lf\|t^{\frac{1}{2}-\frac d{2p}}u(t)\r\|^2_{[L_{p}(\rd)]^d}\\
&\leq C\delta^2.
\end{align*}
Therefore, by the triangle inequality, we deduce
\begin{equation*}
 \lf\|{(\cdot)}^{\frac{1}{2}-\frac d{2p}}\mathcal H\mathbb P(A(u))(\cdot)\r\|_{ L_\infty([0,T];[L_{p}(\rd)]^d)}\leq 2C\delta^2.
\end{equation*}
Arguing as above, and by the H\"older inequality, we have
\begin{align*}
&\lf\|t^{\frac{1}{2}-\frac d{2p}}(\mathfrak Mu(t)-\mathfrak Mv(t))\r\|_{L_\infty([0,T];[L_{p}(\rd)]^d)}\\
\leq &C\delta\lf(\|u-v\|_{{L_{d+2}([0,T];[L_{d+2}(\rd)]^d)}}+\lf\|t^{\frac{1}{2}-\frac d{2p}}(u-v)(t)\r\|_{L_\infty([0,T];[L_{p}(\rd)]^d)}\r).
\end{align*}
Choosing the constant $\delta$ small enough such that $C\delta<\frac18$, we get the estimates \eqref{set close} and \eqref{contraction2}.
Now, by the contraction mapping principle, we can find a solution $$\tilde{u}\in  L_{d+2}([0,T]; [L_{d+2}(\rd)]^d_0)\cap L_{\fz}([0,T];[L_d(\rd)]_0^d)$$  satisfying $$(\cdot)^{\frac{1}{2}-\frac{d}{2p}}\tilde{u}\in L_{\fz}([0,T];[L_{p}(\rd)]_0^d).$$ By the uniqueness of solution in $L_{d+2}([0,T]; [L_{d+2}(\rd)]^d_0)\cap L_{\fz}([0,T];[L_d(\rd)]_0^d),$ we have $u=\tilde u$ on $[0,T].$
As in {\bf Step 4}, we can find a maximal time $T_{u_0,p}$, such that $\tilde{u}$ can be extended to  $$\tilde{u}\in  L_{d+2}([0,T_{u_0,p}); [L_{d+2}(\rd)]^d_0)\cap L_{\fz}([0,T_{u_0,p});[L_d(\rd)]_0^d)$$  satisfying $$(\cdot)^{\frac 12-\frac{d}{2p}}\tilde{u}\in L_{\fz}([0,T];[L_{p}(\rd)]_0^d)$$ for any $0<T<T_{u_0,p}$. Moreover, we have $T_{u_0,p}\le T_{u_0}$ and $u=\tilde u$ on $[0,T_{u_0,p}).$

We can also show that $u\in C((0,T_{u_0,p});[L_{p}(\rd)]_0^d)$ by the fact that for any $0<t_1<t_2<T_{u_0,p},$
\begin{align*}
    \|u(t_2)-u(t_1)\|_{[L_{p}(\rd)]^d}
&\leq\lf\|(H(t_2-t_1)-I)u(t_1)\r\|_{[L_{p}(\rd)]^d}\\
&\ \ \ \ \ +\lf\|
\int_0^{t_2-t_1}H(t_2-t_1-s)\mathbb P\lf(A(u)\r)(t_1+s)\,ds\r\|_{[L_{p}(\rd)]^d}\\
&\leq \|(H(t_2-t_1)-I)u(t_1)\|_{[L_{p}(\rd)]^d}+(t_2-t_1)^{\frac12-\frac d{2p}}\|u\|_{L_{\fz}([t_1,t_2];[L_d(\rd)]^d)}.
\end{align*}
Hence, it suffices to show that $T_{u_0,p}=T_{u_0}$ for all $d<p\leq\infty$.

We proceed to demonstrate that $T_{u_0,d+2}=T_{u_0}.$ Suppose for contradiction that $T_{u_0,d+2}<T_{u_0},$ then we have $$u\in  L_{d+2}([0,T_{u_0,d+2}]; [L_{d+2}(\rd)]^d_0).$$ Moreover, for any $t\in (0,T_{u_0,d+2}),$ we claim that
\begin{align}\label{cliam_29-1}
    \|u(t)\|_{[L_{d+2}(\rd)]^d}\ge \frac{\delta}{4(T_{u_0,d+2}-t)^{\frac{1}{d+2}}}.
\end{align}
By the inequality \eqref{cliam_29-1}, we deduce that
\begin{align*}
    \lf\|u\r\|_{L_{d+2}([0,T_{u_0,d+2}];[L_{d+2}(\rd)]^d)}^{d+2}=\int_{0}^{T_{u_0,d+2}}\|u(t)\|_{[L_{d+2}(\rd)]^d}^{d+2}\,dt\gtrsim \int_{0}^{T_{u_0,d+2}}\frac{1}{\lf(T_{u_0,d+2}-t\r)}\,dt =\fz.
\end{align*}
Obviously, we find that, the above divergent integral is in contradiction with the fact
that $u\in  L_{d+2}([0,T_{u_0,d+2}]; [L_{d+2}(\rd)]^d_0),$ hence we have $T_{u_0,d+2}=T_{u_0}.$
Now, we prove the claim \eqref{cliam_29-1}. By Lemma \ref{heat}, we have
\begin{align*}
    &\|H(\cdot-t)u(t)\|_{L_{d+2}([t,T_{u_0,d+2}];[L_{d+2}(\rd)]^d)}+\lf\|(\cdot-t)^{\frac 1{d+2}}H(\cdot-t)u(t)\r\|_{L_{\infty}([t,T_{u_0,d+2}];[L_{d+2}(\rd)]^d)}\\
    \le &2\lf(T_{u_0,d+2}-t\r)^{\frac{1}{d+2}}\|u(t)\|_{[L_{d+2}(\rd)]^d}.
\end{align*}
Then, we must have $$\|u(t)\|_{[L_{d+2}(\rd)]^d}\ge \frac{\delta}{4\lf(T_{u_0,d+2}-t\r)^{\frac{1}{d+2}}},$$ otherwise \begin{align*}
    \|H(\cdot-t)u(t)\|_{L_{d+2}([t,T_{u_0,d+2}];[L_{d+2}(\rd)]^d)}+\lf\|(\cdot-t)^{\frac 1{d+2}}H(\cdot-t)u(t)\r\|_{L_{\infty}([t,T_{u_0,d+2}];[L_{d+2}(\rd)]^d)}\le \frac \delta 2
\end{align*}
and we can extend $u$ to a time $T_{u_0,d+2}\le T'<T_{u_0}$ as in the previous step, which leads to a contradiction.

Finally, we show $T_{u_0,\fz}=T_{u_0}.$ By contradiction, we assume that this does not hold. Subsequently, by analogous reasoning, we obtain that
$$\|u(t)\|_{[L_{\fz}(\rd)]^d}\ge \frac{\delta}{4\lf(T_{u_0,\fz}-t\r)^{\frac{1}{2}}},$$ which implies that $$\lim_{t\nearrow T_{u_0,\fz}}\|u(t)\|_{[L_{\fz}(\rd)]^d}=\fz.$$ Fix $\lambda\in (0,T_{u_0,\fz}).$ For $t\in [0,\lambda],$ it trivially holds that $\sqrt{\cdot}u(\cdot)\in L_{\infty}([0,\lambda];[L_{\fz}(\rd)]^d_0).$
 Define \begin{align*}
    M:=\lf\|(\cdot)^{\frac 1{d+2}}u(\cdot)\r\|_{L_{\fz}([0,T_{u_0,\fz}];[L_{d+2}(\rd)]^d)}
  \end{align*}
  and
  \begin{align*}
    N:=\lf\|\sqrt \cdot u(\cdot)\r\|_{L_{\infty}([0,\lambda];[L_{\fz}(\rd)]^d)},
\end{align*}
Consider $\varepsilon\in (0,1)$ sufficiently small and $t\in (\lambda,T_{u_0,\fz}).$
Then we have
\begin{align*}
    \sqrt t\|u(t)\|_{[L_\infty(\rd)]^d}&\le \sqrt t\|H(t)u_0\|_{[L_\infty(\rd)]^d}+\sqrt t\lf(\int_0^{(1-\varepsilon)t}+\int_{(1-\varepsilon)t}^t\r) \|H(t-s)\mathbb P(A(u))(s)\|_{[L_\infty(\rd)]^d}\,ds\\
    &\lesssim \|u_0\|_{[L_d(\rd)]^d}+\mathrm{I}'(t)+\mathrm{II}'(t).
\end{align*}
From Proposition \ref{prop:rp}, Lemma \ref{t5.1} and the H\"{o}lder inequality, we deduce that
\begin{align*}
    \mathrm{I}'(t):&=\sqrt t\int_0^{(1-\varepsilon)t}\lf\|H(t-s)\mathbb P(A(u))(s)\r\|_{[L_\infty(\rd)]^d}\,ds\\
    &\lesssim \sqrt t\int_0^{(1-\varepsilon)t}(t-s)^{-\frac 12-\frac{d}{d+2}}\|u(s)\|_{[L_{d+2}(\rd)]^d}^2\,ds\\
    &\le M^2\sqrt t\int_0^{(1-\varepsilon)t}(t-s)^{-\frac 12-\frac{d}{d+2}}s^{-\frac{2}{d+2}}\,ds\\
    &\le \varepsilon^{-\frac 32}M^2
\end{align*}
and
\begin{align*}
    \mathrm{II}'(t):&=\sqrt t\int_{(1-\varepsilon)t}^t\|H(t-s)\mathbb P(A(u))(s)\|_{[L_\infty(\rd)]^d}\,ds\\
    &\lesssim \sqrt t\int_{(1-\varepsilon)t}^t(t-s)^{-\frac 12-\frac{d}{2d+4}}\|u(s)\|_{[L_{d+2}(\rd)]^d}\|u(s)\|_{[L_{\fz}(\rd)]^d}\,ds\\
    &\le M\sqrt t\int_{(1-\varepsilon)t}^t(t-s)^{-\frac 12-\frac{d}{2d+4}}s^{-\frac{1}{d+2}}\|u(s)\|_{[L_{\fz}(\rd)]^d}\,ds\\
    &\le C M\int_{(1-\varepsilon)t}^t(t-s)^{-\frac 12-\frac{d}{2d+4}}s^{-\frac{1}{d+2}}\sqrt s\|u(s)\|_{[L_{\fz}(\rd)]^d}\,ds,
\end{align*}
which, combined with the estimates of $\mathrm{I}'(t)$ and $\mathrm{II}'(t)$, implies that, for some constant $C>0$,
\begin{align*}
    \sqrt t\|u(t)\|_{[L_\infty(\rd)]^d}\le C\|u_0\|_{[L_d(\rd)]^d}+C\varepsilon^{-\frac 32}M^2+CM\int_{(1-\varepsilon)t}^t(t-s)^{-\frac 12-\frac{d}{2d+4}}s^{-\frac{1}{d+2}}\sqrt s\|u(s)\|_{[L_{\fz}(\rd)]^d}\,ds.
\end{align*}
Define $$\psi(t):=\sup_{\eta\le w\le t}\sqrt w\|u(w)\|_{[L_{\fz}(\rd)]^d}$$ with $\eta=(1-\varepsilon)\lambda$. This yields that
\begin{align*}
    \psi(t)&\le \lf(\sup_{0\le w\le\lambda}+\sup_{\lambda\le w\le t}\r)\sqrt w\|u(w)\|_{[L_{\fz}(\rd)]^d}\\
    &\le N+C\|u_0\|_{[L_d(\rd)]^d}+C\varepsilon^{-\frac 32}M^2+CM\psi(t)\int_{1-\varepsilon}^1(1-s)^{-\frac 12-\frac{d}{2d+4}}s^{-\frac{1}{d+2}}\,ds.
\end{align*}
Since the integral $\int_{0}^1(1-s)^{-\frac 12-\frac{d}{2d+4}}s^{-\frac{1}{d+2}}\,ds$ converges, for sufficiently small $\varepsilon$, we have $$\int_{1-\varepsilon}^1(1-s)^{-\frac 12-\frac{d}{2d+4}}s^{-\frac{1}{d+2}}\,ds\le \frac 1{2CM}.$$ Therefore, for any $t\in (\lambda,T_{u_0,\fz}),$
\begin{align*}
    \sqrt t\|u(t)\|_{[L_{\fz}(\rd)]^d}\le\psi(t)\le 2N+2C\|u_0\|_{[L_d(\rd)]^d}+2C\varepsilon^{-\frac 32}M^2<\fz,
\end{align*}
contradicting  $$\lim_{t\nearrow T_{u_0,\fz}}\|u(t)\|_{[L_{\fz}(\rd)]^d}=\fz.$$

We can then deduce that $T_{u_0,p}=T_{u_0}$ for any $d<p\le \fz$ by interpolation. To prove this, notice that, for $d+2\le p\le \fz$ and any $0<T<T_{u_0},$ \begin{align*}
    \lf\|(\cdot)^{\frac 12-\frac d{2p}}u(\cdot)\r\|_{L_{\fz}([0,T];[L_p(\rd)]^d)}\le \lf\|(\cdot)^{\frac 1{d+2}}u(\cdot)\r\|_{L_{\fz}([0,T];[L_{d+2}(\rd)]^d)}^{\frac{d+2}{p}}\lf\|(\cdot)^{\frac 12}u(\cdot)\r\|_{L_{\fz}([0,T];[L_\fz(\rd)]^d)}^{1-\frac{d+2}p}<\fz.
\end{align*} A similar argument with $p=d+2$ shows that if we assume $T_{u_0,p}<T_{u_0}$, then we have
\begin{align}\label{blowup}
   \|u(t)\|_{[L_{p}(\rd)]^d}\ge \frac{\delta}{4\lf(T_{u_0,p}-t\r)^{\frac 12-\frac{1}{2p}}}, \ \ \ \forall d<p<\infty.
\end{align}
Taking $t\rightarrow T_{u_0,p}$ leads to a contradiction.
For every $p\in(d,d+2),$ by an analogous reasoning via interpolation between $L_{\fz}([0,T];[L_{d}(\rd)]^d)$ and $L_{\fz}([0,T];[L_{d+2}(\rd)]^d),$ then we have $T_{u_0,p}=T_{u_0}$.

\textbf{Step 6.} {\bf(The solution $u$ is actually infinitely smooth)} In this step, we will prove that the solution $u$ is infinitely smooth.
That is, for all $n\in\mathbb{N}$ and $0<T_1<T_2<T_{u_0}$, we have
$$\partial_t^{n}u\in \bigcap_{k\in\bN}L_{\fz}([T_1,T_2];[H_d^k(\rd)]_0^d).$$
Since $B^{k/2}_{d,\infty}(\rd)\subseteq H_d^{(k-1)/2}(\rd)$ (see Lemma \ref{bec}),
then we only need to show that
$$\partial_t^{n}u\in \bigcap_{k\in\nn}L_\infty([T_1,T_2];[B^{k/2}_{d,\infty}(\rd)]_0^d)),
\ \ \mathrm{for\ all} \  0<T_1<T_2<T_{u_0}, n\in\mathbb{N}.$$
Firstly, we show the case $n=0$.
 To do this, we use the method of induction for $k$. Since $L_d(\rd)\subset B^0_{d,\infty}(\rd)$ (see Lemma \ref{l2.81}), we see that, when $k=0$, the result is true. Now let us show the desired result for $k+1$ from the induction assumption on the one for $k\ge 0$. That is, we assume $$u\in L_\infty([T_1,T_2];[ B^{k/2}_{d,\infty}(\rd)]_0^d),
 \ \ \ \mathrm{for\ all} \ 0<T_1<T_2<T_{u_0},$$
then we show that
 $$u\in L_\infty([T_1,T_2];[B^{(k+1)/2}_{d,\infty}(\rd)]_0^d), \ \ \ \mathrm{for\ all} \ 0<T_1<T_2<T_{u_0}.$$
Since $u\in L_{\fz}([0,T_{u_0});[L_d(\rd)]_0^d)$,
by Lemma \ref{hom and inhom}, it suffices to show that
\begin{align*}
u\in L_\infty([T_1,T_2];[\dot{B}^{(k+1)/2}_{d,\infty}(\rd)]^d_0),\ \ \mathrm{for\ all} \  0<T_1<T_2<T_{u_0}.
\end{align*}
  For $0<T_0<T_1<t<T_2<T_{u_0}$, $u$ satisfies the following equation,
$$u(t)=H(t-T_0)u(T_0)-\int^t_{T_0}H(t-s)\mathbb P(A(u))(s)\,ds.$$
Therefore, we have
\begin{align*}
\|u(t)\|_{[\dot B^{(k+1)/2}_{d,\infty}(\rd)]^d}
&\leq\lf\|H(t-T_0)u(T_0)\r\|_{[\dot B^{(k+1)/2}_{d,\infty}(\rd)]^d}
+\lf\|\int^t_{T_0}H(t-s)\mathbb P(A(u)(s))\,ds\r\|_{[\dot B^{(k+1)/2}_{d,\infty}(\rd)]^d}\\
&=:\mathrm{I}_1(t)+\mathrm{I}_2(t).
\end{align*}
From the definition of Besov spaces, and Lemma \ref{lem:psi}, it is easy to check that
\begin{align*}
\mathrm{I}_1(t)&=\|H(t-T_0)u(T_0)\|_{[\dot B^{(k+1)/2}_{d,\infty}(\rd)]^d}\\
&=\sup_{j\in\zz} 2^{(k+1)j/2}\lf\|\dot{\triangle}_jH(t-T_0)u(T_0)\r\|
_{[L_{d}(\rd)]^d}\\
&\leq\sup_{j\in\zz} 2^{(k+1)j/2}\lf\|\mathcal{F}^{-1}(e^{-(t-T_0)|\cdot|^2}\varphi_j)\r\|_{L_1(\rr^d)}\lf\|u(T_0)\r\|
_{[L_{d}(\rd)]^d}\\
&\lesssim \sup_{j\in\zz} 2^{(k+1)j/2} e^{-c(t-T_0)2^{2j}}\|u(T_0)\|_{[L_{d}(\rd)]^d}
\\&\lesssim \max\lf\{(t-T_0)^{-\frac{k+1}{4}},1\r\}\cdot\|u(T_0)\|_{[L_{d}(\rd)]^d}.
\end{align*}
For the nonlinear term $\mathrm{I}_2(t)$,  applying Lemmas \ref{lem:psi} and \ref{map:P1}, we find that
\begin{align*}
\mathrm{I}_2(t)&=\lf\|\int^t_{T_0}H(t-s)\mathbb P (A(u))(s)\,ds\r\|_{[\dot B^{(k+1)/2}_{d,\infty}(\rd)]^d}\\
&\leq\int^t_{T_0}\lf\|H(t-s)\mathbb P(A(u))(s)\r\|_{[\dot B^{(k+1)/2}_{d,\infty}(\rd)]^d}ds\\
&\lesssim\int^{t}_{T_0}\max\lf\{\frac 1{(t-s)^{\frac34}},1\r\}\lf\|\mathbb P (A(u))(s)\r\|_{[\dot B^{k/2-1}_{d,\infty}(\rd)]^{d}}\,ds\\
&\lesssim\int^{t}_{T_0}\max\lf\{\frac 1{(t-s)^{\frac34}},1\r\}\lf\|u(s)\otimes u(s)\r\|_{[\dot B^{k/2}_{d,\infty}(\rd)]^{d^2}}\,ds\\
&\lesssim \lf\|u\otimes u\r\|_{L_\infty([T_0,t];[\dot B^{k/2}_{d,\infty}(\rd)]^{d^2})}\\
&\lesssim C(u),
\end{align*}
where \begin{eqnarray*}C(u):=
\lf\{\begin{array}{ll}
\|u\|_{L_\infty([T_0,t];[B^{k/2}_{d,\infty}(\rd)]^d)}\|u\|_{L_\infty([T_0,t];[L_{\infty}(\rd)]^d)},  \ \ \mathrm{when} \ \ k\ge 1;\\[5pt]
\|u\|_{L_\infty([T_0,t];[L_{d}(\rd)]^d)}\|u\|_{L_\infty([T_0,t];[L_{\infty}(\rd)]^d)}, \ \ \ \ \mathrm{when} \ \  k=0.
\end{array}\r.
\end{eqnarray*} The last inequality follows from the fact that $L_d(\rd)\subset B_{d,\infty}^0(\rd)$ and Lemma \ref{McDonald}.

Next, we show that $u$ is indeed the strong solution of the Navier-Stokes equation \eqref{ns1}, that is, $$\lim_{h\to 0}\lf\|\frac{u(t+h)-u(t)}{h}-\Delta_{\theta} u(t)-\mathbb{P}(A(u))(t)\r\|_{[L_d(\rd)]^d}=0$$ for any $t\in (0,T_{u_0}).$

For any $0<T_1<t<T_2<T_{u_0},$ recall the representation:
$$u(t)=H(t-T_1)u(T_1)-\int_{T_1}^t H(t-s)\mathbb{P}(A(u))(s)\, ds.$$ For
 simplicity, set $f(s):=- \mathbb{P}(A(u))(s).$ Then, for any $0<h<T_2-t$, we have
\begin{align*}
    \frac{u(t+h)-u(t)}{h}
    =&\frac{H(t-T_1+h)u(T_1)-H(t-T_1)u(T_1)}{h}+\int_{T_1}^t \frac{H(t-s+h)f(s)-H(t-s)f(s)}{h}\, ds \\
  &+\frac 1h\int_t^{t+h} H(t-s+h)f(s)\, ds
    \\
    =&\frac{H(h)-1}{h}\lf(H(t-T_1)u(T_1)-\int_{T_1}^t H(t-s)\mathbb{P}f(s)\, ds\r)+\frac 1h\int_t^{t+h}H(t-s+h)f(s)\, ds
    \\
    =&\frac{H(h)-1}{h}u(t)+\frac 1h\int_t^{t+h}H(t-s+h)f(s)\, ds.
\end{align*}
For the first term, since $u\in \bigcap_{k\in\bN}L_{\fz}([T_1,T_2];[H_d^k(\rd)]_0^d),$ by Lemma \ref{heat} \rm{(iii)}, we know
\begin{align*}
    \lim_{h\to 0}\lf\|\frac{H(h)-1}{h}u(t)-\Delta_\theta u(t)\r\|_{[L_d(\rd)]^d}=0.
\end{align*}
For the second term, we can write
\begin{align*}
    \frac 1h\int_t^{t+h}H(t-s+h)f(s)\, ds&=\frac 1h\lf(\int_t^{t+h}H(t-s+h)f(s)-f(t)\, ds\r)+f(t).
\end{align*}
Since $f\in \bigcap_{k\in\bN}L_{\fz}([T_1,T_2];[H_d^k(\rd)]_0^d)$ and $f\in C((0,T_{u_0});[L_d(\rd)]_0^d),$ by Lemma \ref{heat}, we have
\begin{align*}
   &\lf\|\frac 1h\int_t^{t+h}H(t-s+h)f(s)-f(t)\, ds\r\|_{[L_d(\rd)]^d}\\
   \le& \frac 1h\int_t^{t+h}\big\|H(t-s+h)f(s)-f(t)\big\|_{[L_d(\rd)]^d}\, ds
   \\ \le& \frac 1h\int_t^{t+h} \big\|H(t-s+h)f(s)-H(t-s+h)f(t)\big\|_{[L_d(\rd)]^d}\, ds\\
   &\ \ \ +\frac 1h\int_t^{t+h}\big\|H(t-s+h)f(t)-H(h)f(t)\big\|_{[L_d(\rd)]^d}\, ds\\
   &\ \ \ +\frac 1h\int_t^{t+h}\big\|H(h)f(t)-f(t)\big\|_{[L_d(\rd)]^d}\, ds\\
   \le& \frac{1}{h}\int_t^{t+h}\big\|f(s)-f(t)\big\|_{[L_d(\rd)]^d}\,ds
    +\frac{1}{h}\int_t^{t+h}\lf\|\Big(H(t-s)-\mathrm{Id}_{\rd}\Big)f(t)\r\|_{[L_d(\rd)]^d}\,ds\\
   &\ \ \ +\lf\|\lf(H(h)-\mathrm{Id}_{\rd}\r)f(t)\r\|_{[L_d(\rd)]^d}.
\end{align*}
Thus, we deduce
\begin{align*}
   &\lim_{h\rightarrow0}\lf\|\frac 1h\int_t^{t+h}H(t-s+h)f(s)-f(t)\, ds\r\|_{[L_d(\rd)]^d}\\
   \le& \lim_{h\rightarrow0}\frac{1}{h}\int_t^{t+h}\big\|f(s)-f(t)\big\|_{[L_d(\rd)]^d}ds
+\lim_{h\rightarrow0}\frac{1}{h}\int_t^{t+h}\lf\|\Big(H(t-s)-\mathrm{Id}_{\rd}\Big)f(t)
\r\|_{[L_d(\rd)]^d}ds\\
&\quad+\lim_{h\rightarrow0}\lf\|\lf(H(h)-\mathrm{Id}_{\rd}\r)f(t)\r\|_{[L_d(\rd)]^d}=0.
\end{align*}
Altogether, we have proven that $\partial_t u\in L_{\fz}([T_1,T_2];[L_d(\rd)]_0^d)$. Thus, we obtain
 \begin{equation*}
   \partial_t u\in L_{\fz}([T_1,T_2];[L_d(\rd)]_0^d),
 \ \ \ \mathrm{for\ all} \ 0<T_1<T_2<T_{u_0},
 \end{equation*}
and $$\partial_t u(t)=\Delta_{\theta} u(t)+ \mathbb{P}(A(u))(t).$$

Now, we employ induction to show that, for all $n\in \mathbb{N}, 0<T_1<T_2<T_{u_0}$, we have $$\partial_t^{n}u\in \bigcap_{k\in\nn}L_\infty([T_1,T_2];[B^{k/2}_{d,\infty}(\rd)]_0^d).$$
We assume that, for a fixed $n_0\in\mathbb{N}$, we have $\partial_t^m u\in \bigcap_{k\in\bN}L_{\fz}([T_1,T_2];[B^{k/2}_{d,\infty}(\rd)]_0^d)$ for all $m\le n_0$ and $0<T_1<T_2<T_{u_0}.$  Note that $\bigcap_{k\in\bN}L_{\fz}([T_1,T_2];[B^{k/2}_{d,\infty}(\rd)]^d)\subset L_{\fz}([T_1,T_2];[L_\infty(\rd)]^d)$. Thus, for all $m\le n_0$ and $0<T_1<T_2<T_{u_0},$ we have
$$\partial_t^{m}u\in L_{\fz}([T_1,T_2];[L_\infty(\rd)]_0^d).$$
Hence, combining the induction assumption, Lemma \ref{bd}, Lemma \ref{map:P1}, the approach utilized in the proof for the scenario where $n=0$ and the fact
\begin{align*}
    \partial_t^{n_0+1}u&=\partial_t^{n_0}(\Delta_{\theta} u+ \mathbb{P}(A(u)))
    =\Delta_{\theta} \partial_t^n u+ \sum_{j=0}^{n_0}\binom{n_0}{j}\mathbb{P}\mathbf{div}\, (\partial_t^ju\otimes \partial_t^{n_0-j}u),
\end{align*}
we obtain
 $$\partial_t^{n_0+1} u\in \bigcap_{k\in\nn}L_\infty([T_1,T_2];[B^{k/2}_{d,\infty}(\rd)]_0^d),
\ \ \mathrm{for\ all} \  0<T_1<T_2<T_{u_0}.$$
Therefore, we complete the proof of {\bf Step 6}.

\textbf{Step 7.} {(\bf $\|u_0\|_{[L_d(\rd)]^d}$ is small enough $\Longrightarrow T=\fz$)}
If $\|u_0\|_{[L_d(\rd)]^d}$ is sufficiently small, by the estimate \eqref{timespace1}, we can take $T=\infty$ in the definition of $ \mathcal{D}_\delta$. Consequently, the preceding arguments remain valid under this condition.
This completes the proof of part (i) of Theorem \ref{t3.1}.

\bigskip

{\bf Part (ii):}
Next, we proceed to prove part (ii) of Theorem \ref{t3.1}. Suppose that $u_0\in[L_2(\rd)]_0^d\cap[L_d(\rd)]_0^d$.
We will prove that the solution in (i) satisfies $u\in{L_2([0,T];[\dot{H}^1(\rd)]_0^d)}$ for some $T>0$.
Introduce a constant $\sigma>0$, whose value will be specified later. Throughout the subsequent steps, we will frequently invoke the Picard contraction principle.
 Define the metric space $\mathcal{D}_\sigma$ as follows:
\begin{align}\label{metric space}
\mathcal{D}_\sigma
:=&\lf\{u\in C([0,T]; [L_d(\rd)]^d)\cap{L_{d+2}([0,T];[L_{d+2}(\rd)]^d)}\r.\\
&\ \ \ \lf.\cap L_{2+4/d}([0,T]; [L_{2+4/d}(\rd)]^d)\cap L_2([0,T]; [\dot{H}^1(\rd)]^d):\|u\|_{\mathcal{D}_{\sigma}}\leq \sigma\r\},\nonumber
\end{align}
where the distance and the norm are given by
\begin{align}\label{metric}
d(u,v):=\|u-v\|_{L_{2+4/d}([0,T]; [L_{2+4/d}(\rd)]^d)}+\|u-v\|_{L_2([0,T]; [\dot{H}^1(\rd)]^d)}+\|u-v\|_{L_{d+2}([0,T]; [L_{d+2}(\rd)]^d)},
\end{align}
and $\|u\|_{\mathcal{D}_{\sigma}}:=\|u\|_{L_{2+4/d}([0,T]; [L_{2+4/d}(\rd)]^d)}
+\|u\|_{L_2([0,T];[\dot{H}^1(\rd)]^d)}+\|u\|_{L_{d+2}([0,T]; [L_{d+2}(\rd)]^d)}$.
Now, we show that the metric space $(\mathcal D_\sigma,\,d)$ is complete. Indeed,
let $\{u_j\}_j$ be a Cauchy sequence in $(\mathcal D_\sigma,\,d)$. Then, there exist $u\in L_{2+4/d}([0,T]; [L_{2+4/d}(\rd)]^d)$ and $\widetilde{u}\in L_2([0,T]; [\dot{H}^1(\rd)]^d)$ such that
\begin{align*}
u_j\rightarrow u \ \ \mathrm{in} \ \ L_{2+4/d}([0,T]; [L_{2+4/d}(\rd)]^d),\ \
\mathrm{as} \ \ j\rightarrow\fz,
\end{align*}
\begin{align*}
u_j\rightarrow \widetilde{u} \ \ \mathrm{in} \ \ L_2([0,T]; [\dot{H}^1(\rd)]^d),\ \
\mathrm{as}  \ j\rightarrow\fz,
\end{align*}
and
\begin{align*}
u_j\rightarrow \overline{u} \ \ \mathrm{in} \ \ L_{d+2}([0,T]; [L_{d+2}(\rd)]^d),\ \
\mathrm{as} \ \ j\rightarrow\fz.
\end{align*}
 It is straightforward to verify that $u=\widetilde{u}=\overline{u}$.
 Consequently, $u\in (\mathcal{D}_{\sigma},\,d)$ and $d(u_j,\,u)\rightarrow0$ as $j\rightarrow\fz$, which imply that $(\mathcal{D}_\sigma,\,d)$ is complete.

Now, let we consider the map $\mathfrak{M}$ defined as:
\begin{align}\label{solution map}
\mathfrak{M}: u\mapsto Hu_0-\mathcal H\mathbb P\lf(A(u)\r).
\end{align}
We claim that $\mathfrak M$ is a strictly contractive mapping on $(\mathcal{ D}_{\sigma},d)$ with $\sigma>0$ such that $C\sigma\leq\frac18$, where $C$ is the  constant appearing in the following. Initially, for $u\in \mathcal{D}_\sigma$, by invoking the triangle inequality alongside the inequalities \eqref{timespace30} and \eqref{timespace40}, we derive
\begin{align}\label{es1}
\|\mathfrak Mu\|_{L_{2+4/d}([0,T]; [L_{2+4/d}(\rd)]^d)}&\leq \|Hu_0\|_{L_{2+4/d}([0,T]; [L_{2+4/d}(\rd)]^d)}+\lf\| \mathcal{H}\mathbb{P}(A(u))\r\|_{L_{2+4/d}([0,T]; [L_{2+4/d}(\rd)]^d)}\\ \nonumber
&\leq \|Hu_0\|_{L_{2+4/d}([0,T]; [L_{2+4/d}(\rd)]^d)}+C\lf\|\mathbb P(A(u))\r\|_{L_{(2+4/d)'}([0,T]; [L_{(2+4/d)'}(\rd)]^d)}\nonumber,
\end{align}
\begin{align}\label{es2}
\|\mathfrak Mu\|_{L_2([0,T]; [\dot{H}^1(\rd)]^d)}
&\leq \| Hu_0\|_{L_2([0,T]; [\dot{H}^1(\rd)]^d)}+\lf\|\mathcal{H}\mathbb P(A(u))\r\|_{L_2([0,T]; [\dot{H}^1(\rd)]^d)}\\ \nonumber
&\leq \| Hu_0\|_{L_2([0,T]; [\dot{H}^1(\rd)]^d)}+C\lf\|\mathbb P(A(u))\r\|_{L_{(2+4/d)'}([0,T]; [L_{(2+4/d)'}(\rd)]^d)} \nonumber
\end{align}
and
\begin{align}\label{es3}
\|\mathfrak Mu\|_{L_{d+2}([0,T]; [L_{d+2}(\rd)]^d)}
\leq\|Hu_0\|_{L_{d+2}([0,T]; [L_{d+2}(\rd)]^d)}
+\|\mathcal{H}\mathbb P(A(u))\|_{L_{d+2}([0,T]; [L_{d+2}(\rd)]^d)}.
\end{align}
Applying  Lemma \ref{t5.1} and the H\"older inequality, we have
\begin{align}\label{nonlinear11}
\lf\|\mathbb P\lf(A(u)\r)\r\|_{L_{(2+4/d)'}([0,T]; [L_{(2+4/d)'}(\rd)]^d)}
\lesssim&\lf\|A(u)\r\|_{L_{(2+4/d)'}([0,T]; [L_{(2+4/d)'}(\rd)]^d)}\\ \nonumber
\lesssim &\|u\|_{L_{d+2}([0,T];[L_{d+2}(\rd)]^d)}\| u\|_{L_{2}([0,T]; [\dot{H}^1(\rd)]^d)}\\
\leq &C\sigma^2.\nonumber
\end{align}
From the inequality \eqref{timespace3},  Lemma \ref{t5.1} and the H\"older inequality, we conclude that
\begin{align}\label{nonlinear2}
\lf\|\mathcal{H}\mathbb P(A(u))\r\|_{L_{d+2}([0,T]; [L_{d+2}(\rd)]^d)}&=\lf\|\mathbb P(A(u))\r\|_{L_{(d+2)/2}([0,T]; [L_{(d+2)/2}(\rd)]^d)}\\\nonumber
&\leq C\|u\otimes u\|_{L_{(d+2)/2}([0,T]; [L_{(d+2)/2}(\rd)]^{d^2})}\\\nonumber
&\leq C\|u\|_{L_{d+2}([0,T]; [L_{d+2}(\rd)]^{d})}^2\\
&\leq C\sigma^2.\nonumber
\end{align}
Furthermore, by the inequalities \eqref{timespace10}, \eqref{timespace20} and \eqref{timespace1}, there exists a $T=T_{\sigma}>0$, which might depend on $\sigma$, such that
\begin{align}\label{es4}
\|Hu_0\|_{L_{2+4/d}([0,T]; [L_{2+4/d}(\rd)]^d)}+\|Hu_0\|_{L_2([0,T]; [\dot{H}^1(\rd)]^d)}+\|Hu_0\|_{L_{d+2}([0,T]; [L_{d+2}(\rd)]^d)}\leq \frac \sigma4.
\end{align}
From these estimates \eqref{es1}--\eqref{es4}, we deduce
\begin{align}\label{close}
&\|\mathfrak Mu\|_{L_{2+4/d}([0,T]; [L_{2+4/d}(\rd)]^d)}+\|\mathfrak Mu\|_{L_2([0,T]; [\dot{H}^1(\rd)]^d)}+\|\mathfrak M u\|_{L_{d+2}([0,T]; [L_{d+2}(\rd)]^d)}\\
\leq& \frac \sigma4+2C^2\sigma^2+C\sigma^2\leq \sigma,\nonumber
\end{align}
which implies that $\mathfrak Mu\in \mathcal{D}_\sigma$ for any $u\in \mathcal{D}_\sigma$.

The strict contraction of $\mathfrak M$ can be demonstrated in a manner analogous to the proof of part (i). Specifically, let $u,v\in \mathcal{D}_\sigma$, by Corollaries \ref{cor:ts} and \ref{cor:Hts} and the boundedness of $\mathbb P$, we get
\begin{align}\label{contraction}
&d(\mathfrak Mu,\,\mathfrak Mv)\\
\nonumber=&\|\mathfrak Mu-\mathfrak Mv\|_{L_{2+4/d}([0,T]; [L_2(\rd)]^d)}+\|\mathfrak Mu-\mathfrak Mv\|_{L_{2}([0,T]; [\dot{H}^1(\rd)]^d)}+\|\mathfrak M u-\mathfrak M v\|_{L_{d+2}([0,T]; [L_{d+2}(\rd)]^d)}\\
\nonumber\leq& C\sigma\Big(\|u-v\|_{L_{2+4/d}([0,T]; [L_{2+4/d}(\rd)]^d)}+\|u-v\|_{L_{2}([0,T];  [\dot{H}^1(\rd)]^d)}+\|u-v\|_{L_{d+2}([0,T]; [L_{d+2}(\rd)]^d)}\Big)\\
\nonumber\leq& \frac12\Big(\|u-v\|_{L_{2+4/d}([0,T]; [L_{2+4/d}(\rd)]^d)}+\|u-v\|_{L_{2}([0,T];  [\dot{H}^1(\rd)]^d)}+\|u-v\|_{L_{d+2}([0,T]; [L_{d+2}(\rd)]^d)}\Big).
\end{align}
 Therefore,
\begin{align}
d(\mathfrak Mu,\mathfrak Mv)\leq \frac12d(u,v).
\end{align}
According to the Picard contraction mapping principle, there exists a unique $u\in \mathcal{D}_\sigma$ such that
\begin{align}\label{mild solution}
u=Hu_0-\mathcal H\mathbb P\lf(A(u)\r).
\end{align}
Furthermore, utilizing the inequalities \eqref{timespace00}, \eqref{timespace40}, and \eqref{nonlinear11}, we derive the following:
\begin{eqnarray*}
 \|u\|_{L_{\infty}([0,T]; [L_{2}(\rd)]^d)}&=&\|\mathfrak Mu\|_{L_{\infty}([0,T]; [L_{2}(\rd)]^d)}\\
 &\leq& \|Hu_0\|_{L_{\infty}([0,T]; [L_{2}(\rd)]^d)}+\lf\| \mathcal{H}\mathbb{P}(A(u))\r\|_{L_{\infty}([0,T]; [L_{2}(\rd)]^d)}\\
 &\overset{\eqref{timespace00}+ \eqref{timespace40}}\lesssim&\|u_0\|_{[L_{2}(\rd)]^d}+\lf\|\mathbb P\lf(A(u)\r)\r\|_{L_{(2+4/d)'}([0,T]; [L_{(2+4/d)'}(\rd)]^d)}\\
&\overset{\eqref{nonlinear11}}\lesssim&\|u_0\|_{[L_{2}(\rd)]^d}+\|u\|_{L_{d+2}([0,T]; [L_{d+2}(\rd)]^d)}\| u\|_{L_{2}([0,T]; [\dot{H}^1(\rd)]^d)}.
\end{eqnarray*}
Consequently, we confirm that $u\in L_{\infty}([0,T]; [L_{2}(\rd)]^d)$. Additionally, as demonstrated in {\bf Part (i)}, we can show $u\in C([0,T]; [L_{2}(\rd)]^d)$.
This concludes the proof of {\bf Part (ii)} of Theorem \ref{t3.1}.

\end{proof}

\section{The proof of Theorem \ref{t3.1x}}\label{s7}

This section is devoted to the proof of Theorem \ref{t3.1x}, whose proof is similar to
Theorem \ref{t3.1}, thus we only briefly outline the main ideas of the proof.
\begin{proof}[Proof of Theorem \ref{t3.1x}]
For part (i), by substituting the nonlinear term $A(u)$ with $S(u)$ and following the same process as in the proof of Theorem \ref{t3.1}(i), we can conclude part (i). For brevity, we omit the details here.

Next, we proceed to demonstrate part (ii). Drawing parallels from the proof of Theorem \ref{t3.1} (ii), for the equation \eqref{ens1}, we establish the existence of a maximal time $T_{u_0}$ and a unique smooth solution $u$ such that
$$ u\in C([0,T_{u_0});[L_2(\rr_{\theta}^d)]^d_0)\cap L_{2}([0,T_{u_0});[\dot{H}^1(\rr_{\theta}^d)]^d_0).$$
Given that $u$ is smooth and $u_0^*=u_0$, we can derive the conservation law \eqref{energy} directly from the mild form \eqref{mild solution}
This concludes the proof of part (ii).

Finally, we prove part (iii). Specifically, when $d=2$, $u_0^*=u_0$, the maximal time $T_{u_0}=\infty$.
Assume, for contradiction, that $T_{u_0}<\infty$. We claim that
\begin{align}\label{contradiction}
\|u\|_{C([0,T_{u_0}); [L_2(\ri)]^2)}+\|u\|_{L_2([0,T_{u_0});[\dot{H}^1(\ri)]^2)}=\infty;
\end{align}
however, this would contradict the conservation law of $u$,
\begin{align*}
\frac{1}{2}\|u(t)\|^2_{[L_2(\ri)]^2}+\int^t_0\|\nabla_{\theta} u(s)\|^2_{[L_2(\ri)]^4}ds=\frac{1}{2}\|u_0\|^2_{[L_2(\ri)]^2}, \;\forall\;0<t<T_{u_0},
\end{align*}
and thus $T_{u_0}$ must be infinite.

We show the claim \eqref{contradiction} now. Suppose instead that\begin{align*}
\|u\|_{C([0,T_{u_0}); [L_2(\ri)]^2)}+\|u\|_{L_2([0,T_{u_0});[\dot{H}^1(\ri)]^2)}<\infty.
\end{align*}
By applying the Gagliardo-Nirenberg inequality (see Lemma \ref{gn1}) and the H\"older inequality, we deduce that
\begin{align*}
\|u\|_{L_4([0,T_{u_0});[L_4(\ri)]^2)}\lesssim \|u\|^{\frac12}_{L_\infty([0,T_{u_0});[L_2(\ri)]^2)}\|u\|^{\frac12}_{L_2([0,T_{u_0});[\dot{H}^1(\ri)]^2)}<\infty.
\end{align*}
Following the proof in Theorem \ref{t3.1}, we can extend $T_{u_0}$ to a larger time, which contradicts the finiteness of $T_{u_0}$.
Thus, we have completed the proof of Theorem \ref{t3.1x}.

\end{proof}

\section{The proofs of Theorem \ref{semiclassical} and Theorem \ref{semi2}}\label{s8}

In this section, we shall prove Theorems \ref{semiclassical} and \ref{semi2}. First of all, we recall some basic definitions and notations. For any $f, g \in \mathcal{S}(\mathbb{R}^d) $ and $\xi\in\rr^d$, we define
\begin{equation*}
  f \ast_\theta g(\xi) := \int_{\mathbb{R}^d}e^{\frac{\mathrm{i}}{2}(\xi, \theta\eta)} f(\xi-\eta)g(\eta) \,d\eta.
\end{equation*}
Then we have
\begin{align*}
 U_\theta(f)U_\theta(g)
 =&\int_{\rr^d}f(s)\lambda_{\theta}(s)\,ds\cdot\int_{\rr^d}g(t)\lambda_{\theta}(t)\,dt\\
 =&\iint_{\rr^{2d}}f(s)g(t)e^{\frac{\mathrm{i}}{2}(s,\theta t)}\lambda_{\theta}(s+t)\,dtds\\
 =&\iint_{\rr^{2d}}f(s-t)g(t)e^{\frac{\mathrm{i}}{2}(s,\theta t)}\lambda_{\theta}(s)\,dtds\\
 =&\iint_{\rr^{2d}}e^{\frac{\mathrm{i}}{2}(s,\theta t)}f(s-t)g(t)\,dt\lambda_{\theta}(s)\,ds\\
  =&U_\theta(f \ast_\theta g),
 \end{align*}
 and
\begin{align}\label{Moyal product}
    \mathcal{F}^{-1}(\hat{f} \ast_\theta \hat{g})(\xi) = \frac{1}{(2\pi)^{d}}\int_{\mathbb{R}^d}\int_{\mathbb{R}^d} f\left(\xi + \frac{\theta}{2}v\right)g(\xi - w)e^{\mathrm{i}(v, w)} dv dw, \quad f, g \in \mathcal{S}(\mathbb{R}^d).
\end{align}
The transformation \( \star_\theta:(f, g) \mapsto \mathcal{F}^{-1}(\hat{f} \ast_\theta \hat{g}) \) is known as the Moyal product. For further details about Moyal product, we refer to \cite{GJM2022,Rieffel}.

By Lemmas \ref{se1} and \ref{l2.a1}, we obtain the following result:

\begin{lemma}\label{Wigner bounds}
    If $ p > 2 $ and $ s > \frac{d}{2} - \frac{d}{p} $, then for $ x \in H^s(\mathbb{R}^d_{\theta}) $, we have
    $$
    \left\| \mathcal{F}^{-1} \circ U_\theta^{-1}(x) \right\|_{L_p(\mathbb{R}^d)} \lesssim \|x\|_{H^s(\mathbb{R}^d_{\theta})}.
    $$
\end{lemma}

\begin{prop}\label{star}
    The Moyal product extends to a continuous bilinear mapping from \( L_2(\mathbb{R}^d) \times L_2(\mathbb{R}^d) \to \mathcal{S}'(\mathbb{R}^d) \). Furthermore, \( \mathcal{F}^{-1}(\hat{f} \ast_\theta \hat{g}) \to fg \) in the distributional sense as \( \theta \to 0 \).
\end{prop}

\begin{proof}
    By Lemmas \ref{l2.a1} and \ref{se1}, for \( f, g \in \mathcal{S}(\mathbb{R}^d) \) and \( s > \frac{d}{2} \), we have
    \begin{align*}
        \left\| \mathcal{F}^{-1}(\hat{f} \ast_\theta \hat{g}) \right\|_{H^{-s}(\mathbb{R}^d)} &= \left\| U_{\theta}(\hat{f} \ast_\theta \hat{g}) \right\|_{H^{-s}(\mathbb{R}^d_{\theta})} \\
        &\leq \left\| U_{\theta}(\hat{f} \ast_\theta \hat{g}) \right\|_{L_1(\mathbb{R}^d_{\theta})} \\
        &\leq \lf\| U_\theta(\hat{f})\r\|_{L_2(\mathbb{R}^d_{\theta})} \lf\| U_\theta(\hat{g}) \r\|_{L_2(\mathbb{R}^d_{\theta})} \\
        &= \| f \|_{L_2(\mathbb{R}^d)} \| g \|_{L_2(\mathbb{R}^d)}.
    \end{align*}
    Hence, the Moyal product extends to a bounded mapping from \( L_2(\mathbb{R}^d) \times L_2(\mathbb{R}^d) \) to \( H^{-s}(\mathbb{R}^d) \).

    Now we show that, for every \( \psi \in \mathcal{S}(\mathbb{R}^d) \),
    \begin{align}\label{esm}
        \lim_{\|\theta\| \to 0} \int_{\mathbb{R}^d} \mathcal{F}^{-1}(\hat{f} \ast_\theta \hat{g}) \psi = \int_{\mathbb{R}^d} fg\psi.
    \end{align}
    Let $$ \mathrm{I}_\theta := \int_{\mathbb{R}^d} \left( \mathcal{F}^{-1}(\hat{f} \ast_\theta \hat{g})(\xi) \psi(\xi) - f(\xi)g(\xi)\psi(\xi) \right) d\xi .$$
    By the equation \eqref{Moyal product}, we obtain
    \begin{align*}
        \mathrm{I}_\theta &=\frac{1}{(2\pi)^{d}} \iiint_{\mathbb{R}^{3d}} \left[ f\left(\xi + \frac{\theta}{2}v\right) - f(\xi) \right] g(\xi - w)\psi(\xi)e^{\mathrm{i}(v, w)} dv dw d\xi \\
        &=\frac{1}{(2\pi)^{d}} \iiint_{\mathbb{R}^{3d}} \left[ f\left(\xi + \frac{\theta}{2}v\right) - f(\xi) \right] \psi(\xi)e^{\mathrm{i}(v, \xi)} g(w)e^{-\mathrm{i}(v, w)} dw dv d\xi \\
        &=\frac{1}{(2\pi)^{d}} \iint_{\mathbb{R}^{2d}} f(\xi) \left[ \psi\left(\xi - \frac{\theta}{2}v\right) - \psi(\xi) \right] \hat{g}(v)e^{\mathrm{i}(v, \xi)} dv d\xi.
    \end{align*}
    By the dominated convergence theorem, we can deduce that
    $$\lim_{\|\theta\| \to 0}\mathrm{I}_\theta = 0,$$
    which implies \eqref{esm}. Therefore, we complete the proof of Proposition \ref{star}.
\end{proof}

Given that \( u \in C([0, T_{u_0}); [L_2(\rd)]^d_0) \) is a solution to equation \eqref{ens1}, it follows that $$\phi_\theta\triangleq \mathcal{F}^{-1} \circ U^{-1}_\theta(u) \in C([0, T_{u_0}); [L_2(\mathbb{R}^d)]^d_0)$$
 is a solution to
\begin{eqnarray}\label{ens2}
\begin{cases}
\partial_t\phi - \Delta \phi + \frac12\mathbb{P}\,\mathbf{div}\,\lf(\left[ \phi \otimes_\theta \phi \right]+\left[ \phi \otimes_\theta \phi \right]^T\r) = 0;\\[5pt]
\mathbf{div}\, \phi = 0;\\[5pt]
\phi(0) = \mathcal{F}^{-1} \circ U
^{-1}_\theta(u_0),
\end{cases}
\end{eqnarray}
where \( \phi \otimes_\theta \phi: = \left( \mathcal{F}^{-1}(\hat{\phi_i} \ast_\theta \hat{\phi_j}) \right)_{1 \le i,j \le d} \).

We can see that \eqref{ens2} is precisely the symmetric quantization of the classical Navier-Stokes equation, where the pointwise product is replaced by the symmetric Moyal product. By the same reason, if  $$u \in C([0, T_{u_0}); [L_2(\rd)]^d_0) \cap L_2([0, T_{u_0}); [\dot{H}^1(\rd)]^d_0)$$
 satisfies the energy identity \eqref{energy11} for all \( T < T_{u_0} \), then
  $$\phi_\theta \in C([0, T_{u_0}); [L_2(\mathbb{R}^d)]^d_0) \cap L_2([0, T_{u_0}); [\dot{H}^1(\mathbb{R}^d)]^d_0) $$
   satisfies
\begin{align*}
    \frac{1}{2}\|\phi_\theta(T)\|^2_{[L_2(\mathbb{R}^d)]^d} + \int^T_0 \lf\|\nabla\phi_\theta(s)\r\|^2_{[L_2(\mathbb{R}^d)]^d} ds = \frac{1}{2}\|\phi(0)\|^2_{[L_2(\mathbb{R}^d)]^d}
\end{align*}
for all \( T < T_{u_0} \). Moreover, in that case we can deduce that $u \in L_\fz([0, T_{u_0}); [L_2(\rd)]^d_0)$ and then $\phi_\theta \in L_\fz([0, T_{u_0}); [L_2(\bR^d)]^d_0)$ by the energy identity.

Similarly, under the conditions of Theorem \ref{semi2}, we also have the asymmetric quantization of the classical Navier-Stokes equation:
\begin{eqnarray}\label{ens22}
\begin{cases}
\partial_t\phi - \Delta \phi + \mathbb{P}\,\mathbf{div}\,\left[ \phi \otimes_\theta \phi \right]= 0;\\[5pt]
\mathbf{div}\, \phi = 0;\\[5pt]
\phi(0) = \mathcal{F}^{-1} \circ U
^{-1}_\theta(u_0),
\end{cases}
\end{eqnarray}
which follows from \eqref{ens} in the same way as \eqref{ens2} from \eqref{ens1}.

\begin{remark}
    {\rm The equation \eqref{ens22} coincides with the noncommutative Navier-Stokes equation in \cite{BSS2014}, and \eqref{ens2} is just the symmetrization of \eqref{ens22}.}
\end{remark}

\bigskip

We first present the proof of Theorem \ref{semiclassical}.

\begin{proof}[Proof of Theorem \ref{semiclassical}]
The existence of $u_\theta$ is ensured by Theorem \ref{t3.1x}.
By the boundedness of the map \( \mathcal{F}^{-1} \circ U_\theta^{-1} \),  the family \( \{\phi_\theta\}_\theta \) forms a bounded subset of \( L_\infty([0,\fz);[L_2(\mathbb{R}^2)]_0^2) \) and thus contains a subsequence \( \{\phi_{\theta_n}\}_{n} \) such that as \( \|\theta_n\| \to 0 \), \( \phi_{\theta_n} \to \phi \) in the weak-\( \ast \) topology as \( n \to \infty \) for some $\phi\in L_\infty([0,\fz);[L_2(\mathbb{R}^2)]_0^2)$. Note that \( \phi_\theta \) is the solution to \eqref{ens2} with initial datum \( \phi_0 \). We will show that \( \phi \) is the solution to the classical Navier-Stokes equation with initial datum \( \phi_0 \).

To prove this, observe that for every \( \psi \in [C_{\mathrm{c}}^\infty([0,\fz) \times \mathbb{R}^2)]^2_0 \), we have $\mathbb P\psi=\psi.$ Hence by testing over $\psi$ on the both sides of \eqref{ens2} and integrating by part, the self-adjointness of $\mathbb P$ yields
\begin{align*}
\int_{[0,\fz)\times\mathbb{R}^2}\partial_t\phi_{\theta_n}\psi&=-\int_{\mathbb{R}^2} \phi_0 \psi(0)-\int_{[0,\fz)\times\mathbb{R}^2}\phi_{\theta_n}\partial_t\psi\\
&=\int_{[0,\fz)\times\mathbb{R}^2}  \Delta \phi_{\theta_n} \psi - \frac 12\mathbb P\mathbf{div}\lf(\left[ \phi_{\theta_n} \otimes_{\theta_n} \phi_{\theta_n} \right]+\left[ \phi_{\theta_n} \otimes_{\theta_n} \phi_{\theta_n} \right]^T\r)\psi \\
&=\int_{[0,\fz)\times\mathbb{R}^2} \phi_{\theta_n} \Delta \psi + \frac 12\lf(\left[ \phi_{\theta_n} \otimes_{\theta_n} \phi_{\theta_n} \right]+\left[ \phi_{\theta_n} \otimes_{\theta_n} \phi_{\theta_n} \right]^T\r) \nabla \psi,
\end{align*}
which implies that
\begin{align*}
\int_{[0,\fz)\times\mathbb{R}^2} \phi_{\theta_n} \partial_t\psi + \phi_{\theta_n} \Delta \psi + \frac 12\lf(\left[ \phi_{\theta_n} \otimes_{\theta_n} \phi_{\theta_n} \right]+\left[ \phi_{\theta_n} \otimes_{\theta_n} \phi_{\theta_n} \right]^T\r) \nabla \psi = -\int_{\mathbb{R}^2} \phi_0 \psi(0).
\end{align*}
Since \( \phi_{\theta_n} \to \phi \) in the weak-\( \ast \) topology, it follows that
\begin{align*}
\int_{[0,\fz)\times\mathbb{R}^2} \phi_{\theta_n} \partial_t\psi \to \int_{\mathbb{R}^2} \phi \partial_t\psi, \quad \int_{[0,\fz)\times\mathbb{R}^2} \phi_{\theta_n} \Delta \psi \to \int_{\mathbb{R}^2} \phi \Delta \psi.
\end{align*}
By Proposition \ref{star}, we also have
\[
\int_{[0,\fz)\times\mathbb{R}^2} \left[ \phi_{\theta_n} \otimes_{\theta_n} \phi_{\theta_n} \right] \cdot \nabla \psi \to \int_{[0,\fz)\times\mathbb{R}^2} (\phi \otimes \phi) \cdot \nabla \psi
\]
and
\[
\int_{[0,\fz)\times\mathbb{R}^2} \left[ \phi_{\theta_n} \otimes_{\theta_n} \phi_{\theta_n} \right]^T \cdot \nabla \psi \to \int_{[0,\fz)\times\mathbb{R}^2} (\phi \otimes \phi)^T \cdot \nabla \psi=\int_{[0,\fz)\times\mathbb{R}^2} (\phi \otimes \phi) \cdot \nabla \psi.
\]
Thus, we obtain
\begin{align*}
\int_{[0,\fz)\times\mathbb{R}^2} \phi \partial_t\psi + \phi \Delta \psi + (\phi \otimes \phi) \cdot \nabla \psi = \int_{\mathbb{R}^2} \phi_0\psi(0),
\end{align*}
which implies that \( \phi \) is a weak solution to the classical Navier-Stokes equation with initial datum \( \phi_0 \).
Furthermore, note that \( \{\phi_{\theta_n}\}_n \) is a bounded subset of \( L_2([0,\fz);[\dot{H}^1(\mathbb{R}^2)]_0^2) \) and therefore contains a subsequence convergent in the weak-$\ast$ topology of \( L_2([0,\fz);[\dot{H}^1(\mathbb{R}^2)]_0^2) \). This allows us to conclude that \( \phi \in L_\fz([0,\fz);[L_2(\mathbb{R}^2)]_0^2) \cap L_2([0,\fz);[\dot{H}^1(\mathbb{R}^2)]_0^2) \).
By the properties of weak-\( \ast \) convergence, we obtain the energy inequality:
\begin{align*}
&\frac{1}{2}\|\phi(t)\|^2_{[L_2(\mathbb{R}^2)]^2} + \int^t_0 \|\nabla\phi(s)\|^2_{[L_2(\mathbb{R}^2)]^4} ds \\
\le &\liminf_{n\to\infty} \left( \frac{1}{2}\|\phi_{\theta_n}(t)\|^2_{[L_2(\mathbb{R}^2)]^2} + \int^t_0 \|\nabla\phi_{\theta_n}(s)\|^2_{[L_2(\mathbb{R}^2)]^4} ds \right)\\
 =&\frac{1}{2}\|\phi_0\|^2_{[L_2(\mathbb{R}^2)]^2}
\end{align*}
for all \( t <\fz \), which implies that \( \phi \) is a Leray-Hopf weak solution, whose definition can be found in e.g. \cite[Theorem 5.1]{fjr72}.
Now we show that \( \phi_\theta \to \phi \) as \( \|\theta\| \to 0 \). For any subsequence of \( \{\phi_\theta\}_\theta \), there exists a convergent subsequence whose limit is \( \tilde{\phi} \). Using the same reasoning as above, we can show that \( \tilde{\phi} \) is also a Leray-Hopf weak solution with initial datum \( \phi_0\). It is well-known that Leray-Hopf weak solution is unique in two spatial dimensions, a proof of which can be found in \cite[Theorem 2.1]{whhg11}, so we conclude that \( \tilde\phi = {\phi}\). The Urysohn subsequence principle implies that \( \phi_\theta \to \phi \) in the weak-\( \ast \) topology of \( L_\infty([0,\fz);[L_2(\mathbb{R}^2)]^2_0) \) as \( \|\theta\| \to 0 \).

\end{proof}

\begin{remark}\label{reason}
{\rm    If we assume that $u_\theta$ in Theorem \ref{semiclassical} is the unique solution to \eqref{ens} with initial datum $u_{\theta,0}=U_\theta\circ\cF(\phi_0),$ via a similar procedure, we can also extract a weak-$\ast$ convergent subsequence $\{\phi_{\theta_n}\}_n$ such that its limit $\phi$ is a weak solution to the classical Navier-Stokes equation with initial datum $\phi_0.$ However, $\phi$ may not be a Leray-Hopf weak solution since $u_\theta$ may fail to satisfy the energy identity. Hence such $\phi$ may not be unique, which means that $\{\phi_\theta\}_\theta$ may not be convergent.}
\end{remark}

Next we give the proof of Theorem \ref{semi2}.

\begin{proof}[Proof of Theorem \ref{semi2}]
It suffices to consider the equation \eqref{ens1}, another one \eqref{ens} can be handled similarly.
We claim that there exists a constant $0<T_{\phi_0}\le\fz$ such that under the assumptions of $\phi_0,$ \(  \exists \,u_\theta \in \mathcal{N}_{T_{\phi_0}}^\theta \) which is the unique solution to \eqref{ens1} with initial datum \( u_{\theta,0}\). From Theorem \ref{t3.1}, one may let $T_{u_{\theta,0}}$ be the maximal existence time with respect to $u_{\theta}$. Then it suffices to show that $0<\inf_{\theta} T_{u_{\theta,0}}\le\fz,$ where $\theta$ is taken over all $d\times d$ antisymmetric matrices. So that we can take $T_{\phi_0}=\inf_{\theta} T_{u_{\theta,0}}.$

If $\phi_0$ satisfies the assumption of Theorem \ref{semi2} \rm{(i)}, then by Lemma \ref{l2.a1}, $u_{\theta,0}$ is uniformly bounded in $[L_{p}(\rd)]^d$ since $\|u_{\theta,0}\|_{[L_p(\rd)]^d}\lesssim \|\cF \phi_0\|_{[L_{p'}(\bR^d)]^d}$ for $p'=d'-\varepsilon$ and $2.$ Interpolating between $d'-\varepsilon$ and $2,$ we also have $\|u_{\theta,0}\|_{[L_d(\rd)]^d}\lesssim \|\cF \phi_0\|_{[L_{d'}(\bR^d)]^d}.$ Moreover, the above estimates are independent of $\theta.$ Then the existence of $T_{u_{\theta,0}}$ will be ensured by Theorem \ref{t3.1}. By the blow up criterion \eqref{blowup}, we have $T_{u_{\theta,0}}\gtrsim (\|\cF\phi_0\|_{[L_{p'}(\bR^d)]^d})^{-\frac{2p}{p-1}}$, which implies that $\inf_{\theta}T_{u_{\theta,0}}>0.$
If $\phi_0$ satisfies the assumption of Theorem \ref{semi2} \rm{(ii)}, then $\|u_{\theta,0}\|_{[L_d(\rd)]^d}\lesssim \|\cF\phi_0\|_{[L_{d'}(\rd)]^d}$ are sufficiently small, which implies the global well-posedness of $u_\theta$ by Theorem \ref{t3.1} \rm{(i)}. Hence $\inf_{\theta}T_{u_{\theta,0}}=\fz.$

Now the proof is similar to that of Theorem \ref{semiclassical}. Indeed, \( u_\theta \) is smooth by Theorem \ref{t3.1}, then using Lemma \ref{Wigner bounds} with \( p = d \) and \( p = d+2 \), we can deduce that \( \phi_\theta \) is a solution to \eqref{ens2} and that \( \{\phi_\theta\}_\theta \) forms a bounded subset of \( \mathcal{N}^0_{T} \) for a fixed $0<T<T_{\phi_0}$. Repeating the process used for \( d = 2 \), we can find a subsequence of \( \{\phi_\theta\}_\theta \) that converges in weak-\( \ast \)  topology to \( \phi \), where \( \phi \in L_{\infty}([0,T);[L_2(\mathbb{R}^d) \cap L_d(\mathbb{R}^d)]_0^d) \cap L_{d+2}([0,T);[L_{d+2}(\mathbb{R}^d)]^d_0) \) is a weak solution to the classical Navier-Stokes equation with initial datum \( \phi_0 \).
It was proven in \cite[Theorem 2.1,Theorem 3.3]{fjr72} that a weak solution in \( L_{d+2}([0,T);[L_{d+2}(\mathbb{R}^d)]^d_0) \), which may not be a Leray-Hopf weak solution, but  is still a mild solution, and thus unique. The remaining part of the proof is therefore analogous to the case \( d = 2 \).
\end{proof}

\noindent \textbf{Acknowledgements.} The authors are partially supported by National Natural Science Foundation of China (No. 12071355, No. 12325105, No. 12031004, and No. W2441002).
L. Wang is partially supported by a grant from the Research Grants Council of the Hong
Kong Administrative Region, China (No. CityU 21309222).
W. Wang is supported by
 China Postdoctoral Science Foundation (No. 2024M754159), and Postdoctoral Fellowship Program of CPSF (No. GZB20230961).

\bigskip

\end{document}